\documentstyle[12pt,newlfont]{amsart}

\newcommand{\nc}{\newcommand}

\nc{\one}{\mbox{\bf 1}}
\nc{\invtensor}{\underset{\leftarrow}{\otimes}}
\nc{\ad}{\operatorname{ad}}
\nc{\rk}{\operatorname{rank}}
\nc{\corank}{\operatorname{corank}}
\nc{\Sym}{\operatorname{Sym}}

\nc{\sym}{\operatorname{sym}}
\nc{\id}{\operatorname{id}}
\nc{\htt}{\operatorname{ht}}
\nc{\Ker}{\operatorname{Ker}}
\nc{\im}{\operatorname{Im}}
\nc{\re}{\operatorname{Re}}
\nc{\sn}{\operatorname{sn}}
\nc{\spn}{\operatorname{span}}

\nc{\Irr}{\operatorname{Irr}}
\nc{\sgn}{\operatorname{sgn}}
\nc{\F}{\operatorname{F}}

\nc{\Soc}{\operatorname{Soc}}
\nc{\Inj}{\operatorname{E}}
\nc{\Hom}{\operatorname{Hom}}
\nc{\End}{\operatorname{End}}
\nc{\supp}{\operatorname{supp}}
\nc{\Card}{\operatorname{Card}}
\nc{\Mod}{\operatorname{Mod}}
\nc{\Ann}{\operatorname{Ann}}
\nc{\Ind}{\operatorname{Ind}}
\nc{\Coind}{\operatorname{Coind}}

\nc{\wt}{\operatorname{wt}}
\nc{\ch}{\operatorname{ch}}
\nc{\Stab}{\operatorname{Stab}}
\nc{\Sch}{{\cal S}\mbox{\em ch}}

\nc{\Spec}{\operatorname{Spec}}
\nc{\Prim}{\operatorname{Prim}}
\nc{\Max}{\operatorname{Max}}
\nc{\Aut}{\operatorname{Aut}}
\nc{\Fract}{\operatorname{Fract}}
\nc{\gr}{\operatorname{gr}}
\nc{\Gr}{\operatorname{Gr}}
\nc{\wdM}{\widetilde{M}}
\nc{\wdV}{\widetilde{V}}
\nc{\wdN}{\widetilde{N}}
\nc{\wdchi}{\widetilde{\chi}}

\nc{\pprv}{{\cal J}}
\nc{\eprv}{{\det{\cal J}}}
\nc{\uu}{{\cal I}}

\nc{\cO}{\operatorname{\cal O}}
\nc{\wdO}{\operatorname{\widetilde{\cal O}}}
\nc{\Ob}{\operatorname{\cal Ob}}
\nc{\Dglie}{\operatorname{{\cal D}glie}}
\nc{\Fin}{\operatorname{{\cal F}in}}
\nc{\cC}{\operatorname{\cal C}}
\nc{\wdC}{\operatorname{\widetilde{\cal C}}}

\nc{\Sg}{{\cal S}({\frak g})}

\nc{\Ug}{\widetilde{\cal U}}
\nc{\Zg}{{\cal Z}({\frak g})}
\nc{\tZg}{{\widetilde{\cal Z}({\frak g})}}
\nc{\Zk}{{\cal Z}({\frak k})}
\nc{\Sh}{{\cal S}({\frak h})}
\nc{\Uh}{{\cal U}({\frak h})}

\nc{\Uk}{{\cal U}({\frak k})}
\nc{\Ag}{{\cal A}({\frak g})}

\nc{\cZ}{\cal Z}
\nc{\cS}{\cal S}
\nc{\cP}{\cal P}
\nc{\cL}{\cal L}
\nc{\cU}{\cal U}
\nc{\cH}{\cal H}
\nc{\cK}{\cal K}
\nc{\cF}{\cal F}
\nc{\fg}{\frak g}
\nc{\CO}{\cal O}
\nc{\fn}{\frak n}
\nc{\fm}{\frak m}
\nc{\fh}{\frak h}
\nc{\ft}{\frak t}
\nc{\fk}{\frak k}
\nc{\fp}{\frak p}
\nc{\fI}{\frak I}
\nc{\veps}{\varepsilon}
\nc{\fsl}{\frak {sl}}

\nc{\dirlim}{\underset{\rightarrow}{\lim}\,} 
\nc{\nen}{\newenvironment}
\nc{\ol}{\overline}
\nc{\ul}{\underline}
\nc{\ra}{\rightarrow}
\nc{\lra}{\longrightarrow}
\nc{\Lra}{\Longrightarrow}
\nc{\Lla}{\Longleftarrow}
\nc{\Llra}{\Longleftrightarrow}
\nc{\thla}{\twoheadleftarrow}
\nc{\hra}{\hookrightarrow}
\nc{\iso}{\overset{\sim}{\lra}}
\nc{\ssubset}{\underset{\not=}{\subset}}

\nc{\Thm}[1]{Theorem~\ref{#1}}
\nc{\Prop}[1]{Proposition~\ref{#1}}
\nc{\Lem}[1]{Lemma~\ref{#1}}
\nc{\Cor}[1]{Corollary~\ref{#1}}
\nc{\Conj}[1]{Conjecture~\ref{#1}}
\nc{\Claim}[1]{Claim~\ref{#1}}
\nc{\Defn}[1]{Definition~\ref{#1}}
\nc{\Exa}[1]{Example~\ref{#1}}
\nc{\Rem}[1]{Remark~\ref{#1}}
\nc{\Note}[1]{Note~\ref{#1}}
\nc{\Quest}[1]{Question~\ref{#1}}
\nc{\Hyp}[1]{Hypoth\`ese~\ref{#1}}

\nen{thm}[1]{\label{#1}{\bf Theorem.\ } \em}{}
\nen{prop}[1]{\label{#1}{\bf Proposition.\ } \em}{}
\nen{lem}[1]{\label{#1}{\bf Lemma.\ } \em}{}
\nen{cor}[1]{\label{#1}{\bf Corollary.\ } \em}{}
\nen{conj}[1]{\label{#1}{\bf Conjecture.\ } \em}{}
\nen{claim}[1]{\label{#1}{\bf Claim.\ } \em}{}

\nen{defn}[1]{\label{#1}{\bf Definition.\ } }{}
\nen{exa}[1]{\label{#1}{\bf Example.\ } }{}

\nen{rem}[1]{\label{#1}{\em Remark.\ } }{}
\nen{note}[1]{\label{#1}{\em Note.\ } }{}
\nen{exer}[1]{\label{#1}{\em Exercise.\ } }{}
\nen{sket}[1]{\label{#1}{\em Sketch of proof.\ } }{}
\nen{quest}[1]{\label{#1}{\bf Question.\ } \em}{}
\nen{hyp}[1]{\label{#1}{\bf Hypoth\`ese.\ } \em}{}
\begin{document}


\title[]{Annihilation Theorem and Separation Theorem for  basic classical
Lie superalgebras}

\author[]{Maria Gorelik}

\address{ 
{\tt email: gorelik@@mpim-bonn.mpg.de} 
}

\thanks{The author was partially supported by TMR Grant No. FMRX-CT97-0100.
Research at MSRI is supported in part by NSF grant DMS-9701755.}

\begin{abstract}
In this article we prove that for a basic classical Lie superalgebra
the annihilator of a strongly typical
Verma module is a centrally generated ideal. 
For a basic classical Lie superalgebra of type I we prove that
the localization of the enveloping algebra by a certain central element
is free over its centre. 
\end{abstract}

\maketitle

\section{Introduction}
\subsection{}
Let $\fg$ be a complex semisimple Lie algebra, $\cU$ be its
universal enveloping algebra and $\Zg$ be the centre of $\cU$.
Consider $\cU$
as a $\fg$-module with respect to the adjoint action. 
Separation Theorem of Kostant (see~\cite{ko}) 
states the existence of a
submodule $H$ of $\cU$ such that the multiplication map
provides the bijection $H\otimes \Zg\iso\cU$.
Moreover the multiplicity of each simple finite dimensional module
$V$ in $H$ is equal to the dimension of its zero weight space.
Such an $\ad\fg$-invariant subspace $H$ is called a {\em harmonic space}.
An easy proof of Separation Theorem was found by Bernstein and
Lunts---see~\cite{bl}.
This theorem is an important ingredient in the proof of
the annihilation theorem of Duflo (see~\cite{d}, 8.4.3)
asserting that the annihilator of a Verma module is generated
by its intersection with $\Zg$.
The annihilation theorem is reproven by Joseph and G.~Letzter.
They also generalize it to the quantum case--- see~\cite{jl},\cite{jnato}.

In this paper we obtain analogous theorems in the case 
of basic classical Lie superalgebras. This was done earlier 
for the completely reducible case--- see~\cite{mu},\cite{gl1}.

Let $\fg=\fg_0\oplus\fg_1$ be a basic classical Lie superalgebra, $\Ug$ be its
universal enveloping superalgebra and $\Zg$ be the centre of $\Ug$.
Let $T$ be a special {\em ghost} element constructed in~\cite{ghost}. 
Call a highest weight module {\em strongly typical}
if it is not annihilated by $T$.

We prove the following version of Annihilation Theorem.
\subsubsection{}
\begin{thm}{intr:AH}
The annihilator of a strongly typical
Verma $\fg$-module $\wdM$ is a centrally generated ideal.
\end{thm}

Moreover, for a strongly typical Verma module $\wdM$,
we describe the quotient $\Ug/\Ann\wdM$ as
an $\ad\fg$-module and show that the natural map
from $\Ug/\Ann\wdM$ to the locally finite part
of $\End_{\Bbb C}(\wdM)$ is bijective.

\subsubsection{}
The proof of~\Thm{intr:AH} goes as follows. As in~\cite{jl}, 
we use the Parthasarathy--Ranga-Rao--
Varadarajan (PRV) determinants. We generalize the notion of 
PRV determinants (see~\cite{prv}) to basic classical Lie superalgebras. 
Our construction is based on the fact that the two-sided ideal
$\Ug T=T\Ug$, considered as  $\ad\fg$-module, 
is injective in a category of locally finite modules.

In the completely reducible case one assigns
to each simple finite dimensional
module $\wdV$  a PRV determinant which is
a polynomial in $\Sh$.

On the contrary, for  the 
non-completely reducible case the lack
of harmonic space forces us to substitute 
the PRV determinant by a set of PRV determinants  
corresponding to the same $\wdV$. However, we do not 
have to calculate these determinants, but only verify that they are 
non-zero. We use these determinants to show
that if the locally finite part $F(\wdM,\wdM)$ of
the endomorphisms $\End_{\Bbb C}(\wdM)$ of
a strongly typical Verma module $\wdM$ has
``a right size''
as an $\ad\fg$-module and the natural
map $\Ug\to F(\wdM,\wdM)$ is surjective, then $\Ann \wdM$
is centrally generated. 

For type I we directly verify both conditions. 
The crucial point in the study of type II case is the 
construction in Section~\ref{apprpairs} of a {\em perfect mate} 
$\chi\in\Max\cZ(\fg_0)$ for each strongly typical $\wdchi\in\Max\Zg$.
We call a maximal ideal $\chi\in\Max\cZ(\fg_0)$
a perfect mate for $\wdchi\in\Max\Zg$
if the following conditions are satisfied.

(i) For any Verma $\fg$-module
annihilated by $\wdchi$, its $\fg_0$-submodule annihilated
by a power of $\chi$ is a Verma $\fg_0$-module. 

(ii) Any $\fg$-module annihilated by $\wdchi$ has
a non-trivial $\fg_0$-submodule annihilated by $\chi$. 

The condition (ii) seems to be difficult to check.
However, it turns out that it is enough to verify (ii)
only for simple highest weight $\fg$-modules.
This is deduced from~\cite{mu2}.

For $\fg={\frak {osp}}(1,2l)$ 
the annihilator
of a  Verma $\fg$-module $\wdM$ is a centrally generated ideal
iff $\wdM$ is strongly typical---see~\cite{gl1}.
In this paper we prove the similar equivalence for the
basic classical Lie superalgebras of  type I.

As it is shown in~\cite{psI}, if $\fg$ has type I
then for any strongly typical $\wdchi\in\Zg$
the algebra $\Ug/(\Ug\wdchi)$ is 
the matrix algebra over $\cU(\fg_0)/(\cU(\fg_0)\chi)$ for a suitable
$\chi\in\Max\cZ(\fg_0)$. As it was pointed out by V.~Serganova
this result implies~\Thm{intr:AH} for type I case.
The opposite implication is also easy (see~\ref{EquivTypeI}).

\subsection{}
Let $\fg$ be a basic classical Lie superalgebra which is not
completely reducible. Then $\Ug$ is not a domain (see~\cite{al}) and  $\Zg$
is not Noetherian (\cite{mu}, 2.8).
However all non-zero central elements are non-zero divisors
and $\Zg$ contains an element 
$z$ such that the localized algebra
$\Zg[z^{-1}]$ is isomorphic to a localization of a polynomial algebra.
One can take $z:=T^2$ where $T$ is the element mentioned above.
The element $T$ is even; it commutes with the even elements of $\Ug$ and 
anticommutes with the odd ones. Moreover, the image of $T$ 
in the symmetric algebra $\Sg$ belongs to the
top exterior  power of $\fg_1$. These properties 
determine $T$ up to a scalar. 

It is easy to show (see~\cite{ghost}, 4.5) that $\Ug$, considered
as an $\ad\fg$-module, does not admit a factorization of the form 
$H\otimes\Zg$. For type I Lie superalgebras,
we prove the following version of Separation Theorem.
\subsubsection{}
\begin{thm}{intr:SP}
For $\fg$  of the type I,  there exists 
$\ad\fg$-invariant subspace $H$ of $\Ug$
such that the multiplication map provides the 
bijection $H\otimes \Zg[T^{-2}]\to\Ug[T^{-2}]$.
\end{thm}

Clearly, $\Zg[T^{-2}]$ coincides with the centre of $\Ug[T^{-2}]$.
As an $\ad\fg$-module, $H$ is injective in an appropriate category
of locally finite modules and for any simple finite dimensional module $\wdV$
one has $\dim \Hom_{\fg}(\wdV,H)=\dim\wdV|_0$ where $\wdV|_0$
is the zero weight space of $\wdV$.
For a basic classical Lie superalgebra of type II, we obtain
a weaker result, namely that
the similar assertions hold if we substitute $z$ by
 a certain subset $S$ of $\Zg$. This set $S$ can be
described in terms of the PRV determinants. 

A natural conjecture is
that one can always choose $S$ equal to $\{T^2\}$.
A possible way to prove this conjecture is to show 
that an irreducible
factor of a PRV determinant is either a factor of Shapovalov 
form or is of the form $(\beta^{\vee}+\beta^{\vee}(\rho))$
for some odd coroot $\beta$. However it is not clear
how to calculate these determinants if $\fg$ is not
completely reducible.

\subsection{Content of the paper.}
In Section~\ref{sectprelim} we recall some facts about the 
basic classical Lie superalgebras. 

In Section~\ref{sectFin} we define
a category $\Fin$ of locally finite $\fg$-modules and provide
some properties of $\Fin$. We also recall the construction
and properties of the element $T$.

In Section~\ref{sectPRV} we investigate the $\fg$-module
structure of $\Ug$ given by the adjoint action. 
We start with studying $\Hom_{\frak g}(\wdV,{\Ug})$ 
for a simple finite dimensional module $\wdV$.
For each $\wdV$ we construct
a central element $z$ such that the localized module
$\Hom_{\frak g}(\wdV,{\Ug}[z^{-1}])=\Hom_{\frak g}(\wdV,{\Ug})[z^{-1}]$
is a free ${\Zg}[z^{-1}]$-module whose rank is equal 
to $\dim\wdV|_0$. In~\ref{sPRV} 
we generalize a notion of PRV determinants to
the case of non-completely reducible Lie superalgebras. 
We also establish properties of these determinants
which are similar to the properties of the original
PRV determinants. In~\ref{ssectsprthm} we show that
for a suitable subset $S$ of $\Zg$ the localized algebra
$\Ug[S^{-1}]$ is free over its centre $\Zg[S^{-1}]$.
We describe the $\ad\fg$-module structure
of the corresponding ``generic harmonic space'' $H$. We show
that one may choose $H$ to be the $\ad\fg$-module
generated by $H'T$ where $H'$ is a certain harmonic space of
$\cU(\fg_0)$.

In Section~\ref{sectSA} we establish a connection
between PRV determinants and the annihilators of simple modules.
We show that if all PRV determinants do not vanish at a point
$\lambda\in\fh^*$ and if a  simple module $\wdV(\lambda)$ 
is strongly typical then its annihilator is a
centrally generated ideal. Moreover, for
a  simple strongly typical Verma module $\wdV(\lambda)$ 
all PRV determinants do not vanish at a point
$\lambda\in\fh^*$ iff the natural map
$\Ug\to F(\wdM,\wdM)$ is surjective and $F(\wdM,\wdM)$
has  a certain nice structure as $\ad\fg$-module
(it should be isomorphic to the ``generic harmonic space'' $H$).

Sections~\ref{prepI},\ref{typeI} are devoted to type I
case.  Preliminary facts are concentrated in Section~\ref{prepI}.
In Section~\ref{typeI} we prove 
that for a suitable $\ad\fg$-stable $H$ the multiplication map
provides an isomorphism $H\otimes\Zg [T^{-2}]\to \Ug [T^{-2}]$.
We also prove that the annihilator of a Verma module
is centrally generated iff this module is strongly typical.

Sections~\ref{apprpairs},\ref{sectannthm}  are devoted to type II
case. 

In Section~\ref{apprpairs} we describe, for
each strongly typical $\wdchi\in\Max\Zg$ its perfect mate
$\chi\in\Max\cZ(\fg_0)$. 
Note that in   type I case for any  strongly typical $\wdchi\in\Max\Zg$ 
and for any $\lambda$ such that $\wdchi$ 
annihilates $\wdM(\lambda)$, the ideal 
$\Ann_{\cZ(\fg_0)} M(\lambda)$ is a perfect mate for $\wdchi$.
This does not hold for type II. For
certain ``generic'' $\wdchi$,   the ideal 
$\Ann_{\cZ(\fg_0)} M(\lambda)$ is a perfect mate if
one chooses $\lambda$ satisfying, apart from $\wdchi\wdM(\lambda)=0$,
also a kind of ``dominance''
condition. For $B(m,n)$ and $G(3)$ all  strongly typical central characters
are generic. For the remaining superalgebras
$D(m,n), D(2,1,\alpha)$ and $F(4)$
we select perfect mates for  non-generic strongly typical
central characters case by case.

In Section~\ref{sectannthm} we prove~\Thm{intr:AH}. The existence of
a perfect mate for each strongly typical 
$\wdchi\in\Max\Zg$ enables us to show 
that for a Verma module $\wdM$ with the central character $\wdchi$
the natural map $\Ug\to F(\wdM,\wdM)$ is surjective and $F(\wdM,\wdM)\cong H$. 
According Section~\ref{sectSA}, these two conditions
imply that the annihilator of $\wdM$ is centrally generated ideal.

In Section~\ref{sectRmV}
we study  the $\fg_0$-structure of Verma $\fg$-modules.

Appendix~\ref{appendix} contains some lemmas used in the main text. 
It also contains alternative proofs of some results 
(non-vanishing of the PRV determinant,\Cor{cordetnn} and~\Thm{fact?}).
These  proofs do not use Kostant's Separation Theorem.
Probably, these proofs may be 
useful in the case when   separation theorem does not hold.

{\em Acknowledgments.} It is a pleasant duty to express the gratitude
to my teacher A.~Joseph. 
Many ideas of this paper owe to his book and to his courses given at 
Weizmann Institute.
The author would like to thank M.~Duflo, E.~Lanzmann,
I.~Musson and I.~Penkov for helpful
discussions and V.~Serganova for her extreme patience
and useful comments. This work was done when the author
was a postdoctoral fellow at MSRI.  The author would like to thank
MSRI and the organizers of the Noncommutative Algebra Program
for their support and hospitality.

\section{Preliminaries}
\label{sectprelim}
In this paper the ground field is ${\Bbb C}$. 
Everywhere in the paper apart from
Sections~\ref{prepI},\ref{typeI},
${\frak g}=\fg_0\oplus\fg_1$ denotes one (unless otherwise
specified, an arbitrary one) of
the basic classical complex Lie superalgebras
${\frak {gl}}(m,n), {\frak {sl}}(m,n),
{\frak {osp}}(m,n), {\frak {psl}}(n,n)$.
Each of these Lie superalgebras possesses the following properties:
it admits a ${\frak g}$-invariant
bilinear form which is non-degenerate on $[\fg,\fg]$
and the even part
${\frak g}_0$ is a reductive Lie algebra.

In Sections~\ref{prepI},\ref{typeI}, ${\frak g}=\fg_{\ol 0}\oplus\fg_{\ol 1}$ 
denotes one of the basic classical complex Lie superalgebras
of type I: ${\frak {gl}}(m,n), 
{\frak {sl}}(m,n),
{\frak {osp}}(2,n), {\frak {psl}}(n,n)$. We shall slightly change our
notations since these superalgebras admit ${\Bbb Z}$-grading
${\frak g}=\fg_{-1}\oplus\fg_0\oplus\fg_1$. The even part
$\fg_{\ol 0}$ coincides with $\fg_0$ and the odd part $\fg_{\ol 1}$
is the sum of two dual $\fg_0$-modules $\fg_{-1}$ and $\fg_1$.

In this section we present the main preliminary facts 
about the structure
of the basic classical complex Lie superalgebras
and their representations, which we shall use in this paper.

\subsection{Conventions}
We denote by ${\Bbb N}^+$ the set of positive integers
and by $\#I$ the number of elements of the given set $I$.
If $A$ is an algebra, $N$ is an $A$-module and
$X,Y$ are subsets of $A$ and $N$ respectively, 
we denote by $XY$ the set of the products $xy$ where $x\in X,y\in Y$.

For a ${\Bbb Z}_2$-homogeneous
element $u$ of a superalgebra denote by $d(u)$ its 
${\Bbb Z}_2$-degree. In all formulas where this notation is used,
$u$ is assumed to be ${\Bbb Z}_2$-homogeneous.

For a superalgebra $\fm$ denote by $\cU(\fm)$ its universal enveloping
superalgebra. Set $\Ug:=\cU(\fg)$ and $\cU:=\cU(\fg_0)$.

In this text all modules are assumed to be left modules.
For given module $N$ we denote by $N^{\oplus r}$ the direct sum of
$n$-copies of $N$. We say that $A$-module $N$ is locally
finite if $\dim Av<\infty$ for all $v\in N$.

The symbol $\wdV$ (resp., $V$) is always used for
a simple $\fg$ (resp., $\fg_0$) module
and the symbol $\wdM$ (resp., $M$) for
a Verma $\fg$ (resp., $\fg_0$) module.

\subsection{Triangular decompositions}
\label{trdec}
Triangular decompositions 
of the superalgebras are defined in~\cite{psg} as follows.
A Lie subsuperalgebra $\fh\subset\fg$ is called a {\em Cartan
subsuperalgebra} if $\fh$ is nilpotent and coincides 
with its centralizer in $\fg$. 
For the basic classical Lie superalgebras the set of
Cartan subalgebras coincides with the set of Cartan subalgebras
of $\fg_0$. Fix a Cartan subalgebra $\fh$;
it acts semisimply on ${\frak g}$:
$${\frak g}:=\displaystyle\oplus_{\mu\in {\frak h}^*} {\frak g}|_{\mu},
\ \ \ {\frak g}|_{\mu}:=\{a\in\fg| \ \forall h\in\fh,\ [h,a]=\mu(a)\}.$$
Denote by $\Delta$ the set of non-zero roots that is the set
$\{\alpha\in{\frak h}^*|\ 
{\frak g}|_{\alpha}\not=0\}\setminus\{0\}$.
An element $h\in\fh$ is called regular if $\re \alpha(h)\not=0$
for all $\alpha\in \Delta$. Any regular element determines
the decomposition $\Delta=\Delta^+\coprod\Delta^-$
where
$$\begin{array}{cc}
\Delta^+:=\{\alpha\in\Delta|\ \re\alpha(h)>0\}, &
\Delta^-:=\{\alpha\in\Delta|\ \re\alpha(h)<0\}.\ 
\end{array}$$
Moreover it determines a decomposition 
$\fg={\fn}^-\oplus\fh\oplus{\fn}^+$
where
$$\begin{array}{cc}
\fn^+:=\displaystyle\oplus_{\alpha\in\Delta^+} {\frak g}|_{\alpha}, &
\fn^-:=\displaystyle\oplus_{\alpha\in\Delta^-} {\frak g}|_{\alpha}.\ 
\end{array}$$
Such decompositions of $\fg$ are called 
{\em triangular decompositions}.
It is clear that $\fn^{\pm}$ are nilpotent Lie subsuperalgebras
of $\fg$.

A Lie subsuperalgebra ${\frak b}\subset\fg$ is called 
a {\em Borel subsuperalgebra} if ${\frak b}=\fh\oplus{\fn}^+$
for some triangular decomposition 
${\frak g}={\fn}^-\oplus\fh\oplus{\fn}^+$. A Borel subsuperalgebra
determines the triangular decomposition; we will add a lower
index to designate the corresponding Borel subsuperalgebra
in the case where the choice of triangular decomposition
is not clear from the context. We
denote by $\Delta({\frak b})$ the set of non-zero roots of ${\frak b}$.
We say that a vector of $v$ of $\fg$-module is ${\frak b}$-primitive
if $[{\frak b},{\frak b}]v=0$ and $\fh v\in {\Bbb C}v$.

A triangular decomposition 
${\frak g}={\fn}^-\oplus\fh\oplus{\fn}^+$
induces the triangular decomposition of the even part
${\frak g}_0={\fn}_0^-\oplus\fh\oplus{\fn}_0^+$.
The group of inner automorphisms of $\fg_0$ acts transitively
on all triangular decompositions of $\fg_0$ and the action
of this group can be extended to $\fg$. 
Hence the theory does not depend from the choice of
a triangular decomposition of $\fg_0$. 
In the sequel we fix a triangular decomposition 
$\fg_0={\fn}_0^-\oplus\fh\oplus{\fn}_0^+$ 
and consider triangular decompositions of $\fg$
which induce this fixed triangular decomposition of $\fg_0$.

\subsubsection{}
Denote by $\Delta_0$ the set of non-zero even 
roots of ${\frak g}$ and by $\Delta_1$ the set of odd 
roots of ${\frak g}$. Set $\Delta^{\pm}_0:=\Delta_0\cap \Delta^{\pm}$
and $\Delta^{\pm}_1:=\Delta_1\cap \Delta^{\pm}$.

Denote by $(-,-)$ a 
$\fg$-invariant bilinear form on $\fg$ which is non-degenerate
on $[\fg,\fg]$ and the induced
$W$-invariant bilinear form on ${\frak h}^*$.
A root $\alpha\in\Delta$ is called {\em isotropic}
if $(\alpha,\alpha)=0$. For a root $\alpha$ 
denote by $\alpha^{\vee}$ the element of $\fh$
satisfying $\alpha^{\vee}(\mu)=(\alpha,\mu)$ for each $\mu\in\fh^*$.

Set
$$\overline{\Delta}^+_{0}\!\!:=
\{\alpha\in\Delta^+_{0}|\ \alpha/2\not\in\Delta^+_{1}\},\ \ \ 
\overline{\Delta}^+_{1}\!\!:=\{\beta\in\Delta^+_{1}|\ 
2\beta\not\in\Delta^+_{0}\}.$$
The set of isotropic roots coincides with
 $\overline{\Delta}^+_{1}$.

Remark that  $\sum_{\alpha\in\Delta^+} n_{\alpha}\alpha=0$
for some $n_{\alpha}\in {\Bbb N}$ implies $n_{\alpha}=0$
for all $\alpha\in\Delta^+$. This allows us to define
the standard partial order relation on ${\frak h}^*$ by
$\lambda\leq \mu \Longleftrightarrow \mu-\lambda\in 
\sum_{\alpha\in\Delta^+} {\Bbb N}\alpha$. One can easily
sees that the minimal (with respect
to this partial order) elements of $\Delta^+$ form
a basis of simple roots.

Denote by $\pi_0$ the basis of simple roots of ${\frak g}_0$
and by $W$ the Weyl group of $\Delta_0$. Denote
by $|W|$ the number of elements in $W$.
For $w\in W$ set $\sn (w):=(-1)^{l(w)}$ where $l(w)$ is the length of 
$w$. For a non-isotropic root $\alpha$ define 
$s_{\alpha}\in\Aut {\frak h}^*$ by setting 
$$s_{\alpha}(\lambda):=\lambda-
2{(\alpha,\lambda)\over(\alpha,\alpha)}\alpha.$$
Evidently $s_{k\alpha}=s_{\alpha}$ and so the subgroup
of $\Aut {\frak h}^*$ generated by the $s_{\alpha}$
coincides with $W$.

Set 
$$\begin{array}{ccc}
\rho_0\!:={1\over 2}\displaystyle\mathop{\sum}_
{\alpha\in \Delta_0^+}\alpha,\ & 
\rho_1\!:={1\over 2}\displaystyle\mathop{\sum}_
{\alpha\in \Delta_1^+}\alpha,\ & 
\rho\!:=\rho_0-\rho_1.
\end{array}$$ 
For a simple root $\alpha$ one has
$2(\alpha,\rho)=(\alpha,\alpha)$.

Define the translated action of $W$ on ${\frak h}^*$ and on 
the symmetric algebra $\Sh$
by the formulas: 
$$w.\lambda:=w(\lambda+\rho)-\rho,\ \ w.f(\lambda):=f(w^{-1}.\lambda),\ \
\forall \lambda\in {\frak h}^*,w\in W.$$

\subsubsection{}
For a $\fg_0$-module $N$ and an element $\mu\in\fh^*$ set
$$N|_{\mu}:=\{m\in M|\ hm=\mu(h)m,\ \forall h\in\fh\}.$$
We shall consider mainly $\fh$-diagonalizable modules that is
satisfying $N=\sum_{\mu\in \fh^*} N|_{\mu}$. Set
$$\Omega(N):=\{\mu\in {\frak h}^*|\ N|_{\mu}\not=0\}.$$
If $\dim N|_{\mu}<\infty$ for each $\mu\in {\fh}^*$ ,
set
$\ch N:=\sum_{\lambda\in {\frak h}^*} (\dim N|_{\mu})e^{\mu}$.

When we use the notation $\Ug|_{\mu}$, the action of $\fg$ on
$\Ug$ is assumed to be the adjoint action.

\subsection{}
\label{cntrgr}
An important property of the basic classical Lie superalgebras
is the existence of a Cartan superantiautomorphism $\sigma$
coming from the supertransposition of matrices. 
Recall that an even linear endomorphism $\iota$ of a Lie 
(resp., associative) superalgebra
is called a {\em  superantiautomorphism} if
$\iota([x,y])=(-1)^{d(x)d(y)}[\iota(y),\iota(x)]$ 
(resp., $\iota(xy)=(-1)^{d(x)d(y)}\iota(y)\iota(x)$) for all 
homogeneous elements $x,y$. The Cartan superantiautomorphism $\sigma$
has the following properties: 
$$\begin{array}{rll}
a) & \sigma^2(g)=(-1)^{d(g)}g,& \forall g\in\fg,\\
b) & \sigma({\frak n}^+)={\frak n}^-,  & \\
c) &  \sigma(h)=h, & \forall h\in \fh. 
\end{array}$$

The restriction of $\sigma$ to $\fg_0$ is a Cartan antiinvolution:
$$\begin{array}{rll}
a) & (\sigma^2)|_{\fg_0}=\id, &\\
b) & \sigma({\frak n}_0^+)={\frak n}_0^-,  & \\
c) &  \sigma(h)=h, & \forall h\in \fh. 
\end{array}$$

\subsubsection{Symmetric algebra.}
Denote by $\cF$ the canonical filtration of $\Ug$ given by
$\cF^k:=\fg^k$. This filtration is $\ad\fg$-invariant
and the associated graded superalgebra
$\Sg$ inherits an  $\ad\fg$-module structure.
The superalgebra $\Sg$ is supercommutative: it is the product
of the symmetric (even) algebra
$\cS(\fg_0)$ and the external superalgebra
$\Lambda\fg_1$.

For $u\in\Ug$ denote by $\gr u$ its image in $\Sg$;
identify $\Uh$ and its image $\Sh$.

\subsubsection{Centre}
\label{syminv}
By definition, the (super)centre $\Zg:=\Ug^{\fg}$.
One has $\gr\Zg=\Sg^{\fg}$. 

Denote by $\gr P$ the projection
$\Sg\to\Sh$ along $\Sg\gr({\frak n}^-+{\frak n}^+)$.
The restriction of $\gr P$ to $\Sg^{\fg}$
provides a monomorphism $\iota:\Sg^{\fg}\to\Sh^{W}$.
As a consequence, all non-zero elements of $\Sg^{\fg}$ (resp., $\Zg$)
are non-zero divisors in $\Sg$ (resp., $\Ug$).
The image of $\iota$ is described 
in~\cite{kch1},~\cite{s1},~\cite{bzv}; $\iota$
is bijective iff $\fg={\frak {osp}}(1,2l)$.

Denote by $P_{\emptyset}$ the projection
$\Sg\to\cS(\fg_0)$ with the kernel
$K:=\sum_{i\geq 1}\cS(\fg_0)\Lambda^i\fg_1$.
It is easy to see that $K$ is $\ad\fg_0$-invariant
and so $P_{\emptyset}$ is an $\ad\fg_0$-map.
Moreover $P_{\emptyset}$ provides a monomorphism
$\Sg^{\fg}\to\cS(\fg_0)^{\fg_0}$ since 
$\gr P=\gr P\circ P_{\emptyset}$ and so
the injectivity of $P_{\emptyset}$ on $\Sg^{\fg}$
follows from the injectivity of $\gr P$.

\subsubsection{Harish-Chandra projection.}
Denote by $\cP$ the Harish-Chandra projection $\Ug\to\Sh$
with respect to the decomposition 
$\Ug=(\Ug{\frak n}^++{\frak n}^-\Ug)\oplus \cU(\fh)$.
The restriction of $\cP$ to $\Ug|_0=\Ug^{\fh}$ 
is an algebra homomorphism.
An element $a\in\Ug^{\fh}$ acts on a primitive vector
of weight $\mu$  by the multiplication
by the scalar $\cP(a)(\mu)$. 

The Harish-Chandra projection 
provides a monomorphism $\Zg\to\Sh^{W.}$.
If $\wdN$ is a $\fg$-module generated by a primitive
vector of weight $\lambda$, then a central
element $z$ acts on $\wdN$ by the multiplication
by the scalar $\cP(z)(\lambda)$.

Call $\wdchi\in\Max\Zg$ {\em a central character }
of a $\fg$-module $N$ if $\wdchi^r N=0$ for $r>>0$.

\subsubsection{}
\label{sigma;cent}
Let $z$ be an element of $\Zg$. Since $z$ has weight zero,
$z=\cP(z)+\sum_i u^{-}_iu^+_i$ where  $u^{-}_i\in \cU({\fh^-})\fn^-$,
$u^{+}_i\in \cU({\fh^+})\fn^+$ for all $i$.
One has $\sigma(u^{-}_iu^+_i)=\pm\sigma(u^+_i)\sigma(u^-_i)\in 
\cU({\fh^-})\fn^-\cU({\fh^+})\fn^+$. Therefore 
$$\sigma(z)=\sigma(\cP(z))+\sum_i \pm\sigma(u^+_i)\sigma(u^{-}_i)\in
\cP(z)+\Ug\fn^+$$
that is $\cP(\sigma(z))=\cP(z)$. Thus
the superantiautomorphism $\sigma$ 
stabilizes the central elements. This implies that
$\cZ(\fg_0),\Zg\subset\Ug^{\sigma}$.

\subsubsection{}
\label{defindcoind}
For a ${\Bbb Z}_2$-graded $\fg_0$-module $L$ 
denote by $\Ind_{\fg_0}^{\fg}L$ a vector space
$\Ug\otimes_{\cU} L$ (here $\Ug$ is considered as a right
$\cU$-module and a left $\Ug$-module through the multiplication)
equipped by the natural structure of a left 
$\Ug$-module.  
Denote by $\Coind_{\fg_0}^{\fg}L$ a vector space
$\Hom_{\cU}(\Ug,L)$ (here $\Ug$ is considered as a left
$\cU$-module) equipped by the following structure of a left 
$\Ug$-module: $(uf)(u'):=f(u'u)$ for any 
$f\in \Hom_{\cU}(\Ug,L),\ 
u,u'\in \Ug$. For a $\fg$-module $\wdN$ and a $\fg_0$-module $L$ one
has the canonical bijections
$$\begin{array}{cl}
\Hom_{\fg_0}(\wdN,L)\iso
\Hom_{\fg}(\wdN,\Coind_{\fg_0}^{\fg}L),& \\
\Hom_{\fg_0}(L, \wdN)\iso
\Hom_{\fg}(\Ind_{\fg_0}^{\fg}L, \wdN). &
\end{array}$$
By~\cite{bf}, 
$\Ind_{\fg_0}^{\fg}L\cong\Coind_{\fg_0}^{\fg}L$.

\subsection{Hopf algebra structure} 
The enveloping algebra $\Ug$ is a supercommutative 
Hopf superalgebra. This means, in particular, that 
the antipode $S$ is a superantiautomorphism of $\Ug$ and that
the comultiplication
${{\Delta}'}:\Ug\to\Ug\otimes\Ug$ is a homomorphism of superalgebras
satisfying the relation $s\circ \Delta'=\Delta'$
where $s$ is a linear map
$s:\Ug\otimes\Ug\to\Ug\otimes\Ug$
given by $s(u_1\otimes u_2):=(-1)^{d(u_1)d(u_2)}(u_2\otimes u_1)$.
The Hopf algebra structure on $\Ug$ is given by 
$$\begin{array}{rl}
\Delta'(g)&=g\otimes 1+1\otimes g,\\
\varepsilon(g)&=0,\\
S(g)&=-g\end{array}$$
for any $g\in {\frak g}$.

The Hopf superalgebra structure on $\Ug$ gives 
a ${\frak g}$-module structure 
on the tensor product $N_1\otimes N_2:=N_1\otimes_{\Bbb C} N_2$
of two ${\frak g}$-modules $N_1,N_2$. The map
$n_1\otimes n_2\mapsto (-1)^{d(n_1)d(n_2)}$
provides the canonical $\fg$-isomorphism $N_1\otimes N_2\to
N_2\otimes N_1$.

\subsubsection{}
Throughout the paper we shall write ``$\ad \fg$-module''
instead ``$\fg$-module with respect to the adjoint action''.

View $\Ug$ as $\fg$-module through the adjoint action given by
$$(\ad g)u:=ug-(-1)^{d(g)d(u)}ug,
\ \ \forall g\in {\fg},u\in{\Ug}.$$
As  $\ad\fg$-module, $\Ug$ is locally finite.

For any $\fg$-modules $N_1,N_2$ view 
$\Hom(N_1,N_2):=\Hom_{\Bbb C}(N_1,N_2)$
as a $\fg$-module through the adjoint action:
$$(\ad g)\psi(v)=g\psi(v)-(-1)^{d(g)d(\psi)}
\psi(gv) \ 
\ \forall g\in {\fg},\psi\in\Hom(N_1,N_2).$$
We denote by $F(N_1,N_2)$ the locally finite part of  $\ad\fg$-module
$\Hom(N_1,N_2)$. Notice that $F(N_1,N_2)$ coincides with 
the  $\ad\fg_0$-locally finite part of $\Hom(N_1,N_2)$, 
since $\Ug$ is a finite extension of $\cU$.

Throughout the paper an action of $\fg$ on ${\Ug}$ 
and on $F(N_1,N_2)$ is assumed,
by default, to be the adjoint action. 

For a $\fg$-module $N$, the natural map 
$\Ug\to \End(N)$, coming from the action of $\Ug$ on $N$,
is an $\ad\fg$-homomorphism and its image lies in $F(N,N)$.

\subsubsection{}
\label{Frrec}
Let $L$ be a finite dimensional  $\fg$-module. Equip
the dual supervector space $L^*$ by $\fg$-module structure
through the antipode $S$:
$$g.f(v):=(-1)^{d(g) d(f)}
f(S(g)v)=(-1)^{d(g) d(f)}f(-gv),\ \forall v\in L, g\in {\frak g}.$$

We shall use the following form of the Frobenius reciprocity
$$\Hom_{\fg}\bigl(L,\Hom(N_1,N_2)\bigr)\cong\Hom_{\fg}(L\otimes N_1,N_2)
\cong\Hom_{\fg}(N_1,N_2\otimes L^*)$$
for any $\fg$-modules $N_1,N_2$ and a finite dimensional
 $\fg$-module $L$.

\subsection{The category $\wdO$ and Verma modules}
\label{wdO}
Denote by $\cO$ the full subcategory of finitely generated
$\fh$-diagonalizable $\fg_0$-modules  which are locally
$\fn^+_0$-finite.
Denote by $\wdO$ the similarly defined category
of $\fg$-modules. 

Since $\cU(\fn^+)$ is finite over $\cU(\fn_0^+)$, 
a $\fg$-module $N$ belongs $\wdO$ iff as 
$\fg_0$-module $N$ belongs to $\cO$. In particular, 
any  module of category $\wdO$ has a finite length.

\subsubsection{}
\label{Vsmod}
For $\lambda\in {\frak h}^*$ denote by ${\Bbb C}_{\lambda}$
a simple ${\frak b}$-module such that
${\frak n}^+v=0$ and $hv=\lambda(h)v$ for any $h\in {\frak h},v\in  
{\Bbb C}_{\lambda}$. Define a Verma module $\wdM(\lambda)$ by setting
$$\wdM(\lambda):={\Ug}\otimes_{\cU({\frak b})}
{\Bbb C}_{\lambda}.$$
The module  $\wdM(\lambda)$ has
a unique simple quotient, which we denote by $\wdV(\lambda)$.
Similarly, denote by $M(\lambda)$ and $V(\lambda)$
respectively, Verma and simple $\fg_0$-modules of the highest weight
$\lambda$.

\subsubsection{}
\begin{defn}{typdef}
A weight $\lambda\in\fh^*$ is called {\em typical } if
$(\lambda+\rho,\beta)\not=0$ for any isotropic $\beta\in\Delta_1$.
\end{defn}

\subsubsection{}
\label{projO}
If $\lambda$ is typical then 
$\Ann_{\Zg} \wdM(\lambda)=\Ann_{\Zg} \wdM(\lambda')$ implies 
$\lambda'\in W.\lambda$--- see~\cite{kch}, Theorem 2. 
In particular, if a typical weight $\lambda$ is a minimal element
in its orbit $W.\lambda$ then $\wdM(\lambda)$ is simple.

On the other hand, if a typical weight $\lambda$ is a maximal element
in $W.\lambda$ then $\wdM(\lambda)$ is projective in $\wdO$.
Indeed, take a short exact sequence
$0\to \wdN'\to \wdN\to \wdM(\lambda)\to 0$
in $\wdO$.
One may assume that $\wdchi^{k}\wdN=0$ where 
$\wdchi:=\Ann_{\Zg} \wdM(\lambda)$ and $k$ is a positive integer.
Then  the weight of a primitive vector of
any simple subquotient of $\wdN$ belongs to $W.\lambda$. 
Since $\lambda$ is a maximal element of $W.\lambda$
one has 
$\Omega(\wdN)\cap\{\lambda+{\Bbb Z}\Delta^+\}
\subseteq \{\lambda-{\Bbb N}\Delta^+\}$.
Thus a preimage of a highest weight vector of $\wdM(\lambda)$
is primitive and the above exact sequence splits. 
 
By a similar argument,  a short exact sequence
$0\to \wdM(\lambda')\to \wdN\to \wdM(\lambda)\to 0$
in  $\wdO$ splits if $\lambda'\not>\lambda$.

\subsubsection{}
\begin{defn}{strtypdef}
A weight $\lambda\in\fh^*$ is called {\em strongly typical } if
$(\lambda+\rho,\beta)\not=0$ for any $\beta\in\Delta_1$.
\end{defn}
\subsubsection{}
Call the highest weight modules $\wdV(\lambda),\wdM(\lambda)$ 
typical (resp., strongly typical) if
$\lambda$ is typical (resp., strongly typical).

Apart from the cases $B(m,n)$ (that is ${\frak {osp}}(2m+1,2n)$)
and $G(3)$, all odd roots of $\fg$ are isotropic and thus
the notions of ``typical'' and ``strongly typical'' coincide.

\subsubsection{}
\label{Ndual}
Take $N\in\wdO$. Equip the graded dual vector space 
$N^{\#}:=\oplus_{\mu\in {\frak h}^*} (N|_{\mu})^*$ 
 by the following $\fg$-module structure
$$u.f(v):=(-1)^{d(f)d(u)}f(\sigma(u)v)), \ \ 
\forall u\in{\Ug},f\in N^{\#},v\in N$$
where $\sigma:\Ug\to\Ug$ is the  superantiautomorphism
defined in~\ref{cntrgr}. One can easily sees that
$N\mapsto N^{\#}$ defines a duality functor on $\wdO$.

Since the restriction of $\sigma$ on ${\frak h}$
is equal to identity, $\ch N^{\#}=\ch N$.
In particular, $\wdV(\lambda)^{\#}\cong \wdV(\lambda)$
and it is isomorphic to 
the socle of $\wdM(\lambda)^{\#}$.

\subsubsection{}
\label{zerointr}
According to~\cite{lm}, for any Zariski dense $S\subseteq \fh^*$
$$\displaystyle{\cap_{\lambda\in S}}\Ann\wdM(\lambda)=0.$$

\subsubsection{}
\label{shapform}
A Verma module is not irreducible iff its highest weight
is  a root of a {\em Shapovalov form}. Shapovalov forms
are polynomials in $\Sh$ indexed by
the weights of $\cU({\frak n}^-)$.
Each polynomial admits a linear factorization
which was established by Kac. The linear
factors of these polynomials are
$$\begin{array}{ll}
\alpha^{\vee}+(\rho,\alpha)-n(\alpha,\alpha)/2 &\ n\in {\Bbb N}^+,
\alpha\in\overline{\Delta}_{0}^+,\\
\alpha^{\vee}+(\rho,\alpha)-n(\alpha,\alpha)/2 &\ n\in 2{\Bbb N}+1,
\alpha\in(\Delta_{1}^+\setminus \overline{\Delta}_{1}^+),\\
\alpha^{\vee}+(\rho,\alpha), &
\alpha\in\overline{\Delta}_{1}^+.
\end{array}$$
Hence $\wdM(\lambda)$ is not irreducible iff 
$(\lambda+\rho,\alpha)= n(\alpha,\alpha)/2$
for a positive root $\alpha$ and a positive integer $n$
which should be odd for odd $\alpha$.

\subsection{Finite dimensional modules}
\label{tpc}
Necessary and sufficient conditions
for $\wdV(\lambda)$ to be finite dimensional are given
in~\cite{kadv}, Theorem 8. One can immediately sees
from these conditions that any typical finite dimensional 
$\wdV(\lambda)$  is strongly typical.
If $\wdV(\lambda),\wdV$ are finite dimensional
satisfying $\Ann_{\Zg} \wdV(\lambda)=\Ann_{\Zg} \wdV$
and $\wdV$ is typical then
$\wdV\cong\wdV(\lambda)$---see~\cite{kch1}, Prop. 2.7.
 
The following character formula of a typical finite dimensional
module is established by Kac (see~\cite{kch1},\cite{kch}): 
\begin{equation}
\label{tpcch}
\ch \wdV(\lambda)= D\displaystyle\sum_{w\in W} \sn (w)
e^{w.\lambda},\ \text{ where } 
D=\displaystyle\prod_{\alpha\in\Delta^+_0}(1-e^{-\alpha})^{-1}
\displaystyle\prod_{\beta\in\Delta^+_1}(1+e^{-\beta}). 
\end{equation}

\subsection{Odd reflections}
\label{oddrefl}
The odd reflections were introduced by I.~Penkov and V.~Serganova---
see~\cite{psg}, 3.1. If $\fg$ is a basic classical
Lie superalgebra, the definition takes the following form. 
Two Borel subsuperalgebras ${\frak b},{\frak b}'\subset {\fg}$
are {\em connected by an odd reflection along $\beta$}
iff $\beta$ is a simple odd isotropic
root of ${\frak b}'$ and 
$$\Delta({\frak b})=\{-\beta\}\cup\Delta({\frak b}')
\setminus\{\beta\}.$$
Note that $\Delta_0({\frak b})=\Delta_0({\frak b}')$
and so ${\frak b}_0={\frak b}'_0$.
As it is shown in~\cite{psg}, 3.1, any  Borel subsuperalgebras 
${\frak b},{\frak b}'\subset {\fg}$ satisfying 
${\frak b}_0={\frak b}'_0$ are connected by a chain of odd reflections.

Let ${\frak b}$ and ${\frak b}'$ be Borel
subsuperalgebras  connected by the odd reflection 
along an odd isotropic root $\beta.$ Then
${\frak b}'={\frak b}\cap {\frak b'}+{\Bbb C}x$
where $x$ is a non-zero element of $\fg_{\beta}$.
If $v$ is a ${\frak b}$-primitive vector such that
$xv=0$ then $v$ is also a ${\frak b}'$-primitive.
This implies  $\wdV_{{\frak b}}(\lambda)=\wdV_{{\frak b}'}(\lambda)$
for $\lambda$ satisfying $(\lambda,\beta)=0$.
If $xv\not=0$ then $xv$ 
is a ${\frak b}'$-primitive since $[x,x]=0$.
This implies $\wdV_{{\frak b}}(\lambda)=
\wdV_{{\frak b}'}(\lambda+\beta)$
and $\wdM_{{\frak b}}(\lambda)=\wdM_{{\frak b}'}(\lambda+\beta)$
for $(\lambda,\beta)\not=0$. Remark that $(\rho_{{\frak b}},\beta)=0$,
since $\beta$ is a simple isotropic root,
and so $(\lambda,\beta)=0$ iff $(\lambda+\rho_{{\frak b}},\beta)=0$.
Taking into account that $\rho_{{\frak b}}=\rho_{{\frak b}'}+\beta$,
one concludes that for Borel subsuperalgebras ${\frak b}$, ${\frak b}'$  
satisfying ${\frak b}_0={\frak b}'_0$ and a ${\frak b}$-typical
weight $\lambda$ one has
$\wdV_{{\frak b}}(\lambda)=\wdV_{{\frak b}'}(\lambda')$
and $\wdM_{{\frak b}}(\lambda)=\wdM_{{\frak b}'}(\lambda')$
where $\lambda+\rho_{{\frak b}}=\lambda'+\rho_{{\frak b}'}$.

\subsection{Completely reducible Lie superalgebras}
Recall that a Lie superalgebra is called completely reducible
if all its finite dimensional modules are
completely reducible. According to the theorem of    
Djokovi\'c and Hochschild (see~\cite{sch}, p. 239), any finite dimensional 
completely reducible Lie superalgebra is a direct sum of 
semisimple Lie algebras
and algebras ${\frak {osp}}(1,2l)$ ($l\geq 1$). 
The superalgebra $\fg=={\frak {osp}}(1,2l)$ has
many features of the semisimple Lie algebra; in particular,
$\Ug$ is a domain and $\Zg$ is a polynomial algebra.

\subsubsection{Separation Theorem}
\label{spr}
In~\cite{ko} Kostant establishes the following theorem which is
called sometimes ``Separation Theorem''.

\begin{thm}{separation} 
Let $\fg$ be a semisimple complex Lie algebra.
There exists an $\ad\fg$-invariant subspace ${H}$ in 
${\cal U}({\frak g})$ such that the multiplication map induces an
isomorphism
${\cal Z}({\frak g})\otimes {H}\iso{\cal U}({\frak g})$. 
Moreover, for every simple
finite dimensional module $V$, $[{H}:{V}]=\dim {V}_0$.
\end{thm}

Such an $\ad\fg$-invariant subspace $H$ is called a 
{\em harmonic space}.

In~\cite{mu} Musson proves the analogous theorem for 
$\fg={\frak {osp}}(1,2l)$.

\subsubsection{Annihilation Theorem}
\label{duflothm}
In~\cite{d} Duflo proves the following theorem.

\begin{thm}{ann}
Let $\fg$ be a semisimple complex Lie algebra. Then for any 
$\lambda\in \fh^*$
$$\Ann M(\lambda)=\cU(\fg) \Ann_{\Zg} M(\lambda).$$
\end{thm}

Let $\fg$ be a semisimple complex Lie algebra and $M$ be a Verma module.
The multiplicity of 
each simple finite dimensional module in $F(M,M)$ 
is equal to $\dim V|_0$ if $M$ is a $\fg_0$-simple Verma module.
Combining the above theorem and~\Thm{separation}, one concludes
that for such $M$ the natural map $\Ug/\Ann M(\lambda)\to F(M,M)$
is an isomorphism. In~\cite{j16}, 6.4 Joseph generalizes this
result to any $\fg_0$-Verma module.

\subsubsection{}
\label{annforosp12l}
For the case $\fg={\frak {osp}}(1,2l)$ 
the following results are obtained 
in~\cite{gl1},\cite{gl2}. 

\begin{thm}{annosp1}
The annihilator of a Verma module $\wdM$ coincides with
$\Ug \Ann_{\Zg} \wdM$
iff $\wdM$ is strongly typical.
\end{thm}

Using Separation Theorem, one concludes
that  for a  strongly typical $\wdM$ one has
$\Ug/(\Ann \wdM)\cong H$.

To describe the non-centrally generated annihilators of Verma
modules, it is convenient to substitute the centre $\Zg$
by the algebra $\tZg$ defined in~\ref{twist}.

\begin{thm}{annosp2}
For any $\lambda\in \fh^*$
$$\Ann \wdM(\lambda)=\Ug \Ann_{\tZg} \wdM(\lambda).$$
If $\lambda$ is not  strongly typical
then the ideal $\Ann_{\tZg} \wdM(\lambda)$ is a maximal ideal
of the algebra $\tZg$.
\end{thm}

The algebra $\Ug$ is free over $\tZg$; more precisely,
$\Ug$ contains $\ad\fg_0$-submodule $K$ such that
the multiplication map induces the $\ad\fg_0$-isomorphism
$K\otimes \tZg\to\Ug$. Moreover, as $\ad\fg_0$-modules
$H\cong K\oplus K$ and $K\cong \Ug/(\Ann \wdM)$
if $\wdM$ is not strongly typical .

\section{Category $\Fin$ and a twisted adjoint action}
\label{sectFin}
\subsection{Category $\Fin$}
\label{catFin}
Denote by $\Fin_0$ the full category of $\fg_0$-modules 
whose objects are sums of  simple finite dimensional modules. Denote by
$\Fin$ the full category of $\fg$-modules whose objects, considered as
$\fg_0$-modules, belong to  $\Fin_0$. Since $\Ug$ is a finite
extension of $\cU({\frak g}_0)$, any module $N\in\Fin$
is locally finite.
In other words, the objects of $\Fin$ are the locally finite 
$\fh$-diagonalizable modules ($\fg_0$ is reductive and so
all locally finite $\fh$-diagonalizable are completely reducible).

Denote by $\Irr$ the set of isomorphism classes of simple 
finite dimensional $\fg$-modules and by $\Irr_0$ 
the set of isomorphism classes of simple 
finite dimensional $\fg_0$-modules.
Note that $\Ug$ considered as $\ad\fg$-module belongs to  $\Fin$.

Throughout this section all modules are objects of $\Fin$.
Everywhere in  the paper, injectivity and projectivity
mean, by default,  injectivity and projectivity
in the category $\Fin$.

\subsubsection{}
\label{bgrsocle}
Recall that the socle $\Soc N$ of the module $N$ is the sum of its 
simple submodules. Since any module in $\Fin$ is locally finite,
it has a finite dimensional submodule and so
a non-trivial socle. 

Recall that a module $N$ is called {\em an essential extension} 
of its submodule $N'$ if for any non-zero
submodule $N''$ of $N$ one has $N'\cap N''\not=0$.
Any module in $\Fin$ is an essential extension of its socle.

\subsubsection{}
\label{socphi}
For a homomorphism $\psi: N\to N'$ denote by $\Soc\psi$ its
restriction to $\Soc N$. If $\Soc\psi$ is a monomorphism then
$\psi$ is also a monomorphism, since $\ker\psi\cap\Soc N=0$
and $N$ is an essential extension of $\Soc N$.

\subsection{Injective and projective objects in $\Fin$}
Recall that an injective envelope $\Inj(N)$ of the module
$N$ is an injective module which 
is an essential extension of $N$.  An injective envelope of
a given module is unique up to isomorphism.

A typical simple finite dimensional
module is injective and projective in $\Fin$ (see~\ref{typdef}).

A standard reasoning of~\cite{mcl}, 3.11 shows that any submodule $N$
of an injective module $E$ has an injective envelope
which is a direct summand of $E$.
In~\Lem{injsoc} we show that any module in $\Fin$ has 
an injective envelope.

\subsubsection{}
\label{definj}
A finite dimensional $\fg_0$-module $L$ 
is projective and injective in $\Fin_0$ and so
$\Coind_{\fg_0}^{\fg}L\cong\Ind_{\fg_0}^{\fg}L$
is projective and injective in $\Fin$--- see~\ref{defindcoind}.

Let $\wdV$ be a simple finite dimensional $\fg$-module
and $V$ be its simple $\fg_0$-submodule.
Then $\wdV$ is a submodule of  the injective module
$\Coind_{\fg_0}^{\fg}V$. Therefore $\wdV$ has an injective
envelope $\Inj(\wdV)$ which is a direct summand of 
$\Coind_{\fg_0}^{\fg}V$. The last is projective and
so $\Inj(\wdV)$ is projective as well.

\subsubsection{}
\begin{lem}{injsoc}
Any module $N\in\Fin$ has an injective envelope $\Inj(N)$.
Moreover $\Inj(N)\cong\Inj(\Soc N)$.
\end{lem}
\begin{pf}
By~\ref{bgrsocle}, $N$ is an essential extension of $\Soc N$.
Therefore, by~\cite{mcl}, 3.11.1, there exists a monomorphism
$\iota:N\to\Inj(\Soc N)$ whose socle is the natural embedding
$\Soc N\to  \Inj(\Soc N)$. Since $\Inj(\Soc N)$ is an essential
extension of $\Soc N$, it is also essential  extension of $N$.
Hence $\Inj(\Soc N)$ is an injective envelope of $N$.
\end{pf}

\subsubsection{}
\begin{lem}{suminj}
Any direct sum of injective modules in $\Fin$ is injective.
\end{lem}
\begin{pf}
Let $E_i, i\in I$ be a collection of injective modules.
One has to check that for any monomorphism $\iota: A\to B$
and any homomorphism $\phi: A\to \oplus_{i\in I} E_i$
there exists a homomorphism $\psi: B\to \oplus_{i\in I} E_i$
such that $\psi\iota=\phi$. 
By a standard reasoning based on Zorn's lemma, it is enough 
to verify the above assertion assuming $B$ being cyclic.
Any cyclic module in $\Fin$ is finite dimensional. If $B$
is finite dimensional, then $A$ is also finite dimensional
and so $\phi(A)$ lies in a finite subsum $\oplus_{i\in J} E_i$
($J$ is  a finite subset of $I$). Since each $E_i$ is injective,
$\oplus_{i\in J} E_i$ is injective and so
there exists $\psi: B\to \oplus_{i\in J} E_i$ such that 
$\psi\iota(a)=\phi(a)$ for all $a\in A$. 
The assertion follows.
\end{pf}

\subsubsection{}
\begin{cor}{cr1}
$$\Inj(\oplus_{i\in I} N_i)=\oplus_{i\in I}\Inj(N_i)$$
\end{cor}

\subsubsection{}
\begin{prop}{injprj}
Any injective module in $\Fin$ is projective.
\end{prop}
\begin{pf}
Take an injective module $N$ and write $\Soc N=\sum_{i\in I} L_i$
where the $L_i$ are simple.  Combining~\Lem{injsoc}
and~\Cor{cr1}, one obtains
$N\cong \oplus_{i\in I}\Inj(L_i)$. By~\ref{definj}, all
$\Inj(L_i)$ are projective so $N$ is also projective.
\end{pf}

\subsection{Twisted adjoint action}
\label{twist}
In this subsection we present some results from~\cite{ghost}
which are used in the sequel.

Define a twisted adjoin action of $\fg$ on $\Ug$ by setting
$$(\ad'g)u=gu-(-1)^{d(g)(d(u)+1)}ug, \ \forall g\in\fg,u\in\Ug.$$
Note that $\ad g=\ad' g$ for $g\in\fg_0$.

The anticentre $\Ag$ is the set of invariants of
$\Ug$ under the twisted adjoint action $\ad'\fg$.
The product of two anticentral elements is central.
For $\fg$ being a basic classical Lie superalgebra,
$\Ag$ is even and so  any anticentral element
commutes with the even elements of $\Ug$ and
anticommutes with the odd ones. Therefore ``the ghost centre''
$$\tZg:=\Zg+\Ag$$
is a commutative subalgebra of $\Ug$. 

For $\fg={\frak {osp}}(1,2l)$,
$\tZg$ is a polynomial algebra
and, moreover, $\tZg=\Zg\oplus T\Zg$
where $T$ is the element defined in~\ref{TUg}.

\subsubsection{}
\label{ad'}
Let $L$ be an $\ad\fg_0$-submodule of $\cU$.
Then the $\ad'\fg$-submodule generated by $L$ in $\Ug$
is isomorphic to the induced module $\Ind_{\fg_0}^{\fg} L$.
As $\ad'\fg$-module $\Ug$ is generated by $\cU$
and so $\Ug\cong\Ind_{\fg_0}^{\fg}\cU$.
Using the isomorphism
$\Ind_{\fg_0}^{\fg}\cU\cong \Coind_{\fg_0}^{\fg}\cU$,
one obtains a linear isomorphism $\cZ(\fg_0)\iso\Ag$ given
by  $z\mapsto (\ad' u)z$ where $u\in \Ug$ is such an element
that $u V(0)$ is the trivial $\fg$-submodule of 
$\Ind_{\fg_0}^{\fg} V(0)$.

An anticentral element $z$ acts on a module generated by
a primitive vector $v$ 
in the following way. It acts as $\cP(z)(\lambda)\id$
on $\Ug_0v$ and as $(-\cP(z)(\lambda)\id)$
on $\Ug_1v$ ($\Ug_0,\Ug_1$ are homogeneous components
of $\Ug$). From~\ref{zerointr} it follows that
the Harish-Chandra projection 
provides a monomorphism $\Ag\to\Sh^{W.}$.
The image of this monomorphism is equal
to $t\Sh^{W.}$ where 
$$t:=\prod_{\beta\in\Delta_1^+}(\beta^{\vee}+(\rho,\beta)).$$
Any non-zero element of $\Ag$ is a non-zero divisor in $\Ug$.

\subsubsection{}
\label{TUg}
Denote by $T$ the element of $\Ag$ such that $\cP(T)=t$.
Remark that a Verma module is strongly typical iff
its annihilator does not contain $T$.

Up to a non-zero scalar, $T$ is equal to the image of $1$ under
the above linear isomorphism $\cZ(\fg_0)\iso\Ag$.
The image of $T$ in the symmetric algebra is a non-zero element
of the one-dimensional vector space $\Lambda^{top}\fg_1$.

Remark that $T\Ag\subset \Zg$ and
the restriction of the Harish-Chandra projection provides
an isomorphism $T\Ag\iso t^2\Sh^{W.}$. 
For $\fg={\frak {osp}}(1,2l)$, this implies
that $\Ag$ is a free  module over $\Zg$ generated by $T$.

It is easy to check  that $(\ad g)(uT)=((\ad' g)u)T$. 
Since $T$ is a non-zero divisor, the ideal $\Ug T$
considered as $\ad\fg$-module is isomorphic to $\Ug$
considered as $\ad'\fg$-module. Hence as $\ad\fg$-module 
$\Ug T\cong\Ind_{\fg_0}^{\fg}\cU$
and, in particular, it is injective in $\Fin$. 
Note that $\Ug T$ is a two-sided ideal.
Remark that, apart from the case
when ${\frak g}$ is completely reducible, 
$\Ug$ itself as $\ad {\frak g}$-module
is neither injective nor projective in $\Fin$: it contains, 
as a direct summand, a trivial representation generated by $1$.

\subsubsection{}
For any $N\in\Fin_0$ one has
$$\dim\Hom_{\fg}(\wdV(\nu),\Coind_{\fg_0}^{\fg} N)=
\dim\Hom_{\fg_0} (\wdV(\nu),N).$$
Therefore
\begin{equation}\label{ind=suminj}
\Ind_{\fg_0}^{\fg} N\cong
\Coind_{\fg_0}^{\fg} N= \oplus_ {\wdV\in\Irr} 
\Inj(\wdV)^{r({\wdV})},\ \ \ 
r({\wdV}):=\dim \bigl(\Hom_{\fg_0} (\wdV,N)\bigr).
\end{equation}

This has the following useful consequence. Let $H$ be a harmonic
space of $\cU({\frak g}_0)$ that is an $\ad\fg_0$-submodule of
$\cU({\frak g}_0)$ such that 
the multiplication map provides an isomorphism
$H\otimes \cZ({\frak g}_0)\to\cU({\frak g}_0)$---see~\ref{spr}.
Then $H\cong \oplus_{V\in\Irr_0}
V^{\oplus \dim V|_0}$. 
Taking into account~\ref{ad'} and~(\ref{ind=suminj}), one obtains

\begin{equation}\label{indH}
(\ad'{\Ug})H\cong\Ind_{\fg_0}^{\fg} H
=\displaystyle\oplus_{V\in\Irr_0}
\Coind_{\fg_{0}}^{\fg} V^{\oplus \dim V|_0}
=\displaystyle\oplus_ {\wdV\in\Irr} \Inj(\wdV)^{\oplus\dim\wdV|_0}
\end{equation}
since $\sum_{V\in\Irr_0}\dim V|_0\cdot
\dim \bigl(\Hom_{\fg_0}(\wdV,V)\bigr)=\dim\wdV|_0$
for any $\wdV\in\Irr$.

\section{The structure $\Ug$ as $\ad\fg$-module.}
\label{sectPRV}
In this section we study the $\ad\fg$-module structure 
of $\Ug$.
We start from the studying $\Hom_{\frak g}(\wdV,{\Ug})$ 
for a simple finite dimensional module $\wdV$.
For each $\wdV$ we construct
a central element $z$ such that the localized module
$\Hom_{\frak g}(\wdV,{\Ug}[z^{-1}])=\Hom_{\frak g}(\wdV,{\Ug})[z^{-1}]$
is a free ${\Zg}[z^{-1}]$-module whose rank is equal 
to $\dim\wdV|_0$. 
In~\ref{sPRV} we define and study PRV determinants for
non-completely reducible basic classical Lie superalgebras. 

In~\ref{ssectsprthm} we show that
for a suitable $S\subset\Zg$ the localized algebra
$\Ug[S^{-1}]$ is free over its centre $\Zg[S^{-1}]$
and that the corresponding harmonic space $H$ (``generic harmonic
space'') is injective (in $\Fin$) $\ad\fg$-module.
Moreover, the multiplicity of a simple finite dimensional module $\wdV$
in $\Soc H$ is equal to $\dim\wdV|_0$.

\subsection{}
The following lemma of A.~Joseph and G.~Letzter
provides a connection between
$\Hom_{\frak g}(\wdV,{\Ug})$ and $\Hom(\wdV|_0,{\Sh})$.

\subsubsection{}
\begin{lem}{2ass}
An $\ad\fg$-submodule $N$ of $\Ug$ lies in the annihilator of 
$\wdV(\lambda)$ iff $\cP(N|_0)(\lambda)=0$.
\end{lem}
\begin{pf}
The proof is the same as in~\cite{jnato}, 7.2.
Remark that $\cP(N)=\cP(N|_0)$.
Let $v_{\lambda}$ be a highest weight vector of 
$\widetilde{V}(\lambda)$ and 
$\widetilde{V}(\lambda)_-\!:=\cU(\fn^-){\frak n}^-v$. 
One has
$$\cP(N)(\lambda)=0\Longleftrightarrow Nv_{\lambda}\subseteq 
{\widetilde{V}(\lambda)}_- .$$ 
In particular, $N \wdV(\lambda)=0$ forces $\cP(N)(\lambda)=0$.

For the inverse implication, assume that $\cP(N|_0)(\lambda)=0$
that is $Nv_{\lambda}\subseteq {\widetilde{V}(\lambda)}_-$.
The $\ad{\frak g}$-invariance of $N$ implies
${\cal U}({\frak n}^-)N=N{\cal U}({\frak n}^-)$ and thus
$$N\widetilde{V}(\lambda)=N{\cal U}({\frak n}^-)v_{\lambda}
={\cal U}({\frak n}^-)Nv_{\lambda}\subseteq {\cal U}({\frak n}^-)
{\widetilde{V}(\lambda)}_-\subseteq 
{\widetilde{V}(\lambda)}_-\ssubset{\widetilde{V}(\lambda)}.$$
The $\ad{\frak g}$-invariance of $N$ implies also
${\cal U}({\frak g})N=N{\cal U}({\frak g})$ and so
$N\widetilde{V}(\lambda)$ is a submodule of $\widetilde{V}(\lambda)$.
Hence $N{\widetilde{V}(\lambda)}=0$ as required.
\end{pf}

\subsubsection{}
\begin{cor}{cor2ass}
Take $\wdV\in\Irr$, a basis $v_1,\ldots v_r$ be  of $\wdV|_0$
and $\lambda\in\fh^*$. For any
$\theta_1,\ldots,\theta_k\in \Hom_{\frak g}(\wdV,{\Ug})$,
the image of the space $\sum_{i=1}^k\theta_i(\wdV)$ 
under the natural map ${\Ug}\to \End(\wdV(\lambda))$
is isomorphic to $\wdV^{\oplus m}$ where $m$ is the rank
of the matrix 
$\bigl(\cP(\theta_j(v_i))(\lambda)\bigr)_{i=1,r}^{j=1,k}$.
\end{cor}

\subsubsection{}
\label{3ass}
Combining~\ref{zerointr} and~\Lem{2ass} one concludes
that for an  $\ad\fg$-submodule $N$ of $\Ug$ 
the equality $\cP(N|_0)=0$ implies $N=0$.
Moreover $N=0$ provided that $\cP(N|_0)(R)=0$ for a 
Zariski dense subset $R$ of ${\frak h}^*$.

\subsection{Notation}
\label{defuu}
Fix $\wdV\in\Irr$ and consider
$\Hom_{\frak g}(\wdV,{\Ug})$ as a $\Zg$-module with respect to the 
action induced by the multiplication.

The Harish-Chandra projection induces the map 
$\Psi:\Hom_{\frak g}(\wdV,{\Ug})\to \Hom_{\frak g}(\wdV|_0,{\Sh})$ 
given by 
$$\Psi(\phi(v)=\cP(\phi(v)),\ \forall 
\phi\in \Hom_{\frak g}(\wdV,{\Ug}), v\in\wdV|_0.$$
This map  is a monomorphism by~\ref{3ass}.
Denote by $\uu(\wdV)$ the image of $\Psi$.

Define the action of $\Zg$ on $\Sh$ by setting
$zp:=\cP(z)p$. This action induces the structure 
of a $\Zg$-module on $\Hom(\wdV|_0,{\Sh})$.
Obviously, $\Psi$ is a $\Zg$-map. Thus the study
of $\Zg$-module structure of $\Hom_{\frak g}(\wdV,{\Ug})$
reduces to the study of $\Zg$-module structure of 
$\uu(\wdV)$. 

The vector space $\Hom(\wdV|_0,{\Sh})$ has the natural structure
of $\Sh$-module. We denote by $\uu(\wdV)\Sh$ the $\Sh$-span
of $\uu(\wdV)$ inside $\Hom(\wdV|_0,{\Sh})$.

\subsection{Examples}

\subsubsection{}
\begin{exa}{}
Let $\fg$ be a completely reducible
simple Lie superalgebra. In this case $\cP(\Zg)={\Sh}^{W.}$.
Separation Theorem (see~\ref{spr})
state the existence of an $\ad{\fg}$-submodule $\cH$ of $\Ug$
such that the multiplication map induces isomorphism 
${\cH}\otimes{\Zg}\iso {\Ug}$. Thus for any $\wdV\in\Irr$,
a basis of the vector space
$\Hom_{\fg}(\wdV,{\cH})$ is a free $\Zg$-basis of 
$\Hom_{\fg}(\wdV,\Ug)$. Consequently,
$\Hom_{\fg}(\wdV,\Ug)$ is a free $\Zg$-module of the rank
$\dim \Hom_{\fg}(\wdV,{\cH})=\dim\wdV|_0$ (see~\ref{spr}).
\end{exa}

\subsubsection{}
\begin{exa}{exaPRVsl21}
Consider the case ${\frak g}={\frak sl}(2,1)$. 
Then $\fh=\spn\{ z,h\}$ where
$z$ is a central element of the reductive algebra 
$\fg_0=\fsl(2)\times {\Bbb C}$ and $h$ is an element of 
the Cartan subalgebra of ${\frak sl}(2)$. 
Then $\Delta_0=\{\pm\alpha\}, \Delta_1=\{\pm\beta;\pm(\alpha+\beta)\}$.
Choose such a triangular decomposition
that $\Delta^+=\{\alpha,\beta,\alpha+\beta\}$.
The subalgebras $\fn_1^{\pm}=\fg_1\cap\fn^{\pm}$ are supercommutative 
and are $\ad\fg_0$-submodules of $\fg$.
Let $\{x_1,x_2\}$ (resp., $\{y_1,y_2\}$) be a basis of 
$\fn_1^+$ (resp., $\fn_1^-$). 

The algebra $\Sh^{W.}$ is a polynomial algebra generated by $z$ and 
$t=(z-h)(z+h+2)$.
The image $\cP(\Zg)$ in $\Sh^{W.}$ is spanned by $1$ and
the elements $\{t^nz^k,n>0,k\geq 0\}$. 

Denote by $v$ a highest weight vector
of $\wdV:=\wdV(\alpha+2\beta)$; one has $\wdV|_0={\Bbb C}y_1y_2v$.
It is easy to see that a highest weight vector of any copy of 
$\wdV$ inside $\Ug$ is of the form 
$ux_1x_2$ where $u\in\cZ(\fg_0)$.  Up to a scalar,
$\cP((\ad y_1y_2) (ux_1x_2))=(z-h)(z+h+2)\cP(u)=t\cP(u)$. 
One has $\cP(\cZ(\fg_0))=\Sh^{W.}$ and so
$\uu(\wdV)=\Hom ({\Bbb C}, \Sh^{W.}t)$.
Since $\Sh^{W.}$ is not free over $\cP(\Zg)$, the $\Zg$-module
$\uu(\wdV)\cong \Hom_{\fg}(\wdV,\Ug)$ is not free.

One might expect from the above example that $\uu(\wdV)$ is stable
with respect to the multiplication on $\Sh^{W.}$.
However $\uu (\wdV(0))=\Hom ({\Bbb C},\cP(\Zg))$ is not
stable with respect to the multiplication on $\Sh^{W.}$ apart from 
the case when $\cP(\Zg)=\Sh^{W.}$.
\end{exa}

\subsection{}
In this subsection we show that the $\Zg$-rank of the module
$\Hom_{\fg}(\wdV,\Ug)$ is not greater than $\dim \wdV|_0$.
The crucial point is~\Lem{lemfractPRV} asserting
that any elements of $\uu(\wdV)$
which are ``linearly dependent'' over $\Sh$ are 
``linearly dependent'' over $\Zg$.
Throughout this subsection $\wdV\in\Irr$ is fixed.

\subsubsection{Basic definitions}
\label{basicdef}
Let $A$ be a commutative domain.
For an $A$-module $N$ define {\em an $A$-rank } 
of $N$ to be the dimension over the field of fractions $\Fract A$
of the localized module $N\otimes_{A}\Fract A$.
Call elements $\theta_1,\ldots,\theta_k\in N$ {\em $A$-linearly
independent } if their images in the localized module 
$N\otimes_{A}\Fract A$ are linearly independent. Call
the elements {\em $A$-linearly
dependent } if they are not $A$-linearly
independent.
Call {\em an $A$-basic system } of $N$ a collection 
$\theta_1,\ldots,\theta_k\in N$ such that the image of 
this collection in the localized module $N\otimes_{A}\Fract A$
forms its $\Fract A$-basis.

Let $L$ be a vector space and
$v_1,\ldots,v_r$ be a basis of $L$. Let
$N$ be an $A$-submodule of $\Hom(L,A)$ and $k$ be 
the $A$-rank of $N$.
The collection $\theta_1,\ldots,\theta_k\in N$ is 
a $A$-basic system of $N$ iff the matrix 
$\bigl(\theta_j(v_i)\bigr)_{i=1,r}^{j=1,k}
\in {\cal Mat}_{(r\times k)}(A)$ has a 
$(k\times k)$ non-zero minor. Call such a non-zero
minor (which is an element of $A$) {\em a minor of a basic system}
$\theta_1,\ldots\theta_k$. 

If $p\in A$ is a minor of a basic system 
$\theta_1,\ldots\theta_k$ then for any
$\theta\in N$ one has $p\theta=\sum_{j=1}^k p_j\theta_j$
for certain (unique) collection $p_1,\ldots,p_k\in A$. In particular,
the localized module $N\otimes_A A[p^{-1}]$ is 
a free $A[p^{-1}]$-module and the images of $\theta_1,\ldots\theta_k$ 
form a basis of this free module.

Recall that $t^2\Sh^{W.}\subseteq\cP(\Zg)$ (see~\ref{TUg}).

\subsubsection{}
\begin{lem}{lemfractPRV}

(i) Elements of $\uu(\wdV)$ are $\Zg$-linearly independent
iff they are $\Sh$-linearly independent.

(ii) If $\theta_1\ldots,\theta_s\in\uu(\wdV)$
are $\Sh$-linearly independent and
$$\sum_{j=0}^s p_j\theta_j=0$$
for some $\theta_0\in\uu(\wdV)$
and $p_0,\ldots,p_s\in\Sh$,
then there exists $z_0,\ldots,z_s\in\Zg$
such that $\sum_{j=0}^s z_j\theta_j=0$
and $\cP(z_0)=t^2q$ where $q$ is a  
maximal $W.$-invariant divisor of $p_0^{|W|}$.
\end{lem}
\begin{pf}
Obviously, elements of $\uu(\wdV)$ are $\Zg$-linearly independent
provided that they are $\Sh$-linearly independent. The inverse
implication follows from (ii) because we always can choose
a minimal (with respect to the inclusion) subset of 
$\Sh$-linearly dependent elements in a set of 
$\Sh$-linearly dependent elements.

Let us prove (ii). 
The equality $\sum_{j=0}^s p_j\theta_j=0$
means that for any $\mu\in\fh^*$ one has
$$\forall v\in\wdV|_0\ \ \ 0=\sum_{j=0}^s p_j(\mu)\theta_j(v)(\mu)=
\sum_{j=0}^s p_j(\mu)\cP((\Psi^{-1}\theta_j)(v))(\mu).$$
In the light of~\Lem{2ass}, this gives
$$\forall\mu\in\fh^* \ \ \ 
\bigl(\sum_{j=0}^s p_j(\mu)\Psi^{-1}(\theta_j)\bigr)(\wdV)
\subset\Ann\wdV(\mu);$$
here the sum  belongs to $\Hom_{\fg}(\wdV,\Ug)$
because $\Psi^{-1}(\theta_j)\in\Hom_{\fg}(\wdV,\Ug)$ and
$p_j(\mu)\in {\Bbb C}$. 

Fix a root $\alpha\in\pi_0$; let $s_{\alpha}\in W$ be the corresponding
reflection. Assume that $\mu\in\fh^*$ is such that 
$\wdV(\mu)=\wdM(\mu)/\wdM(s_{\alpha}.\mu)$. Then, for ``sufficiently
large'' $\mu$, any copy of $\wdV$ inside $\Ug$ which annihilates
$\wdV(\mu)$, annihilates also $\wdM(\mu)$ and so $\wdV(s_{\alpha}.\mu)$.
To be more precise, choose $\omega\in {\Bbb Q}\Delta$ such
that $(\omega,\alpha')\geq 0$ for all $\alpha'\in\Delta^+$ and
$(\omega,\alpha)=1$. Set 
$$R:=\{\mu\in \fh|\ (a)\ n(\mu)\in {\Bbb N}^+, \ (b)\ 
n(\mu)>(-\nu,\omega),\ (c)
\forall\beta\in\Delta^+\setminus\{{\Bbb Q}\alpha\}
\ \ (\mu+\rho,\beta)\not\in {\Bbb Q}\}.$$
where $\nu$ stands for the lowest weight of $\wdV$
and $n(\mu):=2(\mu+\rho,\alpha)/(\alpha,\alpha)$.
It is easy to see that $R$ is Zariski dense in $\fh^*$
and that any element in $R$ is typical. Take $\mu\in R$.
By~\ref{shapform}, $\wdM(\mu)$ is not simple 
because $n(\mu)\in {\Bbb N}^+$. 
If $\wdM(\mu')$ is a subquotient of $\wdM(\mu)$, then
$(\mu'-\mu)\in {\Bbb N}\Delta^-$ and, by~\ref{projO},
$\mu'=w.\mu$ for some $w\in W$. Using the fact that
any $w\in W$ can be written as the product  
$s_{\beta_1}\cdots s_{\beta_m}$ where
$\beta_1,\ldots,\beta_m\in\Delta_0^+$ are
linearly independent, it is easy to deduce
from the condition (c)
that  $(w.\mu-\mu)\not\in {\Bbb N}\Delta^-$ for any 
$w\not=\id,s_{\alpha}$. 
Thus $\wdV(\mu)=\wdM(\mu)/\wdM(s_{\alpha}.\mu)$ and 
$\wdM(s_{\alpha}.\mu)$ is simple.
Observe that $s_{\alpha}.\mu=\mu-n(\mu)\alpha$ and 
so for any $\xi\in\Omega(\wdM(s_{\alpha}.\mu))$
one has $(\xi,\omega)\leq (\mu,\omega)-n(\mu)$.
Thus
$(\mu+\nu)\not\in \Omega(\wdM(s_{\alpha}.\mu))$ due to the condition
(a).

Let $N$ be a copy of $\wdV$ inside $\Ann \wdV(\mu)$.
Let $u_{\nu}$ be a lowest weight vector of $N$
and $v_{\mu}$ be a highest weight vector of $\wdM(\mu)$.
Since $u_{\nu}\wdV(\mu)=0$ the vector $u_{\nu}v_{\mu}$ belongs
to $\wdM(s_{\alpha}.\mu)$. However the latter 
does not have non-zero elements
of weight $(\mu+\nu)$. Thus $u_{\nu}v_{\mu}=0$. 
One has $u_{\nu} \wdM(\mu)=u_{\nu} \cU(\fn^-)v_{\mu}
=\cU(\fn^-)u_{\nu}v_{\mu}$ because  $u_{\nu}$ is a lowest weight vector
that is $(\ad\fn^-)u_{\nu}=0$. Therefore $u_{\nu}\wdM(\mu)=0$.
Since $\Ann \wdM(\mu)$ is $\ad\fg$-stable,
$N\subset \Ann \wdM(\mu)$ 
and, in particular, $N\subset \Ann \wdV(s_{\alpha}.\mu)$.

Hence $\sum_{j=0}^s p_j(\mu)(\Psi^{-1}\theta_j)(\wdV)$
annihilates $\wdV(s_{\alpha}.\mu)$ for any $\mu\in R$.
Then, by~\Lem{2ass}
$$\left(\sum_{j=0}^s (s_{\alpha}.p_j)\theta_j(v)\right)(\mu')=0$$
for any $v\in \wdV|_0$ and any
$\mu'$ such that $s_{\alpha}.\mu'\in R$. The terms 
$\sum_{j=0}^s (s_{\alpha}.p_j)\theta_j(v)$ are polynomials in $\Sh$.
Since $R$ is a Zariski dense subset of $\fh^*$,
one concludes that these polynomials are equal to zero.
Consequently
$$\sum_{j=0}^s (s_{\alpha}.p_j)\theta_j=0.$$
Taking into account that $\sum_{j=0}^s p_j\theta_j=0$
and that $\theta_1,\ldots\theta_s$ are $\Sh$-linearly independent,
one concludes that
$p_j/p_0=(s_{\alpha}.p_j)/(s_{\alpha}.p_0)$ for all 
$j=1,\ldots,s$. Thus $p_j/p_0$ is $W.$-invariant
for all $j=1,\ldots,s$. Then, by~\Lem{invfract}, for each $j=0,\ldots,s$
there exists $q_j\in\Sh^{W.}$ such that $q_j/q=p_j/p_0$.
By~\ref{TUg}, $t^2 {\Sh}^{W.}\subset \cP(\Zg)$ 
and thus  $t^2q_j\in \cP(\Zg)$ for all $j=0,\ldots,s$.
Since $q_j/q=p_j/p_0$ one has
$$(t^2q)\theta_0+\sum_{j=1}^s (t^2q_j)\theta_j=0.$$
This completes the proof.
\end{pf}

\subsubsection{}
\begin{prop}{fractPRV}
Let $\theta_1,\ldots,\theta_s\in \uu(\wdV)$ be a 
$\Sh$-basic system of $\uu(\wdV)\Sh$, $p$ 
be a minor of this system
and $z\in {\Zg}$ be such that $\cP(z)=t^2q$ where $q$ is a  maximal
$W.$-invariant divisor of $p^{|W|}$. Then  the  localized module 
$\Hom_{\fg}(\wdV,{\Ug}[z^{-1}])$ is freely generated over 
${\Zg}[z^{-1}]$ by  $\{\Psi^{-1}\theta_1,\ldots,\Psi^{-1}\theta_s\}$.
\end{prop}
\begin{pf}
Recall that $\Psi$ provides a $\Zg$-isomorphism from 
$\Hom_{\fg}(\wdV,{\Ug})$
onto $\uu(\wdV)$. Therefore the required assertion is equivalent to
the following statement: the localized module 
$\uu(\wdV)\otimes_{\Zg}{\Zg}[z^{-1}]$ is freely generated over 
${\Zg}[z^{-1}]$ by  $\{\theta_1,\ldots,\theta_s\}$. 
In other words, one has to show that for any
$\theta_0\in \uu(\wdV)$ there exists a unique collection
$z_1,\ldots,z_s$ such that 
$z\theta_0+\sum_1^s z_j\theta_j=0$ or,
equivalently, that there exists a unique collection
$a_1,\ldots,a_s\in\cP(\Zg)\subset\Sh$ such that 
$$(t^2q)\theta_0+\sum_1^s a_j\theta_j=0.$$
The uniqueness immediately follows from the fact that
the collection $\theta_1,\ldots,\theta_s$ is a $\Sh$-basic system of
$\uu(\wdV)\Sh$. To prove the existence, 
recall that $\uu(\wdV)\Sh$ is a $\Sh$-submodule of $\Hom(\wdV|_0,\Sh)$.
By~\ref{basicdef}, for any $\theta_0\in \uu(\wdV)$
there exists 
a unique collection $p_1,\ldots,p_s\in {\Sh}$ such that
$p\theta_0+\sum_{j=1}^s p_j\theta_j=0$.
Now the required existence follows from~\Lem{lemfractPRV}.
\end{pf}

\subsubsection{}
\begin{cor}{corfr}
The $\Zg$-rank of $\Hom_{\fg}(\wdV,\Ug)$ is not 
greater than $\dim \wdV|_0$.
\end{cor}

\subsubsection{}
\begin{rem}{frfr}
Set $r=\dim\wdV|_0$. Assume
that $\Hom_{\fg}(\wdV,{\Ug})$ contains the collection
$\theta_1,\ldots,\theta_r$ possessing the following
property: for a certain $\lambda\in {\frak h}^*$
the images of the modules
$\theta_1(\wdV),\ldots,\theta_r(\wdV)$
under the map ${\Ug}\to {\Ug}/\Ann \wdV(\lambda)$
form a direct sum. We claim that $\theta_1,\ldots,\theta_r$
is a $\Zg$-basic system of $\Hom_{\fg}(\wdV,{\Ug})$.

In fact, choose
a basis $v_1\ldots v_r$ of $\wdV|_0$. By~\Cor{cor2ass},
the rank of the matrix 
$\bigl(\cP(\theta_j(v_i))(\lambda)\bigr)_{i=1,r}^{j=1,r}$
is equal to $r$. Hence the rank of the matrix 
$$\bigl(\Psi(\theta_j(v_i))\bigr)_{i=1,r}^{j=1,r}=
\bigl(\cP(\theta_j(v_i))\bigr)_{i=1,r}^{j=1,r}\in Mat_{(r\times r)}{\Sh}$$
is also equal to $r$ and so $\Psi\theta_1,\ldots,\Psi\theta_r$
are $\Sh$-linearly independent in $\uu(\wdV){\Sh}$. 
Since $\uu(\wdV){\Sh}$ is a $\Sh$-submodule of $\Hom(\wdV|_0,{\Ug})$,
its rank is not greater than $r$. Hence 
$\Psi\theta_1,\ldots,\Psi\theta_r$ is a  $\Sh$-basic
system of $\uu(\wdV){\Sh}$. 
Therefore $\theta_1,\ldots,\theta_r$ is a $\Zg$-basic
system of $\Hom_{\fg}(\wdV,{\Ug})$ by~\Prop{fractPRV}.
\end{rem}

\subsection{}
\label{instead}
In this subsection we show that the $\Zg$-rank of $\Hom_{\fg}(\wdV,\Ug)$
is equal to $\dim \wdV|_0$.

\subsubsection{}
Separation Theorem~\ref{spr} claims  the existence
an $\ad\fg_0$-submodule $H$ of $\cU$
such that the multiplication map provides an isomorphism
$H\otimes\cZ({\frak g}_0)\iso \cU$.
In both proofs (\cite{ko},~\cite{bl}) one constructs, actually,
an $\ad\fg_0$-submodule $H'$ of the symmetric algebra 
$\cS({\frak g}_0)$ such that the multiplication map 
provides an isomorphism
$H'\otimes \cS({\frak g}_0)^{{\frak g}_0}\iso \cS({\frak g}_0)$.

\begin{lem}{hrmnew}
Let $H$ be an $\ad\fg_0$-submodule of $\cU$
such that the multiplication map provides an isomorphism
$\gr H\otimes \cS({\frak g}_0)^{{\frak g}_0}\to \cS({\frak g}_0)$.
Then the $\ad\fg$-module $L:=(\ad\Ug)(HT)$ is isomorphic to
$$\displaystyle\oplus_{\wdV\in\Irr} \Inj(\wdV)^{\oplus\dim\wdV|_0}$$
 and the multiplication induces the monomorphism $L\otimes \Zg\to \Ug$.
\end{lem}
\begin{pf}
Fix $H$ satisfying the above condition.
According to~\ref{twist}, $L:=(\ad\Ug)(HT)=(\ad'\Ug)(H)T$
and as $\ad\fg$-module $L$ is isomorphic to the
$\ad'\fg$-module $(\ad'\Ug)(H)\cong \Ind_{\fg_0}^{\fg}H$.
Now the required isomorphism follows from~(\ref{indH}).

Since $L=(\ad'\Ug)(H)T$ and $T$ is a non-zero divisor,
the injectivity of the map $L\otimes \Zg\to\Ug$ is equivalent
to the injectivity of the map $(\ad'\Ug)(H)\otimes \Zg\to\Ug$ 
(both maps are restrictions of the multiplication map).
To check the injectivity of the map $(\ad'\Ug)(H)\otimes \Zg\to\Ug$ 
it is enough to check the injectivity of its ``graded version'':
$$\gr((\ad'\Ug)(H))\otimes \gr\Zg\to\Sg.$$

Recall that $\gr\Zg=\Sg^{\fg}$. 
Let $\{x_i\}_{i\in I}$ be an ordered basis of $\fg_1$.
For any $J\subseteq I$ set $x_J:=\prod_{j\in J} x_j$
where the product is taken with respect to the given order.  
Assume that
$$\displaystyle\sum_{J\subseteq I; k=1,\ldots,s} 
z_{k,J}\gr ((\ad' x_{J}) h_k)=0$$
for some non-zero $h_1,\ldots,h_s\in H$ and 
some $z_{k,J}\in \Sg^{\fg}$. To check the injectivity
one has to show that all elements
$z_{k,J}$ are equal to zero. Suppose not. 
For  $x\in \fg_1, u\in\Ug$ one has $(\ad'x)u=2xu-(\ad x)u$.
This implies
$\gr((\ad' x_{J}) u)=2^{|J|}\gr x_{J}\gr u$
for any $J\subseteq I$ and $u\in\cU$.
For each $J\subseteq I$ denote by $P_J$ the projection
$\Sg\to \cS(\fg_0)\gr x_J$ with the kernel
$N_J:=\sum_{J'\subseteq I,J'\not=J}\cS(\fg_0)\gr x_{J'}$.  
By~\ref{syminv},  the restriction of $P_{\emptyset}$ to $\Sg^{\fg}$
provides a monomorphism $\Sg^{\fg}\to \cS(\fg_0)^{\fg_0}$.
Choose $J\subseteq I$ such that $z_{k,J}\not=0$ for some $k$
and $J$ is a minimal (under inclusion)
subset of $I$ possessing this property. 
Then 
$$\begin{array}{ll}
0&=P_J\bigl(\displaystyle\sum_{J'\subseteq I, k=1,\ldots,s} 
z_{k,J'}\gr  ((\ad' x_{J'}) h_k)\bigr)\\
&=
P_J\bigl(\displaystyle\sum_{J'\subseteq I, k=1,\ldots,s} 
2^{|J'|}z_{k,J'}\gr x_{J'}\gr h_k)\bigr)\\
&=2^{|J|}\gr x_{J}
\displaystyle\sum_{k=1,\ldots,s}P_{\emptyset}(z_{k,J})\gr h_k.
\end{array}$$
Since  $\gr H\otimes \cS(\fg_0)^{\fg_0}\iso \cS(\fg_0)$,
all elements $z_{k,J}$ are equal to zero.
The lemma is proven.
\end{pf}

\subsubsection{}
\begin{cor}{cordetnn}
The $\Zg$-rank of $\Hom_{\fg}(\wdV,\Ug)$
is equal to $\dim \wdV|_0$.
\end{cor}
\begin{pf}
Fix $\wdV\in\Irr$ and set $r:=\dim \wdV|_0$.
Choose $H$ satisfying the assumption of~\Lem{hrmnew}.
Then the vector space
$\Hom_{\fg}(\wdV,L)$ has dimension $r$; let
$\theta_1,\ldots,\theta_r$ be a basis of 
$\Hom_{\fg}(\wdV,L)$.
The injectivity of the map $L\otimes \Zg\to\Ug$, 
induced by the multiplication, implies the $\Zg$-linear independence
of $\theta_1,\ldots,\theta_r$ in $\Hom_{\fg}(\wdV,\Ug)$.
Thus the $\Zg$-rank of $\Hom_{\fg}(\wdV,\Ug)$ is greater than
or equal to $r$.
Comparing with~\Cor{corfr}, one concludes the required assertion.
\end{pf}

\subsection{PRV determinants}
\label{sPRV}
\subsubsection{Definitions}
\label{defprv}
Fix $\wdV\in\Irr$ and set $r:=\dim\wdV|_0$. \Lem{lemfractPRV} (i) 
implies that $\theta_1,\ldots,\theta_r$
is a $\Zg$-basic system of $\Hom_{\fg}(\wdV,{\Ug})$
iff $\Psi\theta_1,\ldots,\Psi\theta_r$ 
is a $\Sh$-basic system of $\uu(\wdV)\Sh$.
By~\ref{basicdef}, the last is equivalent to the condition
that $\det C\not=0$ where
$$C:=\bigl(\cP(\theta_j(v_i))\bigr)\in Mat_{(r\times r)}({\Sh})$$
and $\{v_1,\ldots,v_r\}$ is a basis of $\wdV|_0$.

Call such a matrix {\em a PRV matrix} (corresponding to
the basic system $\theta_1,\ldots,\theta_r$) and 
the determinant $\det C$ {\em a PRV-determinant}.
Note that for different choices of a basis 
$\{v_1,\ldots,v_r\}$ of $\wdV|_0$ 
the PRV-matrix corresponding
to $\theta_1,\ldots,\theta_r$ differ by the multiplication
on an invertible scalar matrix.

Denote by $\pprv(\wdV)$ the set of all PRV-matrices
and by $\eprv(\wdV)$ the set of PRV determinants
(for fixed $\wdV\in\Irr$). One has
$$\eprv(\wdV)=\{\det \bigl(\cP(\theta_j(v_i))\bigr)|\ 
\theta_1,\ldots,\theta_r\in \Hom_{\fg}(\wdV,{\Ug})\}\setminus\{0\}.$$

For each $\nu\in\fh^*$ such
that $\dim\wdV(\nu)<\infty$
set $\eprv (\nu):=\eprv (\wdV(\nu))$.
For $\lambda\in {\frak h}^*$ we write $\eprv(\wdV)(\lambda)=0$
if $(\det C)(\lambda)=0$ for all $C\in \pprv(\wdV)$.

If $\theta_1,\ldots,\theta_r$ is a $\Zg$-basic system 
of $\Hom_{\fg}(\wdV,{\Ug})$, then for any non-zero $z\in\Zg$ the collection
$z\theta_1,\theta_2\ldots,\theta_r$ is also a $\Zg$-basic system.
Consequently, $\eprv(\wdV)$
is closed under the multiplication on the non-zero
elements of $\cP({\Zg})$.

\subsubsection{}
\label{PRVVN}
Let $N$ be an $\ad\fg$-submodule of $\Ug$ such that 
$\dim\Hom_{\fg}(\wdV,N)=\dim\wdV|_0$.
By slightly abuse of notation, we shall denote
by $\eprv(\wdV;N)$ the determinant of the matrix
$\bigl(\cP(\theta_j(v_i))\bigr)$
where $\theta_1,\ldots,\theta_r$ (resp, $v_1,\ldots, v_r$)
is a basis of $\Hom_{\fg}(\wdV,N)$ (resp., $\wdV|_0$).
Note that $\eprv(\wdV;N)$ can be zero.

For different choices of a basis 
$\theta_1,\ldots,\theta_r$, the values of $\eprv(\wdV;N)$ 
differ by a multiplicative constant. Hence
$\eprv(\wdV;N)$ is a polynomial in $\Sh$ defined up
to a non-zero scalar.

\subsubsection{}
\begin{rem}{origdef}
The original definition of PRV determinants given in~\cite{prv}
for the semisimple Lie algebras differ from our definition.
Namely, the PRV determinant $p(V)$ is equal to
the polynomial  $\eprv(V;H)$
where $H$ is a harmonic space (see~\ref{spr}).
This definition works for completely
reducible Lie superalgebras as well.
In our notation $\eprv(V)=\Sh^{W.}p(V)\setminus\{0\}$---
this follows from Separation~\Thm{spr}.
\end{rem}

\subsubsection{}
\label{chuchu}
Take  $\wdV\in\Irr$ and $p\in \eprv(\wdV)$. Let $\theta_1,\ldots,\theta_r$
be a $\Zg$-basic system of $\Hom_{\fg}(\wdV,{\Ug})$
such that the determinant of the corresponding PRV matrix is
equal to $p$. Denote by $q$ the maximal $W.$-invariant
divisor of $p$.
Combining~\Rem{frfr} and~\Prop{fractPRV},
one concludes that $z\in{\Zg}$ satisfying $\cP(z)=t^2q^{|W|}$
possesses the following property: for any 
$\theta\in\Hom_{\fg}(\wdV,\Ug)$ one has
$z\theta=\sum_{i=1}^r z_i\theta_i$ for certain $z_i\in {\Zg}$. 

In particular, for $a:=\cP(z)^r=t^{2r}q^{r|W|}$ one has
$$a\cdot\eprv(\wdV)\subset \cP({\Zg})p\subseteq{\Sh}^{W.}p.$$
On the other hand, $\cP(\Zg)p\subseteq \eprv(\wdV)$
by~\ref{defprv}.

\subsubsection{}
\begin{cor}{nonzeroPRV}
Take $\wdV\in\Irr$. For  any $p\in \eprv(\wdV)$
there exists $a\in {\Sh}^{W.}$ such that 
$$\cP(\Zg\setminus\{0\})p\subseteq \eprv(\wdV)
\subset {\Sh}^{W.}[a^{-1}]p.$$
\end{cor}

\subsubsection{}
\label{corankprv}
Fix $\wdV\in\Irr$. Take a basic system $\theta_1,\ldots,\theta_r$ and
denote by $C$ the corresponding PRV matrix.
Recall that $\det C\not=0$.
For each $\lambda\in {\frak h}^*$ denote by $C(\lambda)$
the complex matrix which is obtained
from $C$ by the evaluation of all entries at $\lambda$. 
Clearly $\det C(\lambda)=(\det C)(\lambda)$.
By~\Cor{cor2ass}, 
$$\corank C(\lambda)=\dim\Hom_{\fg}
\bigl(\wdV, \Ann\wdV(\lambda)\cap\sum_{j=1}^r\theta_j(\wdV)\bigr).$$
In particular, 
$\Ann\wdV(\lambda)\cap \sum_{j=1}^r\theta_j(\wdV)=0$ 
iff $\det C(\lambda)\not=0$.
Then, for $f_{\lambda}$ being the natural map $\Ug\to 
\End(\wdV(\lambda))$, one has
\begin{equation}
\label{imgf2}
\dim\Hom_{\fg}(\wdV, f_{\lambda}(\Soc\Ug))=\max \{\rk C(\lambda)|
C\in \pprv(\wdV)\}.
\end{equation}

Combining~\ref{dimhom} and~(\ref{imgf2}) one obtains

\subsubsection{}
\begin{cor}{corJPRV2}
Assume that $\wdM(\lambda)$ is simple. Then the image of the socle 
of $\Ug$ under the natural map $\Ug\to F(\wdM(\lambda),\wdM(\lambda))$ 
coincides with the socle of $F(\wdM(\lambda),\wdM(\lambda))$
iff $\eprv(\wdV)(\lambda)\not=0$ for any $\wdV\in\Irr$.
\end{cor}

\subsubsection{Change of Borel}
\label{PRVchangeB}
In the definition of PRV matrices we use
the Harish-Chandra projection. Therefore this definition
depends on the choice of triangular decomposition.
In the sequel we will add a lower index to
designate the corresponding Borel subsuperalgebra in
the cases when the choice of triangular decomposition
is not clear from the context.
For instance, $\pprv_{{\frak b}}(\wdV)$ is the set
of PRV matrices of the form
$\bigl(\cP_{{\frak b}}(\theta_j(v_i))\bigr)$
where $\theta_1,\ldots,\theta_r$ is a $\Zg$-basic system
of $\Hom_{\fg}(\wdV,\Ug)$ and $v_1,\ldots,v_r$ is
a basis of $\wdV|_0$.

Let ${\frak b}$ and ${\frak b}'$ be connected by the odd reflection 
along an odd isotropic root $\beta$. By~\ref{oddrefl},
$\wdV_{{\frak b}}(\lambda)=\wdV_{{\frak b}'}(\lambda')$ 
where $\lambda=\lambda'$
if $\lambda$ is such that $(\lambda,\beta)=0$ and $\lambda=\lambda'+\beta$
otherwise. Then for any $\wdV\in \Irr$
the equality~(\ref{imgf2}) implies
$$\begin{array}{c} 
\max \{\rk C(\lambda)|\
C\in \pprv_{{\frak b}_1}(\wdV)\}=
\max \{\rk C(\lambda')|\
C\in \pprv_{{\frak b}_2}(\wdV)\},\\
\eprv_{{\frak b}_1}(\wdV)(\lambda)=0\ \ 
\Longleftrightarrow\ \ 
\eprv_{{\frak b}_2}(\wdV)(\lambda')=0.
\end{array}$$

\subsubsection{}
\begin{rem}{remorder0}
The equality~\ref{imgf2})  has the following consequence: 
for any PRV determinant $p\in\eprv(\wdV)$ 
and any $\lambda\in {\frak h}^*$, the order of zero of 
$p$ at the point $\lambda$
is greater than or equal to the value
$(\dim\wdV|_0- \dim\Hom_{\fg}(\wdV,\End (\wdV(\lambda)))$.

In the case when $\fg$ is completely reducible, 
this property is used for the calculation
of PRV determinants---see~\cite{jl},\cite{jnato},\cite{gl1}.
Recall that, in this case, the PRV determinant 
is a polynomial defined up to a non-zero scalar (see~\Rem{origdef}).
\end{rem}

\subsection{Separation type theorem}
\label{ssectsprthm}
\subsubsection{}
\begin{defn}{genhom}
Call an $\ad\fg$-submodule $L$ of $\Ug$ {\em a generic harmonic space}
if the multiplication map provides an isomorphism
$L\otimes {\Zg}[S^{-1}]\iso {\Ug}[S^{-1}]$ for certain
$S\subset(\Zg\setminus\{0\})$.
\end{defn}

In~\Thm{fact?} we describe all generic harmonic spaces in ${\Ug}$.
Of course, if $L$ is a generic harmonic space 
then $L\otimes {\Zg}[S^{-1}]\iso {\Ug}[S^{-1}]$ for  
$S:=\Zg\setminus\{0\}$.
However it is always possible to find a smaller $S$;
in~\Cor{strS} we describe such $S$ in terms of PRV determinants.
As it is shown in Section~\ref{typeI}, for ${\frak g}={\frak {sl}}(m,n)$ or
${\frak {osp}}(2,2n)$ there is a generic harmonic space in ${\Ug}$
such that one can take $S=\{T^2\}$.

\subsubsection{}
\begin{lem}{LLL}
If an $\ad\fg$-submodule $L$ of $\Ug$ is such that
$$\begin{array}{ll}
(a) &  L\cong \displaystyle\oplus_{\wdV\in\Irr} \Inj(\wdV)^{\oplus\dim\wdV|_0},\\
(b) & \forall \wdV\in\Irr\ \   \ \eprv(\wdV;L)\not=0 \\ 
\end{array}$$
then it is a generic harmonic space.
\end{lem}
\begin{pf}
Set $\fh^*_f:=\{\nu|\ \dim\wdV(\nu)<\infty\}$.
For each $\nu\in \fh^*_f$  choose
a basis $\theta_1^{\nu},\ldots,
\theta_{r(\nu)}^{\nu}$ of $\Hom_{\fg}(\wdV(\nu),L)$.
Recall that $r(\nu)=\dim \wdV(\nu)|_0$. By~\ref{defprv},
the collection $\theta_1^{\nu},\ldots,\theta_{r(\nu)}^{\nu}$
forms a $\Zg$-basic system of $\Hom_{\fg}(\wdV(\nu),\Ug)$.
By~\ref{chuchu}, for a suitable $z(\nu)\in\Zg$
the collection 
$\theta_1^{\nu},\ldots,\theta_{r(\nu)}^{\nu}$ forms
a free $\Zg[z(\nu)^{-1}]$-basis of the localized module
$\Hom_{\fg}(\wdV(\nu),\Ug[z(\nu)^{-1}])$.

Set $S:=\{z(\nu),\ \nu\in\fh^*_f\}$.
Denote by $\psi$ the map $L\otimes \Zg [S^{-1}]\to\Ug [S^{-1}]$
induced by the multiplication.
For any $\nu\in\fh^*_f$ the collection 
$\theta_1^{\nu},\ldots,\theta_{r(\nu)}^{\nu}$ forms
a free $\Zg [S^{-1}]$ basis of the localized module
$\Hom_{\fg}(\wdV(\nu),\Ug [S^{-1}])$.
This means that the restriction 
of $\phi$ on the space
$\sum_{i=1}^{s(\nu)}\theta_i^{\nu}(\wdV)\otimes \Zg [S^{-1}]$
is a monomorphism and its image coincides with  the isotypical
component of $\wdV(\nu)$ in the socle of $\Ug [S^{-1}]$. Then
the restriction of $\phi$ on $\Soc L\otimes \Zg [S^{-1}]$ 
is a monomorphism and 
$\Soc \Ug [S^{-1}]=\phi(\Soc L\otimes  \Zg [S^{-1}])$.
Using the equality
$\Soc L\otimes  \Zg [S^{-1}]=\Soc(L\otimes  \Zg [S^{-1}])$
and~\ref{socphi},
one concludes that $\phi$ is a monomorphism.
From~\ref{suminj}, it follows
that  $L\otimes \Zg [S^{-1}]$ is an injective module in $\Fin$. 
Since the image of $\psi$ contains $\Soc \Ug [S^{-1}]$,
$\psi$ is a bijection. The lemma is proven.
\end{pf}

\subsubsection{}
\begin{cor}{HL}
If $H$ satisfies the assumption of~\Lem{hrmnew} then
$(\ad\Ug)(HT)$ is a generic harmonic space.
\end{cor}
\begin{pf}
By~\Lem{hrmnew}, $L:=(\ad\Ug)(HT)$ fulfills the condition (a)
of~\Lem{LLL}. 

Fix $\wdV\in\Irr$ and choose
a basis $\theta_1^,\ldots,
\theta_{r}$ of $\Hom_{\fg}(\wdV,L)$.
Observe that $r=\dim \wdV|_0$.
The injectivity of the map $L\otimes \Zg\to\Ug$, 
induced by the multiplication, implies the $\Zg$-linear independence
of $\theta_1,\ldots,
\theta_{r}$ in $\Hom_{\fg}(\wdV,\Ug)$.
Since the $\Zg$-rank of the latter is equal to $r$,
the collection $\theta_1,\ldots,\theta_{r}$
forms a $\Zg$-basic system of $\Hom_{\fg}(\wdV,\Ug)$.
Then, by~\ref{defprv}, $\eprv(\wdV;L)\not=0$.
By~\Lem{LLL}, $(\ad\Ug)(HT)$ is a generic harmonic space.
\end{pf}

\subsubsection{}
\begin{thm}{fact?}
An $\ad\fg$-submodule $L$ of $\Ug$ is a generic harmonic space
iff 
$$\begin{array}{ll}
(a) &  L\cong \displaystyle\oplus_{\wdV\in\Irr} \Inj(\wdV)^{\oplus\dim\wdV|_0}\\
\text{and} & \text{one of the following conditions holds}\\
(b) & \forall \wdV\in\Irr\ \   \ \eprv(\wdV;L)\not=0 \\ 
(c) & \text{ the multiplication map provides an embedding }
L\otimes \Zg\to\Ug.
\end{array}$$
\end{thm}
\begin{pf}
By~\Lem{LLL}, $L$ satisfying (a) and (b) is a generic harmonic space.
Arguing as in~\Cor{HL} one concludes that
$L$ satisfying (a) and (b) is also a generic harmonic space.

It remains to show that any generic harmonic space $L$
fulfills the conditions (a)-(c). The 
condition (c) obviously holds. Moreover, for
any $\wdV\in\Irr$ a basis $\theta_1,\ldots,\theta_r$
of $\Hom_{\fg}(\wdV,L)$ is a $\Zg$-basic system
of $\Hom_{\fg}(\wdV,\Ug)$. Therefore $\eprv(\wdV;L)\not=0$
and so (b) holds as well.

To verify the condition (a), let us show that all 
generic harmonic spaces are
pairwise isomorphic as $\ad\fg$-modules. Indeed, let
$L$ and $L'$  be generic harmonic spaces.
Since ${\Ug}$ is countably dimensional, one can choose
the corresponding sets $S,S'\subset{\Zg}\setminus\{0\}$ having
countable number of elements. Take a maximal ideal
$m$ of $\Zg$ such that $m\cap (S\cup S')=\emptyset$.
Then as $\ad\fg$-modules $L\cong \Ug/(m\Ug)\cong L'$.
In~\Cor{HL} we construct a generic harmonic space
satisfying (a).
Hence all generic harmonic spaces satisfy the condition (a).
The theorem is proven.
\end{pf}

\subsubsection{}
\label{receipy}
Now we can formulate the following ``receipt''. For each $\wdV\in \Irr$
fix a $\Zg$-basic system 
$\theta^{\wdV}_1,\ldots,\theta^{\wdV}_r$
of $\Hom_{\fg}(\wdV,{\Ug})$. The module
$\sum_{\wdV\in \Irr;j} T^2\theta^{\wdV}_j(\wdV)$
has an injective envelope inside an
injective $\ad\fg$-module ${\Ug}T$.
This injective envelope is a generic harmonic space by~\Thm{fact?}.

\subsubsection{}
Combining~\ref{chuchu} and the proof of~\Lem{LLL}, one obtains the

\subsubsection{}
\begin{cor}{strS}
Let $L$ be a generic harmonic space and $S\subset \Zg\setminus\{0\}$ 
satisfies the following property:
for any $\wdV\in\Irr$ there exists $s\in S$ 
such that $\cP(s)=t^2q$ where $q$ is a maximal 
$W.$-invariant divisor of $\bigl({\eprv}(\wdV;L)\bigr)^{|W|}$.
Then the multiplication map provides an isomorphism
$L\otimes {\Zg}[S^{-1}]\iso {\Ug}[S^{-1}]$.
\end{cor}

\section{Application to the description of  Verma module annihilators}
\label{sectSA}
In this section we prove~\Thm{cntrgn} which provides a connection
 between PRV determinants and the annihilators of simple modules.
\subsection{}
\begin{prop}{propJPRV1}
For a simple strongly typical Verma module$\wdM(\lambda)$
the following conditions are equivalent:

(i)  For any $\wdV\in\Irr$ one has $\eprv(\wdV)(\lambda)\not=0$.

(ii)The natural map $\Ug\to F(\wdM(\lambda),\wdM(\lambda))$ 
is surjective and 
$$F(\wdM(\lambda),\wdM(\lambda))\cong \displaystyle \oplus_{\wdV\in\Irr} 
\Inj(\wdV)^{\oplus\dim\wdV|_0}.$$
\end{prop}
\begin{pf}
Fix a strongly typical $\lambda$ such that  $\wdM(\lambda)$ is simple.
Denote by $f$ the natural map $\Ug\to F(\wdM(\lambda),\wdM(\lambda))$.

Let us show that (ii) implies (i).
Indeed, if both conditions of (ii) hold then 
$F(\wdM(\lambda),\wdM(\lambda))$ is projective in $\Fin$ and so
$f$ has a left inverse $f^{-1}$.
Denote by $N$ the image of $f^{-1}$. Take $\wdV\in\Irr$
and choose a basis $\theta_1,\ldots,\theta_r$ of $\Hom_{\fg}(\wdV,N)$.
By~\Rem{frfr}, $\theta_1,\ldots,\theta_r$ is a $\Zg$-basic
system of $\Hom_{\fg}(\wdV,{\Ug})$. Denote by $C$ the corresponding
PRV matrix. Since the restriction
of $f$ to $N$ is a monomorphism $\corank C(\lambda)=0$,
by~\ref{corankprv}. Hence $\det C(\lambda)\not=0$ and so
$\eprv(\wdV)(\lambda)\not=0$ as required.

Let us show that (i) implies (ii). Suppose that (i) holds.
By~\Cor{corJPRV2}, $f(\Soc \Ug)$ coincides with the socle of 
$F(\wdM(\lambda),\wdM(\lambda))$. The socle of $\Ug$ is a completely
reducible $\fg$-module and so it contains a submodule $L$ such
that the restriction of $f$
gives an isomorphism $L\iso \Soc F(\wdM(\lambda),\wdM(\lambda))$.
The central element $T^2$ acts  on $\wdM(\lambda)$
by a non-zero scalar $t(\lambda)^2$ (because $\lambda$
is strongly typical) and thus
one can choose $L$ lying in $T^2{\Ug}$.
By~\ref{TUg}, $T{\Ug}$ is an injective $\ad\fg$-submodule
of $\Ug$ and so it contains an injective envelope $\Inj(L)$ of $L$.
The restriction of $f$ to $\Inj(L)$ is a monomorphism because
the restriction of $f$ to $L$ is a monomorphism. 
Therefore $f(\Inj(L))$ is an injective module containing
the socle of $F(\wdM(\lambda),\wdM(\lambda))$. 
Using~\ref{bgrsocle}, one concludes that
the restriction of $f$ provides an isomorphism 
$\Inj(L)\iso F(\wdM(\lambda),\wdM(\lambda))$. 
In particular, $f$ is surjective and  
the $\ad\fg$-module $F(\wdM(\lambda),\wdM(\lambda))$
is injective in $\Fin$. Combining these
facts and~\ref{dimhom}, one obtains (ii).
\end{pf}

\subsection{}
\begin{thm}{cntrgn}
Suppose that $\lambda\in\fh^*$ is strongly typical and
$\eprv(\wdV)(\lambda)\not=0$ for all $\wdV\in\Irr$. Then
$\Ann \wdV(\lambda)=\Ug \Ann_{\Zg} \wdV(\lambda)$.
\end{thm}
\begin{pf}
For each  $\wdV(\nu)\in\Irr$ choose $p_{\nu}\in \eprv(\nu)$
such that $p_{\nu}(\lambda)\not=0$. Let
$\theta_1^{\nu},\ldots,\theta_{r(\nu)}^{\nu}\in 
\Hom_{\fg}(\wdV(\nu),{\Ug})$ be a $\Zg$-basic system
such that the determinant of the PRV matrix corresponding
to this system is equal to $p_{\nu}$.
Let $q_{\nu}$ be the maximal $W.$-invariant divisor of $p_{\nu}^{|W|}$
and $z_{\nu}\in\Zg$ be such that $\cP(z_{\nu})=t^2q_{\nu}$. 
Note that  $z_{\nu}\not\in \Ann_{\Zg} \wdV(\lambda)$  since 
$(t^2p_{\nu})(\lambda)\not=0$. 

Let $S$ be the multiplicative closure 
of the set $\{ T^2;z_{\nu}|\ \wdV(\nu)\in \Irr\}$ .
Denote by $A$ the localization of $\Zg$ by $S$ and by ${\Ug}_A$ 
the localization of $\Ug$ by $S$. 
The action of $\Ug$ on $\wdV(\lambda)$ can be canonically
extended to the action of
the localized algebra ${\Ug}_A$ since
$S\cap \Ann_{\Zg} \wdV(\lambda)=\emptyset$ .
The action of $\Ug$ on $\wdV(\lambda)$ can be canonically
extended to the action of
the localized algebra ${\Ug}_A$ acts on $\wdV(\lambda)$.
Clearly, the ideal $\Ann_A \wdV(\lambda)$ is maximal in $A$. 

Let us show that 
$\Ann_{{\Ug}_A} \wdV(\lambda)={\Ug}_A \Ann_A \wdV(\lambda)$.
Combining~\ref{receipy} and~\Cor{strS}, one concludes the existence
of $\ad\fg$-submodule $H$ of ${\Ug}$ such that
\begin{equation}
\label{socH}
\Soc H=\sum_{\nu: \wdV(\nu)\in\Irr} 
\sum_{i=1}^{r(\nu)}T^2\im\theta^i_{\nu}
\end{equation}
and that the multiplication map induces an isomorphism 
$H\otimes A\iso {\Ug}_A$. To verify
$\Ann_{{\Ug}_A} \wdV(\lambda)={\Ug}_A \Ann_A \wdV(\lambda)$
it is enough to check that $H\cap \Ann \wdV(\lambda)=0$.
Observe that $\eprv(\wdV(\nu);H)=p_{\nu}t^{2r(\nu)}$
and so $\eprv(\wdV(\nu),H)(\lambda)\not=0$ for all $\wdV(\nu)\in\Irr$.
Then, by~\ref{corankprv},
$$(\sum_{i=1}^{r(\nu)}\im\theta^i_{\nu})\cap\Ann \wdV(\lambda)=0.$$
that is $\Soc H\cap\Ann \wdV(\lambda)=0$. Therefore
$H\cap \Ann \wdV(\lambda)=0$ as required. Hence 
$\Ann_{{\Ug}_A} \wdM(\lambda)={\Ug}_A \Ann_A \wdM(\lambda)$.

Take $u\in \Ann_{{\Ug}} \wdV(\lambda)$. Write $u=\sum_{i=1}^m u_iz_i$
where $u_i\in H, z_i\in \Ann_A \wdV(\lambda)$. Choose $z\in S$
such that  $zu_i\in {\Ug}$ and $zz_i\in {\Zg}$ 
for all $i=1,\ldots, m$. Observe that 
$zz_i\in\Ann_{\Zg} \wdV(\lambda)$ for all $i=1,\ldots, m$.
Therefore 
$$z^2u=\sum_{i=1}^m (zu_i)(zz_i)\in {\Ug}\Ann_{\Zg} \wdV(\lambda).$$ 
Recall that $S\cap \Ann_{\Zg} \wdV(\lambda)=\emptyset$ and so 
$(z^2-c)\in \Ann_{\Zg} \wdV(\lambda)$ for a certain non-zero scalar $c$.
Hence $u\in {\Ug}\Ann_{\Zg} \wdV(\lambda)$. This completes the  proof.
\end{pf}

\subsection{PRV determinants and Shapovalov form factorization}
\label{shhh}
Assume that for some $\wdV\in \Irr$ the set
$\eprv(\wdV)$ contains a non-zero polynomial $p$
whose all irreducible factors are either a divisor of a Shapovalov
form (see~\ref{shapform}) or a divisor of $t$.
An irreducible factor of a Shapovalov form takes either a form
$(\alpha^{\vee}-(\rho,\alpha)-n)$ (a first type factor)
for some non-isotropic positive
root $\alpha$ and some $n\in {\Bbb N}^+$ or  
a form $(\alpha^{\vee}-(\rho,\alpha))$
(a second type factor) for some isotropic positive root 
$\alpha$. Note that
the factors of the second type are divisors of $t$. 
Observe that for $n\not=0$ the element
$$s_{\alpha}.(\alpha^{\vee}-(\rho,\alpha)-n)=
-\alpha^{\vee}+(\rho,\alpha)-n$$
is neither a divisor of a Shapovalov form nor a divisor of $t$.
This forces  the maximal $W.$-invariant factor of $p^{|W|}$  to be 
a divisor of $t^m$ for some $m\in {\Bbb N}$.

Then~\ref{fact?}--~\ref{receipy} imply 
\subsubsection{}
\begin{thm}{thmhope1}
Assume that for any $\wdV\in\Irr$ there exists a non-zero element in
$\eprv(\wdV)$ whose any irreducible factor is  
either a divisor of a Shapovalov
form or a divisor of $t$. Then there exists an 
$\ad\fg$-submodule $H$ of $\Ug$ such that
$$H\cong 
\displaystyle\oplus_{\wdV\in\Irr} \Inj(\wdV)^{\oplus \dim\wdV|_0}$$
and the multiplication map induces an isomorphism
$H\otimes\Zg[T^{-2}]\iso \Ug[T^{-2}]$.
\end{thm}

Recall that $\wdM(\lambda)$ is simple iff all 
Shapovalov forms are non-zero
at the point $\lambda$. \Thm{cntrgn} yields
\subsubsection{}
\begin{thm}{thmhope2}
Assume that for any $\wdV\in\Irr$ there exists a non-zero element in
$\eprv(\wdV)$ whose any irreducible factor is  
either a divisor of a Shapovalov
form or a divisor of $t$. Then the annihilator
of $\wdM(\lambda)$ is centrally generated provided 
$\wdM(\lambda)$ being simple strongly typical.
\end{thm}

As it is shown in~\cite{gl1} the condition of the last two theorems
holds for ${\frak g}={\frak {osp}}(1,2n)$.
We will prove that this condition also
holds for  the basic classical Lie superalgebras of type I.
Contrary to the case ${\frak {osp}}(1,2n)$, the proof is not
based on the calculations of PRV determinants.

\section{The basic classical Lie superalgebras of type I.}
\label{prepI}
Throughout this section ${\frak g}$ 
is a basic classical  Lie superalgebras of
type I that is ${\fg}={\frak {gl}}(m,n),{\frak {sl}}(m,n),
{\frak {psl}}(n,n)$ or 
${\fg}={\frak {osp}}(2,n)$. We describe some
common properties of these superalgebras which
are used in the next section. 

\subsection{Notation}
The Lie superalgebras of type I are  ${\Bbb Z}$-graded.
We denote by $\fg_r$ ($r\in {\Bbb Z}$) the corresponding
homogeneous component of $\fg$ and by $\fg_{\ol 0}$ (resp., $\fg_{\ol 1}$)
the even (resp., the odd) part of $\fg$.

One has $\fg_0=\fg_{\ol 0}$ and $\fg_{\ol 1}={\fg}_1+{\fg}_{-1}$.
In particular, both superalgebras $\fg_{\pm 1}$ are supercommutative
and  the exterior algebras
$\Lambda \fg_{\pm 1}$ are naturally embedded into $\Ug$. One has
$\Ug=\Lambda(\fg_{-1})\cU\Lambda(\fg_{1})$.
As $\ad\fg_{0}$-modules $\fg_{\pm 1}$ are irreducible and dual
one to another.

All odd roots are isotropic and so the notion of ``typical''
coincide with the notion of ``strongly typical''.

Extend the above ${\Bbb Z}$-grading on $U({\frak g})$; evidently
$U({\frak g})_{r}=0$ if 
$|r|>\dim\fg_1$. Moreover for $r=\dim\fg_1$
one has $U({\frak g})_{\pm r}=\cU\Lambda^r (\fg_{\pm 1})$.
By default, ``the degree'' of an element of 
$\Ug$ means its degree with respect to this ${\Bbb Z}$-grading.

\subsubsection{Distinguished triangular decomposition}
Set ${\frak n}_d^{{\pm}}:={\frak n}_0^{\pm}+{\frak g}_{\pm 1}$.
The decomposition ${\frak g}={\frak n}_d^-\oplus{\frak h}\oplus 
{\frak n}_d^+$ is called a {\em distinguished } triangular decomposition.
By default, all highest weight modules, positive/negative roots,
Harish-Chandra projection $\cP$ and PRV matrices/determinants
are constructed with respect to the distinguished triangular decomposition.

Note that $(\rho,\alpha)=(\rho_{0},\alpha)$ 
for any $\alpha\in\Delta_{0}$ since $\Delta^+_1=\Omega(\fg_1)$
is $W$-invariant. In particular, $w.\mu=w(\mu+\rho_0)-\rho_0$
for any $w\in W$ and $\mu\in\fh^*$.

The Verma modules constructed using a distinguished triangular
decomposition have the following nice structure:
$\wdM(\mu)=\Ind_{\fg_0\oplus\fn^+_1}^{\fg}M(\mu)$
where the action of $\fn^+_1$ on $M(\mu)$ is assumed to be trivial.

\subsubsection{}
\label{xI}
For $\alpha\in \Delta_{0}^+$ denote by $e_{\alpha}$ 
(resp., $f_{\alpha}$)
an element of the weight $\alpha$ (resp., $-{\alpha}$)
of ${\frak g}$. For $\beta\in \Delta_{1}^+$ 
(resp., $\beta\in \Delta_{1}^-$)
denote by $x_{\beta}$ (resp., $y_{\beta}$) an element of 
the weight $\beta$ of ${\frak g}$. 

Denote by $I$ the set of the positive odd roots with a fixed 
total order. Then $\{x_i\}_{i\in I}$ (resp.,
$\{y_i\}_{i\in I}$) is a basis of $\fg_{1}$ (resp., $\fg_{-1}$).
For any $J\subseteq I$ set  $x_J:=\prod_{i\in J} x_i$, 
$y_J:=\prod_{i\in J} y_i$ where the products are taken with 
respect to the total order on $I$. If we change the order of factors
in the product $y_J$ the result is equal either to  $y_J$ or 
to $(-y_J)$, since $\fg_{\pm 1}$ are supercommutative.
Note that $yy_I=xx_I=0$ for any $y\in\fg_{-1},x\in\fg_1$.

Evidently $x_I\in\Lambda^{top}\fg_1$, $y_I\in\Lambda^{top}\fg_{-1}$ 
and so $x_I,y_I$ are invariant with respect to the adjoint
action of $[\fg_0,\fg_0]$. Moreover $x_Iy_I$ is $\ad\fg_0$-invariant
due to the duality $\fg_1^*\cong\fg_{-1}$.

\subsubsection{}
\label{I+}
Sometimes we will deal with a non-distinguished
triangular decomposition 
${\frak g}={\frak n}^-\oplus{\frak h}\oplus {\frak n}^+$.
We shall use the following notation: ${\frak b}:=\fh+{\frak n}^+$;
$\cP_{\frak b}$ will be the Harish-Chandra projection
with respect to this triangular decomposition, 
$\Delta({\frak b})$ will be the set of non-zero roots of ${\frak b}$
and also 
$$\begin{array}{rl}
I_{+}:&=I\cap \Delta({{\frak b}}),\\
I_{-}:&=I\setminus I_{+},\\
t_{\frak b}:&=\cP_{\frak b}(T)
\end{array}$$
Note that
$$\begin{array}{rl}
t_{\frak b}:&=\prod_{\beta\in\Delta_1\cap\Delta_{\frak b}}
(\beta^{\vee}+(\beta,\rho)).
\end{array}$$

As before, the triangular decomposition of $\fg_0$ is assumed to be fixed
and all triangular decompositions of $\fg$ which we consider
are compatible with this triangular decomposition of $\fg_0$.

\subsubsection{Case ${\frak {psl}}(n,n)$}
\label{pslnn}
This case is rather special. A Cartan algebra of 
$\fg={\frak {psl}}(n+n\varepsilon)$ is ``too small''
and a distinguished triangular decomposition does not fit
the definition  given in~\ref{trdec}.
Moreover the restriction of the
Harish-Chandra projection determined by
the distinguished triangular decomposition
to the zero weight space is not an algebra
homomorphism: for instance, both $y_I,x_I$ have zero weight,
$\cP(y_I)=\cP(x_I)=0$ however $\cP(x_Iy_I)=t$---see~\Cor{PxIyI}.

A possible treatment to this problem is the following ``enlargement
of a Cartan subsuperalgebra''.
The Lie superalgebra $\fg$ is an ideal in  the Lie
superalgebra $\hat {\fg}:={\frak {gl}}(n,n)/({\Bbb C}X)$
where $X$ stands for the identity  matrix.
One has $\fg_{\ol 1}=\hat {\fg}_1$ and $\hat {\fg}_0=\fg_0\times {\Bbb C}z$
where $z$ is a central element of the Lie algebra $\hat {\fg}_0$.

Let ${\hat {\fh}}$ be a Cartan subsuperalgebra of $\hat {\fg}_0$
spanned by $\fh$ and $z$. One can easily sees that
${\hat {\fh}}$ acts semisimply on $\fg$
and a distinguished triangular decomposition of $\fg$
is determined, in a sense of~\ref{trdec},
by a certain regular element $h\in {\hat {\fh}}$.
For $\mu\in {\hat {\fh}}^*$ set 
$\Ug|_{\mu}:=\{u\in\Ug|\ [h,u]=\mu(h)u,\ \forall h\in{\hat {\fh}}\}$.
Then the restriction of the
Harish-Chandra projection determined by
the distinguished triangular decomposition
on $\Ug|_0$ is an algebra homomorphism because 
$\Ug|_0\cap \cU({\fn}^-+\fh){\fn}^-=0$.

It is easy to see that for any weight $\fg$-module $M$
one can extend (not uniquely) its $\fg$-module structure to a  
$\hat {\fg}$-module structure. This implies that~\Lem{2ass}
remains true for the distinguished triangular decomposition
if we define $\Ug|_0$ as above.
For $\ad\fg$-submodules of
$\Ug$ consider the natural $\hat {\fg}$-module structure
coming from the embedding $\Ug$ into $\cU(\hat {\fg})$.
In the sequel, we substitute the categories
$\Irr$ and $\wdO$ for $\fg$ by the same
categories for $\hat {\fg}$ (these  $\hat {\fg}$-categories have the
same sets of objects as their $\fg$-analogues 
but less morphisms). Under these conventions
all propositions of Sections~\ref{sectPRV},~\ref{sectSA}
remain true for $\fg={\frak {psl}}(n,n)$.

\subsection{Useful assertions}
The following lemma is an immediate consequence
of the supercommutativity of ${\frak g}_{\pm 1}$.

\subsubsection{}
\begin{lem}{g1free}
Let $\wdN$ be a $\fg$-module and $N$ be its 
${\frak g}_{0}$-submodule such that $\wdN=\Ug N$.
Then the canonical map $\Ind_{{\frak g}_{0}}^{\fg} N\to \wdN$
is an isomorphism iff for each collection $\{ v_J \}_{J\subseteq I}$ of
elements of $N$ one has 
$$x_I\sum_{J\subseteq I} y_J v_J=0
\ \Longrightarrow\ v_J=0,\ \forall J\subseteq I.$$ 
\end{lem}

Recall that $\wdM(\lambda)$ is typical
iff $t(\lambda)\not=0$ or, equivalently, $T\wdM(\lambda)\not=0$. .

\subsubsection{}
\begin{lem}{ssmod}
A typical Verma module
$\wdM(\lambda)$ contains a simple Verma submodule.
\end{lem}
\begin{pf}
Recall that $(\rho_0,\alpha)=(\rho,\alpha)$ 
for any $\alpha\in\Delta_0$.
Comparing the factorizations of Shapovalov
forms (see~\ref{shapform}) for the Lie superalgebras ${\frak g}$ 
and ${\frak g}_0$, one concludes that 
a typical $\fg$-module $\wdM(\lambda)$ is simple iff
the $\fg_0$-module $M(\lambda)$ is simple. Let $\wdM(\lambda)$
be a  typical Verma module and $v$ be its
highest weight vector.
Then $\cU(\fn^-_0)v\cong M(\lambda)$ contains a
simple submodule $M(\lambda')$.
Therefore there exists $u\in \cU(\fn^-_0)$ such that
${\frak n}_0^+uv=0$.
The $\ad {\frak g}_0$-invariance of ${\frak g}_1$ 
yields ${\frak g}_1(uv)=u\fg_1 v=0$. Thus the vector
$uv$ is primitive and so
$\Ug(uv)$ is a quotient of $\wdM(\lambda')$. The non-zero elements
of $\cU({\frak g}_{0})$ are non-zero divisors in $\Ug$
since the non-zero elements
of $\cS({\frak g}_{0})$ are non-zero divisors in 
$\Sg=\cS({\frak g}_{0})\Lambda(\fg_{\ol 1})$. Thus
$\Ug(uv)\cong\wdM(\lambda')$. Since $M(\lambda')$ is simple,
$\wdM(\lambda')$ is also simple. The lemma is proven.
\end{pf}

\subsection{Element $T$}
For the classical Lie superalgebras of type I the element $T$ 
takes a very simple form given by the following

\subsubsection{}
\begin{lem}{lPxIyI}
Up to a non-zero scalar $T=(\ad' y_Ix_I)1$.
\end{lem}
\begin{pf}
Recall~\ref{TUg}. It is enough to show that 
if $v$ is a generator of a trivial ${\frak g}_{0}$-module
$V(0)$ then $y_Ix_Iv$ spans a trivial ${\frak g}$-submodule of 
$\Ind_{{\frak g}_{0}}^{\frak g}V(0)$. In other words, one
has to verify that $\fg(y_Ix_I)\subset\Ug\fg_0$.
The $\ad {\frak g}_{0}$-invariance of $y_Ix_I$ forces 
${\frak g}_{0}(y_Ix_I)=(y_Ix_I)\fg_0$. Moreover
$\fg_{-1}(y_Ix_I)=0$. Thus it remains to check
that $x_{\beta}y_Ix_I\in\Ug\fg_0$ 
for any $\beta\in\Delta_{1}^+$. 
Setting $J:=I\setminus\{\beta\}$ one has
$$x_{\beta}y_Ix_I=
\pm x_{\beta}y_{-\beta}y_Jx_I=\pm ([x_{\beta},y_{-\beta}]
-y_{-\beta}x_{\beta})y_Jx_I
=\pm([x_{\beta},y_{-\beta}]y_Jx_I-y_{-\beta}((\ad x_{\beta}),y_J)x_I)$$
since $x_{\beta}x_I=0$.

The term $y_Ix_I$ is of the zero weight and so the weight of
$y_Jx_I$ is equal to $\beta$. The term $[x_{\beta},y_{-\beta}]$
lies in $\fh$; one has $[[x_{\beta},y_{-\beta}],x_{\beta}]=0$
since $x_{\beta}^2=0$. Therefore  $\beta([x_{\beta},y_{-\beta}])=0$
and so 
$[x_{\beta},y_{-\beta}]y_Jx_I=y_Jx_I[x_{\beta},y_{-\beta}]\in \Ug\fg_0.$

Let us show that $y_{-\beta}((\ad x_{\beta})y_J)x_I)\in\Ug\fg_0$.
For any $\beta'\in J$ one has $(\ad x_{\beta})y_{-\beta'}\in 
({\fn}^-_0+{\fn}^+_0)$. Since $\Lambda^r({\fg}_{-1})$ is $\ad\fg_0$-invariant
for all $r$, this implies that 
\begin{equation}\label{uu}
(\ad x_{\beta})y_J\in \Lambda^{l-2}({\fg}_{-1})
({\Bbb C}+{\fn}^-_0+{\fn}^+_0)\end{equation}
where $l:=\dim\fg_1$. 
If $y_{J'}\in\Lambda^{l-2}(\fg_{-1})$ has the same weight as
the element $(\ad x_{\beta})y_J$, then 
the set $I\setminus J'$ contains two elements
whose sum is equal to $2\beta$. Since the multiset $I$ contains
exactly one element equal to $\beta$, this implies $\beta\in J'$
that is $y_{-\beta}y_{J'}=0$. Then~(\ref{uu}) gives
$$y_{-\beta}((\ad x_{\beta})y_J)\in y_{-\beta}\Lambda^{l-2}({\fg}_{-1})
({\fn}^-_0+{\fn}^+_0)$$
and so 
$$y_{-\beta}((\ad x_{\beta})y_J)x_I\in \Lambda^{l-1}({\fg}_{-1})
({\fn}^-_0+{\fn}^+_0)x_I=
\Lambda^{l-1}({\fg}_{-1})x_I({\fn}^-_0+{\fn}^+_0)\in\Ug\fg_0$$
since $x_I$ is $\ad([\fg_0,\fg_0]$-invariant.
The assertion follows.
\end{pf}

Take any triangular decomposition $\fg=\fn^-\oplus\fn\oplus\fn^+$ and
retain notation of~\ref{I+}.

\subsubsection{}
\begin{lem}{PbT}
Up to a non-zero scalar one has
$$\cP_{{\frak b}}(x_{I_+}y_{I_+})
\cP_{{\frak b}}(y_{I_-}x_{I_-})=
\prod_{\beta\in\Delta_{1}({\frak b})}
({\beta}^{\vee}+{\beta}^{\vee}(\rho_{{\frak b}})).$$
\end{lem}
\begin{pf}
The right-hand side is equal to 
$t_{\frak b}=\cP_{{\frak b}}(T)$. 
By~\Lem{lPxIyI}, up to a non-zero scalar
$T=(\ad' x_Iy_I)1$. 
Thus one should verify that
$$\cP_{\frak b}((\ad' y_Ix_I)1)=\cP_{{\frak b}}(x_{I_+}y_{I_+})
\cP_{{\frak b}}(y_{I_-}x_{I_-})$$
up to a non-zero scalar. 

Since $\fg_{1}$ is supercommutative,  
$(\ad'x_I)1=2^{\dim\fg_1} x_I$ and 
$$\cP_{\frak b}((\ad' y_Ix_I)1)=2^{\dim\fg_1}\cP_{\frak b}((\ad' y_I) x_I)
=\pm 2^{\dim\fg_1}\cP_{\frak b}((\ad' y_{I_+}\ad'y_{I_-} )x_I).$$
Recall that for  $\alpha\in I_+$ one has $y_{-\alpha}\in {\frak n}^-$ 
and so $\cP((ad' y_{-\alpha})u)=\pm \cP(uy_{-\alpha})$ for
all $u\in {\Ug}$. Therefore, up to a non-zero scalar, 
$$\cP_{\frak b}((\ad' y_Ix_I)1)=\cP_{\frak b}(((\ad'y_{I_-} )x_I)y_{I_+}).$$
Similarly for any $\alpha\in I_-$ one has 
$y_{-\alpha}\in {\frak n}^+$ and so
$\cP_{\frak b}(uy_{-\alpha}y_{I_+})
=\cP_{\frak b}(\pm uy_{I_+}y_{-\alpha})=0$. Therefore
$$\cP_{\frak b}(((\ad'y_{I_-} )x_I)y_{I_+})=
\cP_{\frak b}(y_{I_-}x_Iy_{I_+})=\pm\cP_{\frak b}(y_{I_-}
x_{I_-}x_{I_+}y_{I_+})=\pm\cP_{\frak b}(y_{I_-}x_{I_-})
\cP_{\frak b}(x_{I_+}y_{I_+})$$
where the last equality is a consequence of the fact that
the restriction of the Harish-Chandra projection to 
the zero weight space is an algebra homomorphism.
Hence 
$$\cP_{\frak b}((\ad' y_Ix_I)1)=\pm 2^{\dim\fg_1}
\cP_{{\frak b}}(x_{I_+}y_{I_+})\cP_{{\frak b}}(y_{I_-}x_{I_-}).$$
The lemma is proven.
\end{pf}

\subsubsection{}
\begin{cor}{PxIyI}
Up to a non-zero scalar one has 
$\cP(x_Iy_I)=t=\prod_{\beta\in\Delta_{1}^+}
({\beta}^{\vee}+{\beta}^{\vee}(\rho))$.
\end{cor}

\subsubsection{}
\begin{lem}{atypVerma}
If $\lambda\in {\frak h}^*$ is such that
$\cP_{{\frak b}}(x_{I_+}y_{I_+})(\lambda)=0$ then 
$x_I\wdM_{{\frak b}}(\lambda)=0$. If 
$\cP_{{\frak b}}(y_{I_-}x_{I_-})(\lambda)=0$ then
$y_I\wdM_{{\frak b}}(\lambda)=0$.
\end{lem}
\begin{pf}
Both assertions are similar. To verify the first one,
fix $\lambda$ such that $\cP_{{\frak b}}(x_{I_+}y_{I_+})(\lambda)=0$ 
and set $\wdM:=\wdM_{{\frak b}}(\lambda)$.
Denote by $v$ a highest weight vector of $\wdM$.
One has
$$x_I\wdM=x_I \cU({\frak n}^-_0)
\displaystyle\sum_{J_1\subseteq I_+,J_2\subseteq I_-} x_{J_2}y_{J_1}v=
\cU({\frak n}^-_0)\displaystyle\sum_{J_1\subseteq I_+}x_Iy_{J_1}v$$
since $x_I$ is $\ad{\frak n}^-_0$-invariant and $x_Ix_J=0$ for
$J\not=\emptyset$.
Note that $\Ug_k v=0$ if $k>\# I_-$. Therefore $x_Iy_{J_1}v=0$
if $(\# I-\# J_1)>\# I_-$ that is if $\# J_1<\# I_+$.
Consequently 
$$x_I\wdM=\cU({\frak n}^-_0)x_Iy_{I_+}v=\cU({\frak n}^-_0)x_{I_-}
 x_{I_+}y_{I_+}v=\cU({\frak n}^-_0)x_{I_-}
\cP_{{\frak b}}(x_{I_+}y_{I_+})(\lambda)v=0$$
as required.
\end{pf}

\section{Separation and Annihilation theorems for 
type I case.}
\label{typeI}
Throughout this section ${\frak g}$ is a basic classical  
Lie superalgebras of type I that is 
${\fg}={\frak {gl}}(m,n),{\frak {sl}}(m,n){\frak {psl}}(n,n)$ or 
${\fg}={\frak {osp}}(2,n)$. 
In this section we establish separation and annihilation type
theorems for these algebras.

Retain notation of Section~\ref{prepI}. Until~\ref{outline} 
we will deal only with a distinguished triangular decomposition.
\subsection{}
\begin{prop}{FM}
For a typical Verma module $\wdM$ one has
$$\begin{array}{ll}
\text{(i) }& F(\wdM,\wdM)\cong 
\oplus_{V\in\Irr_0} 
(\Ind_{\fg_{0}}^{\fg} V)^{\oplus\dim V|_0}.\\
\text{(ii)} &
\text{For } H \text{ being a harmonic space of }
\cU({\frak g}_{0})
\text{ the restriction of the natural map }\\
 &  {\Ug}\to F(\wdM,\wdM)
\text{ to } (\ad {\Ug}) (Hx_Iy_I) \text{ is a bijection.}
\end{array}$$
\end{prop}
\begin{pf}
Let $H$ be a harmonic space of $\cU({\frak g}_{0})$.  
Recall that $\dim_{\fg_0}(V,H)=\dim V|_0$ for any $V\in\Irr_0$.

Define the map $\psi:\Ind_{\fg_{0}}^{\fg}(H)\to (\ad {\Ug})(Hx_Iy_I)$ 
by $g\otimes h\mapsto (\ad g)(hx_Iy_I)$. Since $x_Iy_I$ is 
$\ad{\fg}_{0}$-invariant, $Hx_Iy_I\cong H$ as $\ad\fg_0$-modules
and so $\psi$ is a $\fg$-epimorphism. Denote by $f$
the natural map $f: {\Ug}\to F(\wdM,\wdM)$.
To prove the theorem it is enough to show that
the map $f\circ\psi: \Ind_{\fg_{0}}^{\fg}(Hx_Iy_I)\to F(\wdM,\wdM)$
is a bijection.

Let $\mu$ be the highest weight
of $\wdM$ and $v$ be a highest weight vector of $\wdM$. Then
$M:=\cU({\frak g}_{0})v\cong M(\mu)$.
By~\ref{duflothm}, the restriction of the natural map $\cU({\frak g}_{0})\to
F(M,M)$ to $H$ is bijective. 

To prove that $f\circ\psi$ is a monomorphism
recall~\Lem{g1free}. Let  $\{ a_J \}_{J\subseteq I} \subset H$
be such that 
$$b\wdM=0\ \ \text{ where } \ 
b:=\bigl(\sum_{J\subseteq I} (\ad x_Iy_J) (a_Jx_Iy_I)\bigr).
$$
We need to show that $a_J=0$ for all $J\subseteq I$.
One has $(\ad y_J) (a_Jx_Iy_I)=y_Ja_Jx_Iy_I$ since
$\fg_{-1}y_I=0$. Therefore
$$b=(\ad x_I)(\sum_{J\subseteq I} y_J a_Jx_Iy_I).$$
Take $u\in\cU({\frak n}_{0}^-)$.
The equality $\fg_{-1}y_I=0$ implies
$$0=b(y_Iuv)=
\bigl((\ad x_I)(\sum_{J\subseteq I} y_J a_Jx_Iy_I)\bigr)y_Iuv=
\sum_{J\subseteq I} \pm y_J a_Jx_Iy_Ix_Iy_Iuv$$
since $y_Ix_{J'}y_I=0$ for $J'\not=I$ (the degree of $y_Ix_{J'}y_I$
is equal to $r:=\#J'-\# 2I$ and $\Ug_r=0$ for $r<-\#I$).
By~\Cor{PxIyI}, $\cP(x_Iy_I)=t$ and so
$x_Iy_Iuv=ux_Iy_Iv=t(\mu)uv$.
Therefore $b(y_Iuv)=t(\mu)^2\sum_{J\subseteq I} \pm y_J a_J uv$.
Hence $b(y_IM)=t(\mu)^2\sum_{J\subseteq I} \pm y_J a_J M$.
Since $t(\mu)\not=0$ and  $a_JM\subset M$, the equality $b(y_IM)=0$
implies $a_J M=0$ for all $J\subseteq I$. However the $a_J$ belong
to $H$ and $H\cap\Ann M=0$ by~\Thm{duflothm}. Therefore all $a_J$
are equal to zero. 
Hence $f\circ\psi$ is a monomorphism.

Let us show that $f\circ\psi$ is surjective. 
Recall that as $\fg_{0}$-module
$\wdM\cong M\otimes \Lambda {\frak g}_{-1}$.
Take a simple finite dimensional $\fg_{0}$-module $V$.
Using Frobenius reciprocity and~\ref{duflothm},
one obtains 
$$\begin{array}{ll}
\dim\Hom_{\fg_{0}}(V, F(\wdM,\wdM))&=
\dim\Hom_{\fg_{0}}(\wdM\otimes V,\wdM)\\
& =\dim\Hom_{\fg_{0}}(V\otimes\Lambda{\fg}_{-1}\otimes M,
\Lambda\fg_{-1}\otimes M)\\
& =\dim\Hom_{\fg_{0}}\bigl(V\otimes\Lambda\fg_{\ol 1},F(M, M)\bigr)\\
& =\dim\Hom_{\fg_{0}}\bigl(V\otimes\Lambda\fg_{\ol 1},H\bigr)
=\dim\Hom_{\fg_{0}}\bigl(V,\Lambda\fg_{\ol 1}\otimes H\bigr)\\
& =\dim\Hom_{\fg_{0}}(V, \Ind_{\fg_{0}}^{\fg} H)
=[\im (f\circ\psi) :V].
\end{array}$$ 
Since $F(\wdM,\wdM)$ 
is a completely reducible $\ad\fg_{0}$-module and the multiplicity
$[\im (f\circ\psi) :V]$ is finite, this gives
$F(\wdM,\wdM)=\im (f\circ\psi)$ and completes the proof. 
\end{pf}

Retain notation of~\ref{PRVVN}.
\subsubsection{}
\begin{cor}{H'}
Let $H$ be a harmonic space of $\cU({\frak g}_0)$
and $L:=(\ad\Ug)(Hx_Iy_I)$. Then 
$$L\cong\oplus_ {\wdV\in\Irr}\Inj(\wdV)^{\oplus \dim\wdV|_0}$$
and for any $\wdV\in\Irr$ the determinant
$\eprv(\wdV;L)$ admits a linear factorization.
Moreover each linear factor is a factor of a Shapovalov form.
\end{cor}
\begin{pf}
Combining~\Prop{FM} and~(\ref{indH})
one obtains 
$$L\cong\oplus_ {\wdV\in\Irr_f^0}
\Inj(\wdV)^{\oplus\dim\wdV|_0}.$$

Fix $\wdV\in\Irr$ and set $p:=PRV(\wdV;L)$. 
Recall that $\wdM(\mu)$ is simple iff $\mu$ is not a zero
of a Shapovalov form. 
In particular, an atypical Verma module is not simple.
Taking into account~\Cor{cor2ass}, 
one concludes from~\Prop{FM} that $p(\mu)\not=0$ 
provided that $\wdM(\mu)$ is simple.
Thus any zero of $p$
is a zero of a Shapovalov form. Since each Shapovalov
form admits a factorization into linear factors,
this implies that the set of zeroes of $p$ lies in a union
of countably many hyperplanes which correspond to the linear
factors of the Shapovalov forms. Therefore $p$ is a product
of linear factors which are factors of Shapovalov forms.
\end{pf}

\subsection{Separation theorem}
Combining~\Cor{H'},~\Thm{fact?},~\Cor{strS} and~\ref{shhh} one obtains
the following version of Separation Theorem
\subsubsection{}
\begin{thm}{spr1}
Let $H$ be a harmonic space of $\cU({\frak g}_ 0)$. Then 
the multiplication map provides an isomorphism
$\bigl((\ad {\Ug})Hx_Iy_I\bigr)\otimes\Zg[T^{-2}]\iso \Ug[T^{-2}]$.
\end{thm}

\subsection{Annihilation theorem}
In this subsection we prove that
$\Ann\wdM_{\frak b}(\lambda)$ is a centrally generated ideal iff 
$\wdM_{\frak b}(\lambda)$ is typical.

\subsubsection{}
\label{randy}
Combining~\Cor{H'} and~\Thm{thmhope2}, one concludes that
$\Ann\wdM(\lambda)$ is centrally generated if
$\wdM(\lambda)$ is simple. By~\Lem{ssmod}, any typical
Verma module $\wdM(\lambda)$ contains a simple Verma submodule.
This implies that $\Ann\wdM(\lambda)$ is centrally generated
if $\wdM(\lambda)$ is typical.
Using~\ref{oddrefl}, it is easy to
generalize this statement to any Borel subsuperalgebra ${\frak b}$.
Indeed, a typical Verma module $\wdM_{{\frak b}}(\lambda)$ 
is isomorphic to a Verma module $\wdM(\lambda')$. 
Recall that a Verma module $\wdM_{{\frak b}}(\lambda)$ 
is typical iff its annihilator does not contain $T$--- see~\ref{TUg}.
Thus $\wdM(\lambda')$ is typical and so
$\Ann\wdM_{{\frak b}}(\lambda)=\Ann\wdM(\lambda')$ 
is a centrally generated ideal. Hence 
$\Ann\wdM_{{\frak b}}(\lambda)$ is centrally generated
provided that $\wdM_{{\frak b}}(\lambda)$ is typical.

In the rest of this section we will show that 
$\Ann\wdM_{{\frak b}}(\lambda)\not=\Ug\Ann_{\Zg}\wdM_{{\frak b}}(\lambda)$
if $\wdM_{{\frak b}}(\lambda)$ is atypical that is
$t_{\frak b}(\lambda)=0$.

\subsubsection{}
\label{outline}
The proof goes as follows. Take
$\lambda\in {\fh}^*$ such that $t_{\frak b}(\lambda)=0$.
Set ${\wdchi}:=\Ann_{\Zg}\wdM_{{\frak b}}(\lambda)$.
The idea is to find $\wdV\in\Irr$ satisfying the following conditions
$$\begin{array}{ll}
a) & \Hom_{\fg}(\wdV,{\Ug})\not=0,\\
b) & \wdV \text{ is typical },\\
c) & \forall \theta\in \Hom_{\fg}(\wdV,{\Ug}) \ \ 
 \theta(\wdV)\subset \Ann\wdM_{{\frak b}}(\lambda).
\end{array}$$
Assume that $\Ann\wdM_{{\frak b}}(\lambda)$ is a 
centrally generated ideal.
Then the above conditions imply the equality $\Hom_{\frak g}(\wdV,{\Ug})=
{\wdchi}\Hom_{\frak g}(\wdV,{\Ug})$ that yields 
$$\eprv({\frak b},\wdV)\subseteq\cP({\wdchi})\eprv({\frak b},\wdV).$$
Since $\eprv({\frak b},\wdV)$ is a non-zero subset of $\Sh$
and $\cP({\wdchi})$ is a  subset of $\Sh$ which does not
contain non-zero scalars, the last inclusion is impossible.
Thus $\Ann\wdM_{{\frak b}}(\lambda)$
is not a centrally generated ideal.

\subsubsection{}
Retain notation of~\ref{xI}. There exists $z\in\fh$
such that $\ad z$ acts by zero on $\fg_0$ and by identity
(resp., minus identity) on $\fg_1$ (resp., $\fg_{-1}$).
Recall that
$\rho_1={1\over 2}\sum_{\beta\in \Delta_1^+}\beta$ and 
so $2\rho_{1}(z)=\dim {\frak g}_1$.
\subsubsection{}
\begin{lem}{sprv}
There exists $\nu\in {\frak h}^*$ such that 
$$\begin{array}{ll}
a) & \Hom_{\fg}(\wdV(\nu),{\Ug})\not=0,\\
b) & \wdV(\nu)\ \text{ is typical},\\
c') & z(\nu)=\dim {\frak g}_1.
\end{array}$$
\end{lem}
\begin{pf}
Recall that the condition b) is equivalent to the
inequality $t(\nu)\not=0$ where 
$t=\prod_{\beta\in \Delta_{1}^+} (\beta^{\vee}+(\rho,\beta))$.
Observe that $z\not=\beta^{\vee}$ for any $\beta\in \Delta_{1}^+$, since
$\beta^{\vee}(\beta)=0\not=z(\beta)$.
Thus the restriction of the polynomial $t$ on the hyperplane
$$S:=\{\mu\in {\frak h}^*|\ z(\mu)=\dim {\frak g}_1,\ (\Zg\cap\fh)(\mu)=0\}$$ 
is a non-zero polynomial
(one has $\Zg\cap\fh=0$ apart from the case $\fg={\frak {gl}}(m,n)$).

Consider the set 
$$X:=\{\mu\in {\frak h}^*|\ 
\Hom_{{\fg}_0}(V(\mu), \cU({\frak g}_0))\not=0\}.$$
It is easy to see that $X$ is a Zariski dense 
subset of the hyperplane 
$$\{\mu\in {\frak h}^*|\ z(\mu)=0, (\Zg\cap\fh)(\mu)=0\}=S-2\rho_1.$$
Thus $t(\nu'+2\rho_{1})\not=0$ for a certain $\nu'\in X$.
Hence $\nu:=\nu'+2\rho_{1}$
fulfills the conditions b) and c'). Let show that $\nu$ satisfies a).
Since $\nu'\in X$, there exists a copy $V\cong V(\nu)$ inside
$\cU({\frak g}_0)$. Let $u$ be a highest
weight vector of $V$. Then $ux_I$ is a non-zero primitive
(that is $({\frak n}_0^++{\frak g}_1)$-invariant)
element of $\Ug$ having weight $\nu$. Therefore the $\ad\fg$-submodule
generated by $ux_I$ is a finite dimensional quotient of
$\wdM(\nu)$. Since $t(\nu)\not=0$, a unique finite
dimensional quotient of $\wdM(\nu)$ is isomorphic to $\wdV(\nu)$.
Hence $\Hom_{\fg}(\wdV(\nu),{\Ug})\not=0$.
\end{pf}
\subsubsection{}
\begin{rem}{lastrem}
Similarly there exists a typical simple module $\wdV$ with
the lowest weight $\nu$ such that
$\Hom_{\fg}(\wdV,{\Ug})\not=0$ and $z(\nu)=-\dim {\frak g}_1$.
\end{rem}

\subsubsection{}
Fix $\lambda\in {\frak h}^*$ such that
$\cP_{{\frak b}}(x_{I_+}y_{I_+})(\lambda)=0$.
Set ${\wdchi}:=\Ann_{\Zg}\wdM(\lambda)$ and
assume that $\Ann\wdM(\lambda)={\Ug}{\wdchi}$.

Fix $\nu$ satisfying the conditions a)---c') of~\Lem{sprv}. 
We claim that $\nu$ fulfills
the condition c) of~\ref{outline}.

Indeed, let $v$ be a highest weight vector of $\wdV:=\wdV(\nu)$.
Take any $\theta\in\Hom_{\frak g}(\wdV,{\Ug})$.
By the assumption c')
$zv=(\dim\fg_1)v$ that is $[z,\theta(v)]=(\dim\fg_1)\theta(v)$.
Therefore $\theta(v)$ has degree $\dim\fg_1$ that is
$\theta(v)\in \cU({\frak g}_0) x_I$. \Lem{atypVerma} asserts that
$x_I\wdM_{{\frak b}}(\lambda)=0$ and so
$\theta(v)\in\wdM_{{\frak b}}(\lambda)=0$.
Since $\Ann\wdM_{{\frak b}}(\lambda)$ is $\ad\fg$-stable,
it contains $\theta(\wdV)$.

By the assumption $\Ann\wdM(\lambda)={\Ug}{\wdchi}$ and
so  the element $\theta(v)$
can be written in a form $\theta(v)=\sum z_iu_i$ where all $z_i$ belong
to $\wdchi$. Without loss of generality we can assume that each 
element $u_i$ has the same weight and the same central character 
as $\theta(v)$ with respect to the adjoint action of $\Zg\subset {\Ug}$ on
$\Ug$. Since $\wdV$ is typical, $\wdV$ is a
unique, up to isomorphism, cyclic module in $\Fin$ with
this central character. Therefore $\ad{\fg}$-span of each $u_i$
is isomorphic to $\wdV$ that is $u_i\in \theta_i(\wdV)$
for certain $\theta_i\in\Hom_{\frak g}(\wdV,{\Ug})$.
Since $u_i$ has weight $\nu$ and $v$ spans $\wdV|_{\nu}$,
one can assume that $u_i=\theta_i(v)$.
Then $\theta(v)=\sum z_i\theta_i(v)$ and so $\theta=\sum z_i\theta_i$.
Hence
$$\Hom_{\frak g}(\wdV,{\Ug})={\wdchi}\Hom_{\frak g}(\wdV,{\Ug}).$$
Using the fact that 
$\cP_{\frak b}(z\theta)=\cP_{\frak b}(z)\cP_{\frak b}(\theta)$ for any 
$z\in{\Zg},\theta\in\Hom_{\fg}(\wdV(\nu),{\Ug})$, one concludes
$\eprv_{\frak b}(\wdV(\nu))\subseteq\cP({\wdchi})\eprv_{\frak b}(\wdV(\nu))$. 
Since $\eprv_{\frak b}(\wdV(\nu))$ is a non-zero subset
of a polynomial algebra $\Sh$ and $\cP({\wdchi})$ does not contain
non-zero scalars, the last inclusion is impossible. Hence
Hence $\Ann\wdM(\lambda)\not={\Ug}{\wdchi}$.

In the case $\cP_{{\frak b}}(y_{I_-}x_{I_-})(\lambda)=0$ the proof
is similar:  one may choose $\wdV$ as in~\Rem{lastrem}.

Hence we prove that for any triangular decomposition
the following theorem holds.

\subsubsection{}
\begin{thm}{anntypeI}
The annihilator of a Verma module  is a centrally
generated ideal iff this module is typical.
\end{thm}

\subsection{}
\label{EquivTypeI}
Denote by $\Max A$ the set of maximal ideals of an algebra $A$.
In~\cite{psI}, the following theorem is proven.

\subsubsection{}
\begin{thm}{ps}
For any $\wdchi\in\Max\Zg$ not containing $T^2$,
there exists $\chi\in \Max\cZ(\fg_0)$
such that the algebra $\Ug/(\Ug\chi)$ is the matrix algebra
over $\cU(\fg_0)/(\cU(\fg_0)\chi_0)$.
\end{thm}

One can deduce   from this theorem that the annihilator
of a typical Verma module is centrally generated.
On the other hand, one can deduce~\Thm{ps} 
from~\Prop{FM} and~\ref{randy}.
In fact, take $\lambda\in\fh^*$ such that $\wdchi \wdM(\lambda)=0$.
Then $t(\lambda)\not=0$ and so
$$\begin{array}{l}
\Ug/(\Ug\wdchi)=F(\wdM(\lambda),\wdM(\lambda))=
F(M(\lambda)\otimes \Lambda \fg_1^-,M(\lambda)\otimes \Lambda \fg_1^-)=\\
F(M(\lambda),M(\lambda))\otimes \End(\Lambda\fg_1^-)=
\cU(\fg_0)/(\cU(\fg_0)\chi)
\otimes \End(\Lambda\fg_1^-)\end{array}$$
where $\chi:=\Ann_\cZ(\fg_0) M(\lambda)$. Hence $\Ug/(\Ug\chi)$
is a matrix algebra over $\cU(\fg_0)/(\cU(\fg_0)\chi)$.

\section{Perfect pairs}
\label{apprpairs}
In this section we find  for each maximal ideal of $\Zg$ not
containing $T^2$ {\em a perfect mate} which is
a maximal ideal in $\cZ(\fg_0)$ possessing certain properties.

\subsection{}
Denote by $\Max A$ the set of maximal ideals of an algebra $A$
and by $A-\Mod$ the full category of left $A$-modules.
For an $A$-module $N$ and $\chi\in\Max A$ set
$$N_{\chi}:=\{v\in N|\ \chi^rv=0,\ \forall r>>0\}.$$

Call a maximal ideal of $\Zg$ {a $\fg$-central character}. 
For $\wdchi\in \Max \Zg$ (resp., $\chi\in \Max\cZ(\fg_0))$
set
$$\Ug_{\wdchi}:=\Ug/(\Ug\wdchi),\ \ \ \cU_{\chi}:=\cU/(\cU\chi).$$
We canonically identify the $\Ug_{\wdchi}$-modules
and the (non-graded) $\Ug$-modules annihilated by $\wdchi$.

\subsubsection{}
\begin{defn}{strtyp}
Call $\wdchi\in\Max\Zg$ {\em strongly typical} if $T^2\not\in\wdchi$.
\end{defn}

{\em Throughout this section} $\wdchi$ {\em denotes a strongly typical
central character}.

\subsubsection{}
For a fixed triangular decomposition set
$$W(\wdchi):=\{\lambda\in\fh^*|\ \wdchi\wdM(\lambda)=0\}.$$
Recall that $\cP(\Zg)\supset \cP(T^2)\Sh^{W.}$ and so
$W(\wdchi)$ forms a single $W.$-orbit.

Remark that  $\wdchi$ is strongly typical iff 
any $\lambda\in W(\wdchi)$ 
is strongly typical.

For any pair $(\wdchi,\chi)\in \Max\Zg\times\Max\cZ(\fg_0)$
there is a functor 
$\Psi_{\wdchi,\chi}:\Ug_{\wdchi}-\Mod\to \cU_{\chi}-\Mod$
given by $\wdN\mapsto\wdN_{\chi}$.
One of our goals is to prove that  any strongly typical
$\fg$-central character $\wdchi$ has a ``perfect mate''
$\chi\in\Max\cZ(\fg_0)$ such that the above functor provides
the equivalence of categories. For the basic classical
Lie superalgebras of type I this is proven in cite~\cite{psI}. 
In~\cite{geq}, we will prove it for the basic classical
Lie superalgebras of type II.

For  type I case, it turns out that for any $\lambda\in W(\wdchi)$
$\chi:=\Ann_{\cZ(\fg_0)} M(\lambda)$  
is a perfect mate for $\wdchi$
(provided that $\wdchi$ is strongly typical)--- see also~\ref{EquivTypeI}.
This does not hold for type II case.

It is easy to see that if the functor $\Psi_{\wdchi,\chi}$
provides the equivalence of the categories, then
$\chi$ has, at least, the following properties: for a 
Verma $\Ug_{\wdchi}$-module $\wdM$ the $\fg_0$-module
$\wdM_{\chi}$ is a Verma $\fg_0$-module and
for any $\Ug_{\wdchi}$-module $\wdN$ the $\fg_0$-module
$\wdN_{\chi}$ is non-zero. We call $\chi$ {\em a mate of}
$\wdchi$ if it possesses the first property
and {\em a perfect mate} if it possesses both properties.
As we shall show in~\cite{geq}, these two properties really
ensure the equivalence of categories
$\Ug_{\wdchi}-\Mod$ and $\cU_{\chi}-\Mod$ provided that $\wdchi$
is strongly typical.

A pair $(\wdchi,\chi)\in \Max\Zg\times\Max\cZ(\fg_0)$
is called {\em a perfect pair} if $\chi$ 
is a perfect mate  for $\wdchi$.

The goal of this section is to find  a perfect mate 
for any strongly typical $\wdchi$.
This is done in the following way.
\Lem{dist-appr} gives  a combinatorial criterion on $\chi$ to be 
a perfect mate for $\wdchi$.
In~\ref{rootsystems} we fix, for each 
basic classical Lie superalgebra of type II
a triangular decomposition, in terms 
of which we  will describe a perfect mate $\chi$
for $\wdchi$. In~\ref{genericcase}
we consider $\wdchi$ satisfying a certain genericity
condition. For such a ``generic'' $\wdchi$ we show
that for a suitable $\lambda\in W(\wdchi)$ 
the $\fg_0$-central character of $M(\lambda)$ 
is a perfect mate for $\wdchi$. The case when
$\wdchi$ is not generic is treated in~\ref{nongencase}.
It turns out that for $\fg$ of the types $B(m,n)$, $G(3)$
any strongly typical $\wdchi$ is ``generic''.
The remaining basic classical
Lie superalgebras of type II are treated
case by case.

\subsection{Notation}
\subsubsection{}
\label{depth}
We say that a $\Ug$-module $\wdN$ has {\em a finite support }
$\supp_{\Zg}\wdN=\{\wdchi_1,\ldots,\wdchi_k\}$ if
for any $v\in\wdN$ there exist $r_1,\ldots,r_k\in {\Bbb N}^+$
such that $\prod_i\wdchi_i^{r_i}v=0$. In this case,
$$\wdN=\oplus_i \wdN_{\wdchi_i}$$
and each $\wdN_{\wdchi_i}$ is canonically isomorphic to the localization
of the module $\wdN$ by the maximal ideal $\wdchi_i$.
If $\wdN$ has a finite support and $0\to\wdN'\to\wdN\to \wdN''\to 0$
is an exact sequence then, for any $\wdchi'\in\Max\Zg$, the sequence
$0\to\wdN'_{\wdchi'}\to\wdN_{\wdchi'}\to \wdN''_{\wdchi'}\to 0$
is also exact.

For a $\fg$-module $\wdN$ and a maximal ideal $\wdchi'\in\Max\Zg$
define  $\wdchi'$-{\em depth} of $\wdN$ to be 
a minimal $r\in {\Bbb N}\cup\infty$
such that $(\wdchi')^r \wdN_{\wdchi}=0$.

We adopt the similar notation for $\cZ(\fg_0)$ and $\cU$-modules.

\subsubsection{}
\label{Gamma}
Denote by $\Gamma$ the set of subsets of $\Delta^+_1$.
For $\gamma\subseteq \Delta^+_1$ set
$$|\gamma|:=\sum_{\beta\in \gamma}\beta.$$

Define the action of the Weyl group $W$ on $\Gamma$
by setting
$$w_*\gamma:=w\bigl(\gamma\cup -(\Delta^+_1\setminus\gamma)\bigr)
\cap \Delta^+_1.$$
Then
$$|w_*\gamma|=w(|\gamma|-\rho_1)+\rho_1.$$

\subsubsection{}
\label{Vermafilt}
As a $\fg_0$-module, a Verma module $\wdM(\lambda)$ has a filtration
$0=M_0\subset M_1\subset\ldots\subset M_r=\wdM(\lambda)$
such that the set of factors $M_{i+1}/M_i$ coincides
with the multiset $\{M(\lambda-\gamma):\gamma\in \Gamma\}$--- 
see~\cite{mu}, 3.2.

It is easy to check that for any $w\in W,\gamma\in\Gamma$
\begin{equation}
\label{w*}
w.\lambda-|w_*\gamma|+\rho_0=w(\lambda-|\gamma|+\rho_0).
\end{equation}
Therefore the $\fg_0$-central characters of
$M(w.\lambda-|w_*\gamma|)$ and $M(\lambda-|\gamma|)$
coincide. Thus for fixed $\wdchi$
the multiset of $\fg_0$-central characters 
of $\{M(\lambda-\gamma):\gamma\in \Gamma\}$ does not
depend on the choice of $\lambda\in W(\wdchi)$.

\subsubsection{}
\begin{defn}{dist}
Call $\chi\in\Max \cZ(\fg_0)$
{\em a mate} for $\wdchi\in \Max \Zg$ if
for some $\wdM(\lambda)$ satisfying $\wdchi\wdM(\lambda)=0$,
the $\fg_0$-module $\wdM(\lambda)_{\chi}$ is isomorphic to a
Verma module $M(\lambda-|\gamma|)$ for some
$\gamma\in\Gamma$.
\end{defn}

One can easily deduce from~\ref{Vermafilt} that
if $\chi$ is mate for $\wdchi$ then for {\em any}
$\wdM(\lambda)$ satisfying $\wdchi\wdM(\lambda)=0$,
the $\fg_0$-module $\wdM(\lambda)_{\chi}$ is isomorphic to a
Verma module $M(\lambda-|\gamma|)$ for some
$\gamma\in\Gamma$. In particular,
if $\chi$ is mate for $\wdchi$ then for
any $\wdM$ satisfying $\wdchi\wdM=0$, the
$\chi$-depth of $\wdM$ is equal to $1$.

\subsubsection{}
\begin{defn}{appr}
Let us call $\chi\in\Max \cZ(\fg_0)$
{\em a perfect mate} for $\wdchi\in \Max \Zg$
if it is a mate for $\wdchi$
and for any non-zero $\Ug_{\wdchi}$-module $\wdN$ one has $\wdN_{\chi}\not=0$.
\end{defn}

\Defn{dist} is given in terms of
the category $\wdO$. However,~\ref{oddrefl} implies that if
for some triangular decomposition $\chi$ is a mate for $\wdchi$,
then it is  a mate 
also with respect to all other choices of  triangular decomposition.

\subsection{}
Throughout this subsection a strongly typical $\fg$-central character $\wdchi$
is assumed to be fixed.
Our first goal is to reformulate~\Defn{appr}
in terms of category $\wdO$.

\subsubsection{}
\begin{lem}{lemfinsupp}
There exist $\chi_1,\ldots,\chi_k\in\Max\cZ(\fg_0)$
and $r_1,\ldots,r_k\in {\Bbb N}^+$ such that
for any $\Ug_{\wdchi}$-module $\wdN$ one has
$$\wdN=\displaystyle\oplus_{1}^k \wdN_{\chi_i}$$
and $\chi^{r_i}\wdN_{\chi_i}=0$.
\end{lem}
\begin{pf}
Theorem 2.5 of~\cite{mu} implies that for any
$x\in\cZ(\fg_0)$ there exist $z_0,\ldots,z_l\in\Zg$ (where $l=\#\Gamma$)
such that $\sum_0^{l} z_ix^i=0$ and $z_l=T^2$
(notice that 2.1 of~\cite{mu}
contains a misprint; to correct it one
has to substitute $g$ by $g^2$ in 2.1 and in Theorem 2.5). 
Therefore for any $x\in\cZ(\fg_0)$ there
exists $c_0,\ldots,c_{l-1}\in{\Bbb C}$ 
such that $x^l+\sum_0^{l-1} c_ix^i\in (\cZ(\fg_0)\cap\Ug\wdchi)$.
Consequently, the ideal 
$(\cZ(\fg_0)\cap\Ug\wdchi)$ has a finite codimension 
in $\cZ(\fg_0)$ and so there
exist $\chi_1,\ldots,\chi_k\in\Max\cZ(\fg_0)$
and $r_1,\ldots,r_k\in {\Bbb N}^+$ such that
$$(\cZ(\fg_0)\cap\Ug\wdchi)\supseteq \prod_{i=1}^k\chi_i^{r_i}.$$
Then, for any $\Ug_{\wdchi}$-module $\wdN$
$$\prod_{i=1}^k\chi_i^{r_i}\subseteq \Ann_{\cZ(\fg_0)}\wdN.$$
The assertion follows.
\end{pf}

\subsubsection{}
\begin{cor}{corgenby}
If $\chi$ is a perfect mate for $\wdchi$ then 
$$\wdN=\Ug\wdN_{\chi}$$
for any $\Ug_{\wdchi}$-module $\wdN$.
\end{cor}
\begin{pf}
Since $\supp_{\cZ(\fg_0)}\wdN$ is finite, one has
$(\wdN/(\Ug\wdN_{\chi}))_{\chi}=0$.
Hence $\wdN/(\Ug\wdN_{\chi})=0$ because 
$\chi$ is a perfect mate for $\wdchi$.
\end{pf}

\subsubsection{}
\begin{lem}{anncnt}
If $\chi\in\Max\cZ(\fg_0)$ is a mate  for 
$\wdchi\in \Max \Zg$ and
$\wdV(\lambda)_{\chi}\not=0$ for all $\lambda\in W(\wdchi)$, then
$\chi$ is a perfect mate.
\end{lem}
\begin{pf}
Suppose that $\wdV(\lambda)_{\chi}\not=0$
for all $\lambda\in W(\wdchi)$. One has to verify
that $\wdN_{\chi}\not=0$ for any non-zero $\Ug_{\wdchi}$-module $\wdN$.
Since $\Ug$ is Noetherian, $\wdN$ has a simple subquotient $\wdN'$.
\Lem{lemfinsupp} implies that $\wdN_{\chi}\not=0$ provided 
$\wdN'_{\chi}\not=0$.
Hence it is enough to check that $\wdN_{\chi}\not=0$ 
for any simple $\Ug_{\wdchi}$-module $\wdN$. 

Take a simple  $\Ug_{\wdchi}$-module $\wdN$. 
The ideal $\Ann \wdN$  is  primitive
and so, according to~\cite{mu2},  $\Ann \wdN=\Ann\wdV$
where $\wdV$ is a simple highest weight module. Obviously,
$\wdchi\wdV=0$  and so, by our assumption, $\wdV_{\chi}\not=0$. 
As a $\fg_0$-module, $\wdV$ has a finite length. Therefore 
$$\Ann_{\cZ(\fg_0)} \wdN=
\Ann_{\cZ(\fg_0)}\wdV=\chi^{r_1}\chi_2^{r_2}\ldots\chi_k^{r_k}$$
where 
$\{\chi,\chi_2,\ldots,\chi_k\}=\supp_{ \cZ(\fg_0)} \wdV,\ 
r_1,\ldots,r_k\in {\Bbb N}^+$.
Then $\wdN':=\chi_2^{r_2}\ldots\chi_k^{r_k}\wdN\not=0$ and 
$\chi^{r_1} \wdN'=0$. Hence $\wdN_{\chi}\not=0$.
The lemma is proven.
\end{pf}

\subsubsection{}
\begin{lem}{dist-appr}
Take $\lambda\in W({\wdchi})$, $\gamma\in\Gamma$
and set $\chi:=\Ann_{\cZ(\fg_0)} M(\lambda-|\gamma|)$.
$$\begin{array}{ll}
(i) & \forall \gamma'\in \Gamma\setminus\{\gamma\}\ \ \ 
(\lambda-|\gamma'|+\rho_0)\not\in W(\lambda-|\gamma|+\rho_0)\ \ 
\Longleftrightarrow \ \ 
\chi \text{ is a mate for } \wdchi\\
(ii) &\text{If }\chi \text{ is a mate for }\wdchi \text{ and }
\ \Stab_W(\lambda-|\gamma|+\rho_0)\subseteq \Stab_W(\lambda+\rho)\ 
\text{ then }\chi\text{ is a perfect mate}.
\end{array}$$

\end{lem}
\begin{pf}
The equivalence (i) follows from~\ref{Vermafilt}.

For (ii) recall~\Lem{anncnt}. Suppose that $\chi$ is a mate for $\wdchi$
but it is not a perfect mate. Then 
$\wdM(\lambda)_{\chi}=M(\lambda-|\gamma|)$ and
$\wdV(w.\lambda)_{\chi}=0$ for some $w\in W$. The equality~(\ref{w*}), 
implies that for any $y\in W$ 
$$\wdM(y.\lambda)_{\chi}=M(y.\lambda-|y_*\gamma|).$$
Therefore $\wdV(y.\lambda)_{\chi}$ is a quotient of 
$M(y.\lambda-|y_*\gamma|)$.

The module $\wdV(w.\lambda)$ is a homomorphic
image of $\wdM(w.\lambda)$.
Denote the kernel of this homomorphism by $\wdN$.
The module $\wdN$ has a finite length and the factors
of its Jordan-G\"older series have the form $\wdV(\mu)$
for some $\mu\in W.\lambda$ satisfying $\mu<w\lambda$.
Since $0=\wdV(w.\lambda)_{\chi}=(\wdM/\wdN))_{\chi}$
one concludes that the $\fg_0$-module 
$\wdN_{\chi}=\wdM_{\chi}=M(w.\lambda-|w_*\gamma|)$
has a finite filtration such that each  factor
of this filtration is a quotient of 
$M(y.\lambda-|y_*\gamma|)$ for some $y\in W$
satisfying $y.\lambda<w.\lambda$.
Therefore
$$w.\lambda-|w_*\gamma|=y.\lambda-|y_*\gamma|$$ 
for some $y\in W$ satisfying $y.\lambda<w.\lambda$. By~(\ref{w*}),
the above equality is equivalent to
the condition $y^{-1}w\in \Stab_W (\lambda+\rho_0-|\gamma|)$.
However, $y.\lambda<w.\lambda$ implies
$y^{-1}w\not \in \Stab_W (\lambda+\rho)$. 
Thus $\Stab_W (\lambda+\rho_0-|\gamma|)\not
\subseteq \Stab_W (\lambda+\rho)$ as required. 
\end{pf}

\subsection{}
Throughout this subsection a strongly typical $\fg$-central character
$\wdchi$ is assumed to be fixed.

\subsubsection{Right version}
\label{leftright}
One may repeat the above reasonings for the right $\fg$-modules.
Suppose now that $(\wdchi,\chi)$ is a perfect pair
in the ``left sense''. We claim that it is 
a perfect pair in the ``right sense''.

Indeed, the superantiautomorphism $\sigma$ (see~\ref{cntrgr})
provides the duality between the left and the right $\Ug$-modules
given by $N\mapsto N^{\sigma}$:
$$v^{\sigma}.a:=(-1)^{d(a)d(v)}\sigma(a).v,\ \ a\in\Ug,v\in N.$$
By~\ref{sigma;cent}, $\sigma$ stabilizes the elements
of $\Zg$; thus the restriction of the above duality
gives the duality between the left and the right 
$\Ug_{\wdchi}$-modules.
Since $\sigma$ also stabilizes the elements of $\cZ(\fg_0)$,
in order to show that $(\wdchi,\chi)$ is a perfect pair in 
the ``right sense'',
it is enough to check that $\chi$ is  mate for $\wdchi$
in the ``right sense''.

Take a minimal element $\lambda\in W(\wdchi)$. 
Then $\wdM:=\wdM(\lambda)$
is a simple Verma module and $\wdM^{\sigma}$ 
is a right simple module which is a Verma module
with respect to a suitable
triangular decomposition of $\fg$. Since $\chi$ is a mate for
$\wdchi$, $\wdM_{\chi}$ is a Verma $\fg_0$-module, say $M$.
Then 
$$_{\chi}\wdM^{\sigma}:= \{v\in \wdM^{\sigma}|\ v.\chi^r=0, \ \ r>>0\}=
(\wdM_{\chi})^{\sigma}=M^{\sigma}$$
since $\sigma$ stabilizes the elements of $\cZ(\fg_0)$.
Therefore $_{\chi}\wdM^{\sigma}$
is  a (right) $\fg_0$-submodule of a (right) Verma $\fg$-module 
and, at the same time, it is dual to a Verma $\fg_0$-module. 
Since a  Verma $\fg$-module is $\cU(\fn_0^-)$-
torsion-free, $_{\chi}\wdM^{\sigma}$ is  a (right) 
Verma $\fg_0$-module.
Hence $\chi$ is a mate for $\wdchi$ in the ``right sense''.

\subsubsection{}
\label{rightleft}
Let $\chi$ be a perfect mate for $\wdchi$ and
$L$ be a non-zero $\Ug$-bimodule such that $\wdchi.L=L.\wdchi=0$.
Let us show that
$$_{\chi}L_{\chi}:=\{v\in L|\ \chi^r.v=v.\chi^r=0,\ \ r>>0\}$$
is non-zero.

Indeed, according to~\Lem{lemfinsupp}, 
there exist pair-wise distinct $\fg_0$-central characters
$\chi_1,\ldots,\chi_k$
and  positive integers $r_1,\ldots,r_k$ such that 
$\prod_1^k\chi_i^{r_i}.L=0$. Recall that $L_{\chi}\not=0$
since $(\wdchi,\chi)$ is a perfect pair. Thus $\chi=\chi_i$ for a
certain $i$; we can assume that $\chi=\chi_1$. 
One has $\chi_1^{r_1}+\prod_2^k\chi_i^{r_i}=\cZ(\fg_0)$
since $\chi_1,\ldots,\chi_k\in \Max \cZ(\fg_0)$ are pair-wise distinct.
Therefore
$$0\not=L_{\chi}=\prod_2^k\chi_i^{r_i}.L$$
and it is a right submodule of $L$. Clearly
$$_{\chi}L_{\chi}=\{v\in L_{\chi}|\ v.\chi^r=0,\ \ r>>0\}.$$
Using~\ref{leftright}, one obtains $_{\chi}L_{\chi}\not=0$ as required.

\subsection{}
\label{rootsystems}
By~\ref{appr}, to prove the existence of a perfect mate
for a strongly typical $\fg$-central character,
one can choose any triangular decomposition.
Below we describe, for each basic classical Lie superalgebra of type II,
a triangular decomposition we choose for the proof.
We use notation of Kac--- see~\cite{kadv}, 2.5.4 for all
cases except $D(1,2,\alpha)$ where we use $\delta,\epsilon_1,\epsilon_2$
instead $\epsilon_1,\epsilon_2,\epsilon_3$.
The chosen triangular decomposition always corresponds to 
the first ``simplest system of simple roots'' in~\cite{kadv}.

$$\begin{array}{ll}
B(m,n): \ &  \fg_0={\frak {so}}(2m+1)\times{\frak {sp}}(2n)\\
& \pi=\{\delta_1-\delta_2,\ldots,\delta_n-\epsilon_1,
\epsilon_1-\epsilon_2,\ldots,\epsilon_{m-1}-\epsilon_m,\epsilon_m\},\\
& \Delta_0^+=\{\delta_i\pm\delta_j;2\delta_i\}_{1\leq i<j\leq n}\cup
\{\epsilon_i\pm\epsilon_j,;
\epsilon_i\}_{1\leq i<j\leq m},\\
& \Delta_1^+= \{\delta_i;\delta_i\pm\epsilon_j\}
_{1\leq i\leq n,1\leq j\leq m},\\
& \rho_1=(m+{1\over 2})(\delta_1+\delta_2+\ldots+\delta_n)
\end{array}$$
where $\{\delta_1-\delta_2,\ldots,2\delta_n\}$ is a system
of simple roots of ${\frak {sp}}(2n)$ and 
$\{\epsilon_1-\epsilon_2,\ldots,\epsilon_{m-1}-\epsilon_m,\epsilon_m\}$
is a system of simple roots of ${\frak {so}}(2m+1)$.

$$\begin{array}{ll}
D(m,n):\ &  \fg_0={\frak {so}}(2m)\times{\frak {sp}}(2n)\\
& \pi=\{\delta_1-\delta_2,\ldots,\delta_n-\epsilon_1,
\epsilon_1-\epsilon_2,\ldots,\epsilon_{m-1}-\epsilon_m,
\epsilon_{m-1}+\epsilon_m\},\\
& \Delta_0^+=\{\delta_i\pm\delta_j,2\delta_i\}_{1\leq i<j\leq n}\cup
\{\epsilon_i\pm\epsilon_j\}_{1\leq i<j\leq m},\\
& \Delta_1^+= \{\delta_i\pm\epsilon_j\}_{1\leq i\leq n,1\leq j\leq m},\\
& \rho_1=m(\delta_1+\delta_2+\ldots+\delta_n)
\end{array}$$
where $\{\delta_1-\delta_2,\ldots,2\delta_n\}$ is a system
of simple roots of ${\frak {sp}}(2n)$ and 
$\{\epsilon_1-\epsilon_2,\ldots,\epsilon_{m-1}-\epsilon_m,
\epsilon_{m-1}+\epsilon_m\}$
is a system of simple roots of ${\frak {so}}(2m)$.

$$\begin{array}{ll}
D(2,1,\alpha):\ &  
\fg_0={\frak {sl}}(2)\times{\frak {sl}}(2)\times{\frak {sl}}(2),\\
& \pi=\{\delta_1+\epsilon_1+\epsilon_2,-2\epsilon_1,-2\epsilon_2\},\\
& \Delta_0^+=\{2\delta_1;-2\epsilon_1,-2\epsilon_2\},\\
& \Delta_1^+= \{\delta_1\pm\epsilon_1\pm\epsilon_2\},\\
& \rho_1=2\delta_1.\\
 & \\
F(4):\ &  
\fg_0={\frak {sl}}(2)\times {\frak {so}}(7),\\
& \pi=\{{1\over2}(\epsilon_1+\epsilon_2+\epsilon_3+\delta_1),-\epsilon_1,
\epsilon_1-\epsilon_2 ,\epsilon_2-\epsilon_3\},\\
& \Delta_0^+=\{\delta_1\}\cup
\{-\epsilon_i, 1\leq i\leq 3; \pm\epsilon_i-\epsilon_j,1\leq i<j\leq 3\},\\
& \Delta_1^+= \{{1\over 2}
(\delta_1\pm\epsilon_1\pm\epsilon_2\pm\epsilon_3)\},\\
& \rho_1=2\delta_1
\end{array}$$
where $\delta_1$ is a simple root of ${\frak {sl}}(2)$ and
$\{-\epsilon_1,
\epsilon_1-\epsilon_2 ,\epsilon_2-\epsilon_3\}$ is a system
of simple roots of ${\frak {so}}(7)$.

$$\begin{array}{ll}
G(3):\ &  
\fg_0={\frak {sl}}(2)\times G_2,\\
& \pi=\{\delta_1+\epsilon_1,\epsilon_2,
\epsilon_3-\epsilon_2\},\\
& \Delta_0^+=\{2\delta_1\}\cup
\{\epsilon_2, \epsilon_3,-\epsilon_1, \epsilon_3-\epsilon_2,
\epsilon_3-\epsilon_1,\epsilon_2-\epsilon_1\},\\
& \Delta_1^+= \{\delta_1;\delta_1\pm\epsilon_i, i=1,2,3\},\\
& \rho_1={7\over 2}\delta_1
\end{array}$$
where $\epsilon_1+\epsilon_2+\epsilon_3=0$, 
$\{\epsilon_2,\epsilon_3-\epsilon_2\}$ is a system
of simple roots of $G(2)$ and
$2\delta_1$ is a simple root of ${\frak {sl}}(2)$.

The restriction of the non-degenerate invariant
bilinear form $(-,-)$ on $\fg_0$ is a non-degenerate invariant
bilinear form. Thus
$$(\delta_i,\delta_{i'})=(\delta_i,\epsilon_j)=0 \ \ \forall i\not=i'$$
and 
$$(\epsilon_j,\epsilon_{j'})=0  \ \ \forall j\not=j'$$
if $\fg$ is not of the type $G(3)$.

Till the end of this section $\fg$ will denote one (unless otherwise
specified, an arbitrary one) of the basic classical Lie superalgebras
of type II with the fixed triangular decomposition described above.

\subsubsection{More notation}
For all above root systems denote by $n$ the number
of $\delta$'s and by $m$ the number of $\epsilon$'s
(that is $n=1$ for $D(1,2,\alpha), F(4),G(3)$ and $m=2$
for $D(1,2,\alpha)$, $m=3$ for $ F(4),G(3)$).
Remark that for all above systems
$$\rho_1=p\sum_{i=1}^n\delta_i$$ 
for a certain scalar $p\in {1\over 2}{\Bbb N}^+$.

Denote by $W_1$ (resp., $W_2$)
the subgroup of $W$ which acts on $\{\delta_i\}_{i=1}^n$
(resp., on $\{\epsilon_j\}_{j=1}^m$). Then $W=W_1\times W_2$.
Remark that $W_1$ is always the group of signed permutations
of $\{\delta_i\}_{i=1}^n$: for $B(m,n)$ and $D(m,n)$
it is the Weyl group corresponding to ${\frak {sp}}(2n)$
and for $D(1,2,\alpha), F(4),G(3)$
it is the Weyl group corresponding to ${\frak {sl}}(2)$.

For $\mu\in {\fh}^*$ write 
$$\mu=\sum_{i=1}^n k_i\delta_i+\sum_{j=1}^m l_j\epsilon_j$$
and set $\mu_{\delta_i}:=k_i, \mu_{\epsilon_j}:=l_j$.

Say that $\wdchi$ is {\em generic} if there exists
$\lambda\in W({\wdchi})$ such that 
$(\lambda+\rho)_{\delta_i}\not=0$
for $i=1,\ldots,n$. Remark that $T^2\not\in\wdchi$
implies that $\wdchi$ is generic if $\delta_1\in \Delta_1$
that is in the cases $B(m,n)$ and $G(3)$.

Define a lexicographic order on ${\Bbb C}$ by setting $c_1>c_2$
if $\re c_1>\re c_2$ or $\re c_1=\re c_2$ and $\im c_1>\im c_2$.

One can easily sees that for any $\gamma\in\Gamma$ and
$i=1,\ldots,n$ one has $0\leq |\gamma|_{\sigma_i}\leq 2p$  and, 
moreover, $\gamma=\emptyset$
iff $|\gamma|_{\sigma_1}=\ldots=|\gamma|_{\sigma_n}=0$.
We will use many times the following easy lemma.

\subsubsection{}
\begin{lem}{nonsigned}
Let $d;a_1,\ldots,a_s;r_1,\ldots,r_s$ be complex numbers
such that $d,a_1,\ldots,a_s>0$ and
 $0\leq r_1,\ldots,r_s\leq 2d$.
Suppose there exists a signed permutation $\sigma$
which maps $(a_1+d,\ldots,a_s+d)$ to
$(a_1+d-r_1,\ldots,a_s+d-r_s)$.
Then $\sigma$ is a usual (non-signed) permutation
and $r_1=\ldots=r_s=0$.
\end{lem}
\begin{pf}
The permutation $\sigma$ is non-signed because $a_i+a_j>0\geq -2d+r_j$
and so $a_i+d\not=-(a_j+d-r_j)$ for any indexes $i,j$.
Since $\sigma$ is a usual permutation, 
$$\sum_i^s (a_i+d)=\sum_1^s (a_i+d-r_i)$$
that is $\sum r_i=0$. The inequalities $r_1,\ldots,r_s\geq 0$
imply the assertion.
\end{pf}

\subsection{Generic case}
\label{genericcase}
Fix a generic central  character $\wdchi$ .
Since $W_1$ acts by signed permutation
on $\{\delta_i\}_{i=1}^n$,   there exists
$\lambda\in W(\wdchi)$ such that 
$(\lambda+\rho)_{\delta_i}>0$ for all $i=1,\ldots,n$.
Fix such a $\lambda$ and let us show that the 
$\fg_0$-character of $M(\lambda)$ is a perfect mate for $\wdchi$.

To verify that the $\fg_0$-central character of $M(\lambda)$
is a mate for $\wdchi$, recall~\Lem{dist-appr} (i).
Suppose that
$$(\lambda+\rho_0)=w(\lambda+\rho_0-|\gamma|)$$
for some $\gamma\in\Gamma, w\in W$ .
Write 
$$ \lambda+\rho=\sum_1^n k_i \delta_i+\sum_1^m l_j\epsilon_j,\ \ \ 
|\gamma|=\sum_1^n s_i \delta_i+\sum_1^m r_j\epsilon_j.$$
Recall that $k_i>0$ and $0\leq s_i\leq 2p$ for
$i=1,\ldots,n$. One has
$$\lambda+\rho_0=\lambda+\rho+\rho_1=\sum (k_i+p)
\delta_i+\sum l_j\epsilon_j$$
and
$$\lambda+\rho_0-|\gamma|=\sum (k_i+p-s_i)
\delta_i+\sum (l_j-r_j)\epsilon_j.$$
Write $w=w_1w_2$ where $w_1\in W_1,w_2\in W_2$. Then
$$\sum(k_i+p)\delta_i=w_1\bigl(\sum (k_i+p-s_i)\delta_i\bigr)$$
By~\Lem{nonsigned}, this implies $s_1=\ldots=s_n=0$
and thus $\gamma=\emptyset$. 
Hence the $\fg_0$-central character of $M(\lambda)$
is a mate  for $\wdchi$.

One has
$$\Stab_W(\lambda+\rho_0)=
\Stab_W\bigl(\sum (k_i+p)\delta_i+\sum l_j\epsilon_j\bigr)=
\Stab_W\bigl(\sum k_i\delta_i+\sum l_j\epsilon_j\bigr)
=\Stab_W(\lambda+\rho)$$
since $p;k_1,\ldots,k_n>0$.
Using~\Lem{dist-appr} (ii), one concludes that 
the $\fg_0$-central character of 
$M(\lambda)$ is a perfect mate for $\wdchi$.

\subsection{The case when $\wdchi$ is not generic.}
\label{nongencase}
Suppose that $\fg$ is of the type $D(m,n), D(1,2,\alpha)$ or $F(4)$
and $\wdchi$ is not generic.

\subsubsection{Case $D(m,n)$.}
Since $\wdchi$ is not generic, for any $\lambda\in W(\wdchi)$
there exists $i\in\{1,\ldots,n\}$ such that
$(\lambda+\rho)_{\delta_i}=0$. On the other hand, $T^2\not\in\wdchi$
implies $(\lambda+\rho,\beta)\not=0$ for any $\beta\in \Delta_1$.
In particular, $(\lambda+\rho,\delta_i+\epsilon_j)\not=0$ and so
$(\lambda+\rho)_{\epsilon_j}\not=0$ for $j=1,\ldots,m$.
Taking into account that $W_1$ acts on $\{\delta_i\}_1^n$
by signed permutations and $W_2$ acts on $\{\epsilon_j\}_1^m$
by signed permutations changing even number of signs,
one concludes the existence of $\lambda\in W(\wdchi)$ such that
$$\lambda+\rho=\sum_{i=1}^n k_i\delta_i+\sum_{j=1}^m l_j\epsilon_j$$
where
$$\begin{array}{l}
k_1\geq k_2\geq\ldots\geq k_{n-d}>k_{n-d+1}=k_{n-d+2}=\ldots=k_n=0, \\
l_1\geq l_2\geq\ldots\geq l_{m-1}>0,\ 
l_{m-1}\geq l_{m},\ l_{m-1}\geq -l_{m},
\ l_{m}\not=0.
\end{array} 
$$
Fix $\lambda$ as above and set
$$\begin{array}{c}
\gamma_d:=\{\delta_i-\epsilon_j:\ n-d<i\leq n, 1\leq j\leq m-1\}
\cup \{\delta_i-\sn(l_m)\epsilon_m:\ n-d<i\leq n\}
\end{array}$$
where $\sn(l_m):=1$ if $l_m>0$ and $\sn(l_m):=1$ if $l_m<0$.

Let us show that the $\fg_0$-central character of $M(\lambda-|\gamma_d|)$
is a perfect mate for $\wdchi$. To verify that it is a mate
suppose that  
$(\lambda-\gamma_d+\rho_0)=w(\lambda-\gamma+\rho_0)$
for some $\gamma\in \Gamma,w\in W$.
Write 
$$ |\gamma|=\sum_1^n s_i \delta_i+\sum_1^m r_j\epsilon_j.$$
Observe  that $s_i\in \{0,1,\ldots,2m\}$ for
$i=1,\ldots,n$.
One has
\begin{equation}
\label{l-gd}\begin{array}{rl}
\lambda-|\gamma|_d+\rho_0&=\lambda+\rho+(\rho_1-|\gamma_d|)\\
&                         =
\sum_1^{n-d} (k_i+m)\delta_i+\sum_1^{m-1} (l_j+d)\epsilon_j+
(l_m+\sn(l_m)d)\epsilon_m\end{array}
\end{equation}
and
$$\lambda+\rho_0-|\gamma|=\sum_1^{n-d} (k_i+m-s_i)
\delta_i+\sum_{n-d+1}^n (m-s_i)\delta_i+\sum_1^m (l_j-r_j)\epsilon_j.$$
Write $w=w_1w_2$ where $w_1\in W_1,w_2\in W_2$.
Then
$$\sum_1^{n-d}(k_i+m)\delta_i=w_1\bigl(\sum_1^{n-d} (k_i+m-s_i)
\delta_i+\sum_{n-d+1}^n (m-s_i)\delta_i\bigr).$$
For any indexes $i,i'$ such that
$1\leq i\leq n-d< i'\leq n$ 
one has $k_i+m>m\geq \pm (m-s_{i'})$. Therefore $w_1=w'_1w''_1$
where $w'_1$ (resp., $w''_1$) is a signed permutation
of $\{\delta_i\}_1^{n-d}$ (resp., of $\{\delta_i\}^n_{n-d+1}$).
Then $\sum_1^{n-d}(k_i+m)\delta_i=w'_1\bigl(\sum_1^{n-d} (k_i+m-s_i)
\delta_i\bigr)$ and so $s_1=\ldots=s_{n-d}=0$ by~\Lem{nonsigned}.
On the other hand, $0=w''_1\bigl(\sum_{n-d+1}^n (m-s_i)\delta_i\bigr)$
gives  $s_{n-d}=\ldots=s_{n}=m$.

It is easy to see that $s_1=\ldots=s_{n-d}=0$ 
implies $r_j\in \{0,\pm1,\ldots,\pm d\}$
for all $j=1,\ldots,m$. 

Suppose that $l_m>0$. Then
$$w_2\bigl(\sum_1^m (l_j-r_j)\epsilon_j\bigr)=
\sum_1^m (l_j+d)\epsilon_j$$
and, by~\Lem{nonsigned}, $r_j=-d$ for all $j=1,\ldots,m$.
Thus $|\gamma|=|\gamma_d|$. It is easy to check that this implies
the required equality $\gamma=\gamma_d$.

Suppose that $l_m<0$. Then
$$w_2\bigl(\sum_1^m (l_j-r_j)\epsilon_j\bigr)=
\sum_1^{m-1} (l_j+d)\epsilon_j+(l_m-d)\epsilon_m$$
where $w_2$ is a signed permutation of $\{\epsilon_j\}_1^m$. Then
$$w'_2\bigl(\sum_1^{m-1} (l_j-r_j)\epsilon_j+(r_m-l_m)\epsilon_m\bigr)=
\sum_1^{m-1} (l_j+d)\epsilon_j+(-l_m+d)\epsilon_m$$
for another signed permutation $w'_2$ of $\{\epsilon_j\}_1^m$.
Since $l_1,\ldots,l_{m-1},-l_m>0$, \Lem{nonsigned} gives
$r_1=\ldots=r_{m-1}=d$ and $r_m=-d$. 
Thus again $|\gamma|=|\gamma_d|$ and so $\gamma=\gamma_d$.
Hence the $\fg_0$-central character of $M(\lambda-|\gamma|_d)$
is a mate for $\wdchi$.

One has
$$\begin{array}{l}
\Stab_W (\lambda+\rho)=\Stab_W(\bigl(\sum_1^{n-d} k_i\delta_i
+\sum_1^{m} l_j\epsilon_j\bigr)\\
=\Stab_W(\bigl(\sum_1^{n-d} (k_i+m)\delta_i
+\sum_1^{m-1} (l_j+d)\epsilon_j+(l_m+\sn(l_m)d)\epsilon_m \bigr)=
\Stab_W (\lambda-|\gamma_d|+\rho_0)\end{array}$$
since $m;d;k_1,\ldots,k_n;l_1,\ldots,l_{m-1}>\sn(l_m)l_m>0$.
Using~\Lem{dist-appr} (ii), one concludes that 
the $\fg_0$-central character of $M(\lambda-|\gamma|_d)$
 is a perfect mate for $\wdchi$. 

\subsubsection{Case $D(2,1,\alpha)$}
Set $\delta:=\delta_1$.
Since $\wdchi$ is not generic, $(\lambda+\rho)_{\delta}=0$
for any $\lambda\in W(\wdchi)$.
On the other hand, $T^2\not\in\wdchi$
implies $(\lambda+\rho,\beta)\not=0$ for any $\beta\in \Delta_1$.
Since $\Delta_1^+=\{\delta\pm\epsilon_1\pm\epsilon_2\}$, one obtains
$(\lambda+\rho)_{\epsilon_j}\not=0$ for $j=1$ or $j=2$.
Recall that $W\cong {\Bbb Z}_2\times{\Bbb Z}_2\times{\Bbb Z}_2$ where
the copies of ${\Bbb Z}_2$ act by sign on $\delta; \epsilon_1;\epsilon_2$
respectively. 
Therefore there exists $\lambda\in W(\wdchi)$ such that
$$\lambda+\rho=l_1\epsilon_1+l_2\epsilon_2$$
where $l_1,l_2\geq 0$ and at least one of $l_1,l_2$ is non-zero.
Fix such a $\lambda$. Both cases $l_1>0$ and  $l_2>0$ are similar so
we can assume $l_1\not=0$. Set
$$\gamma_d:=\{\delta-\epsilon_1-\epsilon_2;\delta-\epsilon_1+\epsilon_2\}.$$
Let us show that the $\fg_0$-central character of $M(\lambda-|\gamma_d|)$
is a perfect mate for $\wdchi$.

One has $\rho_1-|\gamma_d|=2\epsilon_1$ and
$$\lambda-|\gamma_d|+\rho_0=\lambda+\rho+\rho_1-|\gamma_d|=
(l_1+2)\epsilon_1+l_2\epsilon_2.$$
Assume that $(\lambda-|\gamma_d|+\rho_0)\in W (\lambda-|\gamma|+\rho_0)$
for some $\gamma\in\Gamma$. Then  
$(\lambda-|\gamma|+\rho_0)_{\epsilon_1}=\pm (l_1+2)$ 
that is $|\gamma|_{\epsilon_1}=l_1\pm (l_1+2)$.
Taking into account that $l_1>0$ and that
$|\gamma|_{\epsilon_1}\in \{0,\pm 1,\pm 2\}$, one concludes
$|\gamma|_{\epsilon_1}=-2$. This implies
$\gamma=\gamma_d$. Hence the 
$\fg_0$-central character of $M(\lambda-|\gamma_d|)$
is a mate for $\wdchi$.
Since $l_1>0$ 
$$\Stab_W(\lambda+\rho)=\Stab_W(l_1\epsilon_1+l_2\epsilon_2)
=\Stab_W((l_1+2)\epsilon_1+l_2\epsilon_2)=\Stab_W(\lambda-|\gamma_d|+\rho_0).$$
By~\Lem{dist-appr}, the $\fg_0$-central character of $M(\lambda-|\gamma_d|)$
is a perfect mate for $\wdchi$.

\subsubsection{Case $F(4)$}
Set $\delta:=\delta_1$.
Since $\wdchi$ is not generic, $(\lambda+\rho)_{\delta}=0$
for any $\lambda\in W(\wdchi)$.
On the other hand, $T^2\not\in\wdchi$
implies $(\lambda+\rho,\beta)\not=0$ for any $\beta\in \Delta_1$.
Since $\Delta_1^+=\{{1\over 2}
(\delta\pm\epsilon_1\pm\epsilon_2\pm\epsilon_3\}$, one obtains
$(\lambda+\rho)_{\epsilon_j}\not=0$ for some $j\in\{1,2,3\}$.
Recall that $W=W_1\times W_2$ where $W_1\cong {\Bbb Z}_2$
acts by sign on $\delta$ and $W_2$
is the group of signed permutations of $\{\epsilon_j\}_1^3$.
Therefore there exists $\lambda\in W(\wdchi)$ such that
$$\lambda+\rho=l_1\epsilon_1+l_2\epsilon_2+l_3\epsilon_3$$
where 
$$l_1\geq l_2\geq l_3\geq 0\ \ \&\ \  l_1\not=0.$$
Fix such a $\lambda$. Set
$$\begin{array}{l}
\gamma_d:=\{{1\over 2}(\delta-\epsilon_1\pm\epsilon_2\pm\epsilon_3)\} 
\text{ if } l_1>l_2,\\
\gamma_d:=\{{1\over 2}(\delta-\epsilon_1-\epsilon_2-\epsilon_3);
{1\over 2}(\delta+\epsilon_1-\epsilon_2-\epsilon_3);
{1\over 2}(\delta-\epsilon_1+\epsilon_2-\epsilon_3);
{1\over 2}(\delta-\epsilon_1-\epsilon_2+\epsilon_3)\} 
\text{ if } l_1=l_2. \end{array}$$

We show below that the $\fg_0$-central character of $M(\lambda-|\gamma_d|)$
is a perfect mate for $\wdchi$. 

Indeed, assume that $(\lambda-|\gamma_d|+\rho_0)=w(\lambda-|\gamma|+\rho_0)$
for some $\gamma\in\Gamma, w\in W$. 
Write  $|\gamma|=s\delta+\sum_{1}^3 r_j\epsilon_j$ and
$w=w_1w_2$ where $w_1\in W_1,w_2\in W_2$. One has
$$\lambda-|\gamma|+\rho_0=(2-s)\delta+ \sum_{1}^3 (l_j-r_j)\epsilon_j.$$

Suppose $l_1>l_2$. Then
$$|\gamma_d|=2\delta-2\epsilon_1\ \ \text{ that is }\ \ 
\lambda-|\gamma_d|+\rho_0=
(l_1+2)\epsilon_1+l_2\epsilon_2+l_3\epsilon_3.$$
The equality $(\lambda-|\gamma_d|+\rho_0)=w(\lambda-|\gamma|+\rho_0)$
implies $s=2$. Thus $\gamma$ contains 4 elements and so
$r_1,r_2,r_3\in\{0;\pm 1;\pm 2\}$. Then $l_1+2>l_j-r_j$
for $j=1,2$ and $l_1+2>-(l_j-r_j)$ for $j=1,2,3$. 
Thus $w_2(\epsilon_1)=\epsilon_1$
and so $r_1=-2$. Since $\gamma$ contains 4 elements,
this implies $\gamma=\gamma_d$. Hence 
the $\fg_0$-central character of $M(\lambda-|\gamma_d|)$
is a mate for $\wdchi$.

Suppose $l_1=l_2$. Then
$$|\gamma_d|=2\delta-\epsilon_1-\epsilon_2-\epsilon_3 \ \ \text{ that is }\ \ 
\lambda-|\gamma_d|+\rho_0=\sum_{1}^3 (l_j+1)\epsilon_j.$$

Again $\gamma$ contains 4 elements and so
$r_1,r_2,r_3\in\{0;\pm 1;\pm 2\}$. Assume that $w_2$ is not a usual
permutation that is $w_2(\epsilon_j)=-\epsilon_{j'}$
for some $j,j'$. Then $l_{j'}+1=-(l_j-r_j)$ that is
$r_j-1=l_{j'}+l_j>0$ because $l_1=l_2>0$ and $l_3\geq 0$.
Then $r_j=2$ and, consequently, $r_{j''}=0$ for all $j''\not=j$.
In particular, $w_2$ can change a sign of at most one
of $\epsilon_1,\epsilon_2,\epsilon_3$. Therefore
$l_1+1=l_2+1=l_{j''}-r_{j''}$ for some $j''\not=j$.
Then $l_1+1=l_{j'}$ that contradicts to the inequality $l_1\geq l_3$.
Thus $w_2$ is a usual permutation. 
This implies $\sum_1^3 (l_j+1)=\sum_{1}^3 (l_j-r_j)$
that is $r_1+r_2+r_3=-3$. Taking into account that
$r_j\in \{0;\pm 1;\pm 2\}$ and that the equality
of the form $r_j=\pm 2$ implies $r_{j'}=0$
for $j'\not=j$, one concludes $r_1=r_2=r_3=-1$
that is $|\gamma|=|\gamma_d|$. It is easy to check
that this gives the required equality 
$\gamma=\gamma_d$. Hence 
the $\fg_0$-central character of $M(\lambda-|\gamma_d|)$
is a mate for $\wdchi$.

If $l_1>l_2\leq l_3\leq 0$ then
$$\Stab_W(\lambda+\rho)=\Stab_W(l_1\epsilon_1+l_2\epsilon_2+l_3\epsilon_3)
=\Stab_W ((l_1+2)\epsilon_1+l_2\epsilon_2+l_3\epsilon_3)=
\Stab_W(\lambda-|\gamma_d|+\rho_0).$$
If $l_1=l_2>0$ then $l_3\not=0$ since $(\lambda+\rho,
\delta+\epsilon_1-\epsilon_2+\epsilon_3)\not=0$.
Therefore
$$\begin{array}{rl}
\Stab_W(\lambda+\rho)&=\Stab_W(l_1\epsilon_1+l_2\epsilon_2+l_3\epsilon_3)\\
&=\Stab_W((l_1+1)\epsilon_1+(l_2+1)\epsilon_2+(l_3+1)\epsilon_3)=
\Stab_W(\lambda-|\gamma_d|+\rho_0).\end{array}$$
Hence the $\fg_0$-central character of $M(\lambda-|\gamma_d|)$
is a perfect mate for $\wdchi$.

\section{Annihilation Theorem}
\label{sectannthm}
This section is devoted to the proof of Theorem~\ref{anntypeII}.

\subsection{}
\begin{lem}{soctypV}
Let $\wdM$ be a strongly typical Verma module
and $v\in\wdM$ be a primitive vector.
Then $\Ug v$ is a Verma module. In particular,
$\wdM$ contains a simple Verma submodule.
\end{lem}
\begin{pf}
Set $\wdchi:=\Ann_{\Zg} \wdM$. Since $\wdM$ is strongly typical,
$\wdchi$ is also strongly typical and so for $\wdchi$ there exists
a perfect mate $\chi\in\Max\cZ(\fg_0)$--- see Section~\ref{apprpairs}.
Let $\wdN$ be a submodule of $\wdM$ generated by a primitive
vector. Then $\wdN$ is a quotient of
a Verma module $\wdM'$ and 
$\wdchi\wdM'=0$. By~\ref{dist}, $M':=\wdM'_{\chi}$ is
a Verma $\fg_0$-module. Since $\chi$ is a perfect mate for $\wdchi$,
the $\fg_0$-module $\wdN_{\chi}$ is a non-zero quotient of $M'$.
Taking into account that
$\wdN_{\chi}\subset\wdM$ is torsion-free over $\cU({\frak n}^-_0)$,
one concludes $\wdN_{\chi}=M'$. Thus $(\wdM'/\wdN)_{\chi}=0$
and so $\wdN=\wdM'$ is a Verma module. 

Recall that $\wdM$ has a finite length and so
it contains a simple submodule $\wdN$. A highest weight vector of $\wdN$
is primitive. Hence $\wdN$ is a simple Verma module.
\end{pf}

\subsection{}
\begin{lem}{lemGKuse}
Let $\wdM$ be a strongly typical Verma module and $\wdM'$ 
be a simple submodule of $\wdM$. Then the natural maps
$F(\wdM,\wdM)\to F(\wdM',\wdM)$ and $F(\wdM',\wdM')\iso F(\wdM',\wdM)$
are $\fg$-bimodule isomorphisms.
\end{lem}
\begin{pf}
Denote by $\iota$ the natural map
$F(\wdM,\wdM)\to F(\wdM',\wdM)$ and by $\iota'$
the natural map $F(\wdM',\wdM')\to F(\wdM',\wdM)$ (both
maps are induced by the embedding $\wdM'$ to $\wdM$).
Both maps are $\fg$-bimodule homomorphisms.

To show that $\iota'$ is a bijection and that $\iota$ is a monomorphism,
we  use the following reasoning which is essentially the same
as in~\cite{j16}.

By~\cite{j16}, $F(N_1,N_2)=0$
if  the GK-dimension of any simple quotient
of $N_1$ is greater than GK-dimension of $N_2$
or if the GK-dimension of $N_1$ is less
than GK-dimension of  any simple submodule
of $N_2$. Moreover, the GK-dimension
of a Verma module is equal to the GK-dimension of 
the algebra $\cU({\frak n}^-)$.

Let $v$ be a highest
weight vector of $\wdM$ and $u\in \cU({\frak n}^-)$ be
such that $uv$ is a highest weight vector of $\wdM'$.
Since $\wdM'$ is a Verma module, $u$ is a non-zero divisor
in $\cU({\frak n}^-)$. This implies
$$GK\dim\wdM'=GK\dim\wdM>GK\dim\wdM/\wdM'.$$
The only simple quotient of $\wdM'$ is 
$\wdM'$ itself; thus the inequality $GK\dim\wdM'>GK\dim\wdM/\wdM'$
implies  $F(\wdM',\wdM/\wdM')=0$.
Therefore $\iota'$ is an epimorphism. Obviously, $\iota'$ is a monomorphism.
Hence $\iota'$ is bijective.
By~\Lem{soctypV}, any simple submodule of
$\wdM$ is a Verma module. Thus the GK-dimension
of $\wdM/\wdM'$ is less than  the GK-dimension
of any simple submodule of $\wdM$ and so
$F(\wdM/\wdM',\wdM)=0$. Consequently, $\iota$
is a monomorphism.

It remains to verify the surjectivity of $\iota$.
Since $\iota'$ is bijective, it is enough
to check the surjectivity of the composed map 
$\iota'':=(\iota')^{-1}\circ\iota: F(\wdM,\wdM)\to F(\wdM',\wdM')$.
Denote by $N$ the cokernel of the map $\iota''$.
This is a $\fg$-bimodule and its left and right annihilators
contain $\wdchi:=\Ann_{\Zg}\wdM$. Let $\chi\in\Max\cZ(\fg_0)$
be a perfect mate for $\wdchi$. Using notation of~\ref{rightleft},
one has
$$_{\chi}F(\wdM,\wdM)_{\chi}\iso F(\wdM_{\chi},\wdM_{\chi}),\ \ \
_{\chi}F(\wdM',\wdM')_{\chi}\iso F(\wdM'_{\chi},\wdM'_{\chi}).$$
Since  $\chi$ is a mate for $\wdchi$, both
$M:=\wdM_{\chi}$ and $M':=\wdM'_{\chi}$ are Verma $\fg_0$-modules.
For any $V\in\Irr_0$ one has
$\Hom_{\fg_0}(V,F(M,M))==\dim V|_0=\Hom_{\fg_0}(V ,F(M',M'))$---
see~\cite{j16}, 6.4. Taking into account the injectivity
of $\iota''$, one concludes that $\iota''(F(M,M))=F(M',M')$
and so $_{\chi} N_{\chi}=0$. Thus $N=0$ by~\ref{rightleft}.
Hence $\iota,\iota'$ are isomorphisms.
The lemma is proven.
\end{pf}

Using~\Lem{soctypV} one obtains the

\subsubsection{}
\begin{cor}{GKuse}
Let $\wdM, \wdM'$ be strongly typical Verma modules
and $\wdM'$ be a submodule of $\wdM$. Then
the $\fg$-bimodules $F(\wdM,\wdM)$ and $F(\wdM',\wdM')$
are isomorphic.
\end{cor}

\subsection{}
\begin{prop}{locfinpart}
If $\wdM$ is  a strongly typical projective (in $\wdO$) Verma module then
$$F(\wdM,\wdM)\cong\displaystyle{\oplus_{\wdV\in\Irr}}
\Inj(\wdV)^{\oplus\dim\wdV|_0}.$$
\end{prop}
\begin{pf}
Step 1. Let us show that
for any $\fg_0$-modules $N_1,N_2$ the following $\ad\fg$-modules
$$L:=F(\Ind_{\fg_0}^{\fg} N_1,\Ind_{\fg_0}^{\fg} N_2),\ \ \ 
L':=\Coind_{\fg_0}^{\fg}(F(N_1,N_2)\otimes \Lambda\fg_1)$$
are isomorphic.

Indeed, using the Frobenius reciprocity, one obtains
$$\begin{array}{rl}
\Hom_{\fg_0}(V,F(\Ind_{\fg_0}^{\fg} N_1,\Ind_{\fg_0}^{\fg} N_2))&\cong
\Hom_{\fg_0}(V,F(N_1,N_2)\otimes\Lambda\fg_1\otimes\Lambda\fg_1)\\
& \cong\Hom_{\fg_0}(V,\Coind_{\fg_0}^{\fg}(F(N_1,N_2)\otimes \Lambda\fg_1))
\end{array}$$
for any $V\in\Irr_0$. Hence $L\cong L'$  as $\ad\fg_0$-modules.
Note that $\dim\Hom_{\fg_0}(V,L)<\infty$ for all  $V\in\Irr_0$.

On the other hand, for any $\wdV\in\Irr$ one has
$$\begin{array}{rl}
\Hom_{\fg}(\wdV,F(\Ind_{\fg_0}^{\fg} N_1,\Ind_{\fg_0}^{\fg} N_2))&\cong
\Hom_{\fg}(\Ind_{\fg_0}^{\fg} N_1,\Ind_{\fg_0}^{\fg} N_2\otimes \wdV^*)
\\
&\cong\Hom_{\fg_0}(N_1, \Ind_{\fg_0}^{\fg} N_2\otimes \wdV^*)\cong
\Hom_{\fg_0}(N_1, N_2\otimes\Lambda\fg_1\otimes \wdV^*)\\
&\cong\Hom_{\fg_0}(\wdV, F(N_1,N_2)\otimes \Lambda\fg_1)\\
&\cong\Hom_{\fg}(\wdV,\Coind_{\fg_0}^{\fg}(F(N_1,N_2)\otimes 
\Lambda\fg_1)).
\end{array}$$
Hence $\Soc L\cong\Soc L'$. 
Since $L'$ is injective in $\Fin$, it contains a submodule 
isomorphic to $L$.
Taking into account the $\ad\fg_0$-isomorphism $L\cong L'$, one concludes
that $L\cong L'$  as $\ad\fg$-modules.

Step 2. Let $\chi$ be a perfect mate of $\Ann_{\Zg}\wdM$.
Then $M:=\wdM_{\chi}$ is a Verma $\fg_0$-module and
$\wdM=\Ug M$ by~\Cor{corgenby}. Thus $\wdM$
is a quotient of $\Ind_{\fg_0}^{\fg} M$. Since $\wdM$
is projective, it is a direct summand of $\Ind_{\fg_0}^{\fg} M$
Then the $\ad\fg$-module $F(\wdM,\wdM)$
is a direct summand of the   $\ad\fg$-module
$F(\Ind_{\fg_0}^{\fg} M,\Ind_{\fg_0}^{\fg} M)$.
The last is isomorphic to $\Coind_{\fg_0}^{\fg}(F(M,M)\otimes\Lambda\fg_1)$
and so is injective in $\Fin$. Hence the $\ad\fg$-module $F(\wdM,\wdM)$
 is injective in $\Fin$. 

By~\Lem{soctypV}, $\wdM$ contains 
a simple Verma submodule $\wdM'$. Combining~\Cor{GKuse} and~\ref{dimhom}
one obtains
$$\Soc F(\wdM,\wdM)\cong\Soc F(\wdM',\wdM')
\cong\displaystyle{\oplus_{\wdV\in\Irr}}
\wdV^{\oplus\dim\wdV|_0}.$$
Now the injectivity  of $F(\wdM,\wdM)$ implies
the required assertion.
\end{pf}

\subsection{}
\begin{prop}{propann}
If $\wdM$ is  a strongly typical Verma module then
the natural map $\Ug\to F(\wdM,\wdM)$ is surjective.
\end{prop}
\begin{pf}
Denote by $N$ the cokernel of the natural map $f:\Ug\to  F(\wdM,\wdM)$.
This is a $\fg$-bimodule and its left and right annihilators
contain $\wdchi:=\Ann_{\Zg}\wdM$. Let $\chi\in\Max\cZ(\fg_0)$
be a perfect mate for $\wdchi$. Using notation of~\ref{rightleft}
one has $_{\chi}F(\wdM,\wdM)_{\chi}\iso F(\wdM_{\chi},\wdM_{\chi})$.
Since  $\chi$ is a mate for $\wdchi$,
$M:=\wdM_{\chi}$ is a Verma $\fg_0$-module.
As a $\fg_0$-module, $\wdM$ has a finite length and so
$$\Ann_{\cZ(\fg_0)} \wdM=\chi\prod_{i=1}^k\chi_i^{r_i}$$
where
$\{\chi,\chi_1,\ldots,\chi_k\}=\supp_{\cZ(\fg_0)}\wdM$ and
$r_1,\ldots,r_k\in {\Bbb N}^+$.
Any element of $\prod_{i=1}^k\chi_i^{r_i}$
annihilates $\sum_{i=1}^k \wdM_{\chi_i}$ and acts by scalar
on $\wdM_{\chi}$. Taking into account that
the natural map $\cU\to F(M,M)$ is surjective (see~\cite{j16}, 6.4),
one concludes that
$$f(\cU\prod_{i=1}^k\chi_i^{r_i})=_{\chi}\!\!F(\wdM,\wdM)_{\chi}.$$
and thus $_{\chi} N_{\chi}=0$. By~\ref{rightleft},
$N=0$ as required.
\end{pf}

\subsection{}
\begin{thm}{anntypeII}
For a strongly typical Verma $\fg$-module $\wdM$
$$\Ann\wdM=\Ug\Ann_{\Zg}\wdM.$$
\end{thm}
\begin{pf}
Set $\wdchi:=\Ann_{\Zg}\wdM$. Let $\wdM''$ be a projective
(in $\wdO$) Verma module such that $\wdchi\wdM''=0$
and let $\wdM'$ be a simple submodule of $\wdM''$ (by~\ref{soctypV},
$\wdM'$ is a Verma module). Combining~\ref{GKuse} and~\ref{locfinpart},
one obtains
$$F(\wdM',\wdM')\cong\displaystyle{\oplus_{\wdV\in\Irr}}
\Inj(\wdV)^{\oplus\dim\wdV|_0}.$$
Summarizing~\ref{propJPRV1},\ref{cntrgn} and~\ref{propann},
one concludes that $\Ann\wdM'=\Ug\wdchi$ and that the natural map
$f':\Ug/(\Ug\wdchi)\to  F(\wdM',\wdM')$ is bijective.

Denote by $f$ the natural map $\Ug/(\Ug\wdchi)\to  F(\wdM,\wdM)$
and by $p$ the composition map 
$f\circ (f')^{-1}: F(\wdM',\wdM')\to F(\wdM,\wdM)$.
Obviously $p$ is a $\Ug$-bimodule map. 
By~\ref{propann}, $p$ is surjective. Let us show that
$p$ is bijective. 

Let $\chi$ be a perfect mate for $\wdchi$;
set $M:=\wdM_{\chi}, M':=\wdM'_{\chi}$. Using notation of~\ref{rightleft},
one has
$$_{\chi} F(\wdM,\wdM)_{\chi}=F(M,M),\ \
\ \  _{\chi} F(\wdM',\wdM')_{\chi}=F(M',M')$$
and so $_{\chi} F(\wdM,\wdM)_{\chi},_{\chi} F(\wdM',\wdM')_{\chi}$
are isomorphic $\ad\fg_0$-modules. 
It is easy to see that 
$$_{\chi} (\im p)_{\chi}=p( _{\chi} F(\wdM',\wdM')_{\chi}).$$
Therefore
$p(_{\chi} F(\wdM',\wdM')_{\chi})=_{\chi} F(\wdM,\wdM)_{\chi}$ 
since $p$ is surjective. Taking into
account that the multiplicity
of each simple $\fg_0$-module $V$ in 
$F(M,M)$ is finite, one
concludes that the restriction of $p$ to  $_{\chi} F(\wdM',\wdM')_{\chi}$
is a monomorphism. Thus
$_{\chi} (\ker p)_{\chi}=0$ and $\ker p=0$ by~\ref{rightleft}.
This means that $p$ is bijective
and so $f:\Ug/(\Ug\wdchi)\to  F(\wdM,\wdM)$ is bijective as well.
The assertion follows.
\end{pf}

\subsection{}
\begin{cor}{corannII}
For a strongly typical central character $\wdchi$
$$
\Ug/(\Ug{\wdchi})\cong\displaystyle{\oplus_{\wdV\in\Irr}}
\Inj(\wdV)^{\oplus \dim\wdV|_0}.
$$
\end{cor}

\section{Remark about Verma modules.}
\label{sectRmV}
In this section we study the $\fg_0$-structure of Verma $\fg$-modules.

\subsection{}
Retain notation of Section~\ref{apprpairs}.
Take strongly typical
$\wdchi\in\Max\Zg$ and $\lambda\in W(\wdchi)$.
Recall that as an $\fg_0$-module  $\wdM(\lambda)$ 
has a finite filtration with the factors 
$\{M(\lambda-|\gamma|), \gamma\in\Gamma\}$.
Therefore
\begin{equation}\label{anninZ0}
\Ug\wdchi\cap\cZ(\fg_0)=\Ann\wdM(\lambda)\cap\cZ(\fg_0)=
\prod_i \chi_i^{r_i}
\end{equation}
where the $\chi_i$ are the pairwise
distinct elements of the multiset 
$\{\Ann_{\cZ(\fg_0)}M(\lambda-|\gamma|),\gamma\in\Gamma\}$
and the $r_i$ are positive integers.
Thus for any $\Ug_{\wdchi}$-module $\wdN$ 
$$\wdN=\oplus_i \wdN_{\chi_i},\ \ \ \chi_i^{r_i}\wdN_{\chi_i}=0.$$

Suppose that $\wdM(\lambda)=\oplus_{\gamma\in\Gamma} M(\lambda-\gamma)$
and so $\cZ(\fg_0)$ acts semisimply on $\wdM(\lambda)$.
If, in addition, $\lambda$ is strongly typical then~(\ref{anninZ0})
implies that $\cZ(\fg_0)$ acts semisimply on any $\Ug_{\wdchi}$-module
for $\wdchi=\Ann_{\Zg}\wdM(\lambda)$.

\subsection{}
\begin{prop}{simpleg0}
For $\wdM(\lambda)$ being a simple Verma module
the following conditions are equivalent
$$\begin{array}{ll}
(i) & \wdM=\oplus_{\gamma\in\Gamma} M(\lambda-|\gamma|)\\
(ii) & \forall\gamma,\gamma'\in\Gamma\ \ 
(|\gamma|-|\gamma'|)\in {\Bbb Z}\Delta_0\setminus\{0\}
\Longrightarrow (\lambda-|\gamma|+\rho_0)\not\in 
W(\lambda-|\gamma'|+\rho_0)\\
(iii) & M(\lambda-|\gamma|)\text{ is simple for any } \gamma\in\Gamma.
\end{array}$$
\end{prop}
\begin{pf}
The implication $(ii)\Longrightarrow (i)$ follows from the fact
that any exact sequence $0\to M(\mu)\to N\to M(\mu) \to 0$ in $\cO$
splits (see, for instance,~\ref{projO}).

Let us verify the implication $(i)\Longrightarrow (iii)$.
Recall~\ref{Ndual}. Since $\wdM(\lambda)$ is simple it is isomorphic
to $\wdM(\lambda)^{\#}=\oplus_{\gamma\in\Gamma} 
M(\lambda-|\gamma|)^{\#}$.
The module $\wdM(\lambda)$ is $\cU({\frak n}^-_0)$-torsion-free;
thus all $M(\lambda-|\gamma|)^{\#}$ are
also $\cU({\frak n}^-_0)$-torsion-free. This forces (iii).

The implication $(iii)\Longrightarrow (ii)$ follows
from the fact that if $M(\mu),M(\mu')$ are simple Verma $\fg_0$-modules
with the same $\fg_0$-central character
and $(\mu-\mu')\in {\Bbb Z}\Delta_0$ then $\mu=\mu'$.
This fact can be deduced from~\cite{j}, A.1.14 and A.1.1 (vii).

Hence the conditions (i),(ii) and (iii) are equivalent
provided that $\wdM(\lambda)$ is simple.
\end{pf}
\subsection{}
\begin{cor}{corg0str}
If $\wdM(\lambda)$ contains a simple typical Verma submodule
then the conditions (i) and (ii) are equivalent.
\end{cor}
\begin{pf}
The implication $(ii)\Longrightarrow (i)$ follows from the same argument
as in~\Prop{simpleg0}. To verify the implication $(i)\Longrightarrow (ii)$,
assume that $\wdM(\lambda)=\oplus M_i$ is 
a direct sum of Verma $\fg_0$-modules and that
$$\lambda-|\gamma'|+\rho_0=w(\lambda-|\gamma''|+\rho_0)$$
for some $w\in W$ and $\gamma',\gamma''\in\Gamma$
satisfying $(|\gamma''|-|\gamma'|)\in {\Bbb Z}\Delta_0\setminus\{0\}$.
Take $w'\in W$ such that $\wdM(w'.\lambda)$
is a simple submodule of $\wdM(\lambda)$.
In the notation of~\ref{Gamma} one
has
$$\begin{array}{rl}
w'.\lambda-|w'_*\gamma'|+\rho_0&=w'(\lambda-|\gamma'|+\rho_0)=
w'w(\lambda-|\gamma''|+\rho_0)\\
&=w'w(w')^{-1}(w'.\lambda-|w'_*\gamma''|+\rho_0)\end{array}$$
by~(\ref{w*}).  Furthermore
$|w'_*\gamma'|-|w'_*\gamma''|=w(|\gamma''|-|\gamma'|)\in 
{\Bbb Z}\Delta_0\setminus\{0\}$.
Thus the condition (ii) does not hold for a simple
Verma module $\wdM(w'.\lambda)$. Therefore $\wdM(w'.\lambda)$
has a $\fg_0$-submodule $N$ which is isomorphic to a non-splitting extension
of $M(\mu)$ by $M(\mu')$ for some $\mu,\mu'\in 
\{\lambda-|\gamma|:\gamma\in\Gamma\}$. Since $\wdM(w'.\lambda)$
is a submodule of $\wdM(\lambda)$, $N$ is a $\fg_0$-submodule of
$\wdM(\lambda)$. Since $N$ is indecomposable, it is isomorphic to a 
submodule of a Verma module $M_i$. However  the Gelfand-Kirillov dimension
of any proper quotient
of a Verma $\fg_0$-module is strictly less than
the Gelfand-Kirillov dimension of $M(\mu')$---see~\cite{j16}.
This gives the required contradiction.
\end{pf}
\begin{rem}{}
Recall that a strongly typical Verma module contains 
a simple strongly typical Verma submodule.
\end{rem}

\subsection{}
\label{cond}

Denote by $\Gamma_0$ (resp., $\Gamma_1$)
the set of subsets of $\Delta_1^+$ containing an even (resp., an odd)
number of elements. Take an arbitrary $\lambda\in \fh^*$ and 
fix a ${\Bbb Z}_2$-grading on a Verma module
$\wdM(\lambda)$ in such a way that a highest weight vector becomes
even. As a $\fg_0$-module,  $\wdM=\wdM_0\oplus\wdM_1$.
Let us show that  each
$\wdM(\lambda)_{j}$ has a finite filtration with factors 
$\{M(\lambda-\gamma)| \gamma\in\Gamma_j\}$ ($j=0,1$).

Set
$$\begin{array}{l}
Q_0(\pi):=\{ \sum_{\alpha\in\Delta} k_{\alpha}\alpha|\ 
k_{\alpha}\in{\Bbb Z},
\ \sum_{\alpha\in\Delta_1} k_{\alpha}\ \text{ is even}\},\\
Q_1(\pi):=\{ \sum_{\alpha\in\Delta} k_{\alpha}\alpha|\
 k_{\alpha}\in{\Bbb Z},
\ \sum_{\alpha\in\Delta_1} k_{\alpha}\ \text{ is odd}\}.
\end{array}$$
Note that both $Q_0(\pi),Q_1(\pi)$ are $W$-stable.
We claim that $Q_0(\pi)\cap Q_1(\pi)=\emptyset$ (for 
$\fg={\frak {psl}}(n,n)$  we substitute $\fh$ by
$\hat {\fh}$---see~\ref{pslnn}). Indeed, if
$\fg$ is of the type I then $\fh\cap \cZ(\fg_0)$
contains an element $z$ such that $z(\alpha)=1$ for any
$\alpha\in\Delta_1^+$. Therefore 
$$\sum_{\alpha\in\Delta} k_{\alpha}\alpha(z)=
\sum_{\alpha\in\Delta_1^+} (k_{\alpha}-k_{-\alpha})$$
and so $\mu(z)$ is an even (resp., an odd) integer for
$\mu\in Q_0(\pi)$ (resp., $\mu\in Q_1(\pi)$). 
Let $\fg$ be of the type II. Retain notation of~\ref{rootsystems}. 
One can immediately sees that for $\fg\not= F(4)$ the sum
$\sum_1^n \mu_{\delta_i}$ is an even (resp., an odd) integer
for $\mu\in Q_0(\pi)$ (resp., $\mu\in Q_1(\pi)$).
For the remaining case $\fg=F(4)$, $\mu_{\delta_1}$ is integer 
if $\mu\in Q_0(\pi)$ and belongs to the set ${\Bbb Z}+{1\over 2}$
if $\mu\in Q_1(\pi)$. This implies our claim.

The weights of $\wdM(\lambda)_0$ (resp., $\wdM(\lambda)_1$)
belong to the set $(\lambda-Q_0(\pi))$ (resp., $(\lambda-Q_1(\pi))$).
The weights of $M(\lambda-|\gamma|)$ belong to  $(\lambda-Q_j(\pi))$
iff $\gamma\in\Gamma_j$ ($j=0,1$). This proves that 
$\wdM(\lambda)_{j}$ has a finite filtration with the factors 
$\{M(\lambda-\gamma)| \gamma\in\Gamma_j\}$ ($j=0,1$).

Using arguments as above, it is easy to show
that for $\fg$ of the type II 
$$Q_0(\pi)={\Bbb Z}\Delta_0$$
and for $\fg$  of the type I 
$$Q_0(\pi)\cap \{\mu\in\fh^*|\ \mu(z)=0\} ={\Bbb Z}\Delta_0.$$
This implies that the condition (ii) of~\Prop{simpleg0}
is equivalent to the condition
\begin{equation}
\label{cond*}
\begin{array}{llll}
\forall \gamma',\gamma''\in \Gamma_j & 
(\lambda-|\gamma'|+\rho_0)\in W(\lambda-|\gamma''|+\rho_0) &
\Longrightarrow & |\gamma'|=|\gamma''|
\end{array}\end{equation}
for $j=0,1$.

\subsection{}
\begin{rem}{}
The case $\fg={\frak {osp}}(1,2l)$  was treated 
in~\cite{mu}, 3.7 and~\cite{gl2}, 7.2.
\end{rem}

\section{Appendix}
\label{appendix}
This section contains some lemmas used
in the main text. We also give an alternative proof of
the fact that the $\Zg$-rank of $\Hom_{\fg}(\wdV,\Ug)$
is equal to $\dim \wdV|_0$ for any $\wdV\in\Irr$. 
This alternative proof does not use Separation Theorem~\ref{spr}.

\subsection{}
In this subsection we  prove some lemmas which were used
in the main text.

\subsubsection{}
\label{comb}
Recall that a simple $\fg_0$-module $V(\lambda)$ is finite dimensional
iff $\lambda+\rho_0>w(\lambda+\rho_0)$ for any $w\in W, w\not=\id$.
In principle, the similar fact does not
hold for simple finite dimensional $\fg$-modules.
For instance,
there are triangular decompositions such that the corresponding
$\rho$ is equal to $0$ and so $w.0=0$ for all $w\in W$ even though
$\wdV(0)$ is one-dimensional. 

However, if $\wdV(\lambda)$ is finite dimensional and
strongly typical then $\lambda>w.\lambda$ for all $w\in W, w\not=\id$. 
This can be checked in the following way.
Fix a strongly typical weight $\nu$ such that
$\wdV(\nu)$ is finite dimensional. Write
the character formula~(\ref{tpcch}) in the form
\begin{equation}
\label{chari}
D'\ch\wdV(\nu)=\sum_{w\in W}\sn (w) e^{w.\nu},\ \text{ where } 
D':=(\ch\wdM(0))^{-1}.
\end{equation}
Combining the facts that  $\Stab_W (\nu+\rho)$
is generated by the reflections it contains and 
that $\ch\wdV(\nu)\not=0$, one obtains $w.\nu=w'.\nu$ iff $w=w'$.
A strongly typical Verma module $\wdM(w.\nu)$ has a finite filtration
with the factors of the form $\wdV(w'.\nu)$ where $w'.\nu\leq w.\nu$.
Therefore $\ch\wdM(w.\nu)=\sum_{w'} a_{w,w'}\ch\wdV(w'.\nu)$
where $(a_{w,w'})$ form ``an upper triangular matrix'' that is
$a_{w,w'}=0$ if $w'.\nu\not\leq w.\nu$ and $a_{w,w}=1$.
Consequently, 
$$D'\ch\wdV(w.\nu)=\sum_{w'} b_{w,w'}e^{w'.\nu}
$$
where $b_{w,w'}=0$ if $w'.\nu\not\leq w.\nu$ and $b_{w,w}=1$.
Comparing the last equality with~(\ref{chari}), one obtains
$w'.\nu <\nu$ for all $w'\in W, w'\not=\id$.

\subsubsection{}
\begin{lem}{typtsprod}
Assume that $N_1,N_2$ are finite dimensional ${\frak h}$-diagonalizable
$\fg$-modules and
all simple subquotients of $N_1\otimes N_2$ are typical. Then 
$N_1\otimes N_2$ is a completely reducible module.  
\end{lem} 
\begin{pf}
Since central characters separate 
non-isomorphic typical finite dimensional
modules, it is enough to show that $N_1\otimes N_2$ does not contain
non-trivial extensions $\wdV$ by $\wdV$ for any 
typical finite dimensional
module $\wdV$. Let $N$ be such an extension. 
Then the highest weight subspace
of $N$ is two dimensional and admits a basis $\{v_1,v_2\}$ 
such that $v_1$ is primitive and
$v_1\in {\Ug} v_2$. Then
$v_1\in {\Uh} v_2$ and  so the action of ${\frak h}$ on $N$
is not semisimple. Hence $N_1\otimes N_2$ does not contain submodule
isomorphic to $N$. The assertion follows.
\end{pf}

\subsubsection{}
\begin{lem}{lemfdim}
For any simple finite dimensional module $\wdV$ and any 
$\mu\in {\frak h}^*$, the set of  
$\lambda\in{\frak h}^*$ such that $\wdV(\lambda)$ is finite dimensional
and 
$$\dim\Hom_{\fg}(\wdV,\Hom(\wdV(\lambda-\mu),\wdV(\lambda))
=\dim\wdV|_{\mu}$$
is Zariski dense in ${\frak h}^*$.
\end{lem} 
\begin{pf}
Fix a simple finite dimensional module $\wdV$.
Let $R$ be the subset of $\fh^*$ consisting of the weights
$\lambda$ such that
$\wdV(\lambda)$ is a typical finite dimensional module and
the tensor product $\wdV^*\otimes \wdV(\lambda)$ is the direct sum
of typical simple modules. Let us show that $R$
is Zariski dense in ${\frak h}^*$.

Indeed, take a  finite dimensional $\fg_0$-module
$V(\lambda)$. The induced module
$\Ind_{\frak g_0}^{\frak g} V(\lambda)$ 
has a simple submodule $\wdV(\lambda')$ which is finite dimensional.
The $\fg_0$-module 
$V(\lambda')$ is a $\fg_0$-submodule of $\wdV(\lambda')$
and so is a $\fg_0$-submodule of $\Ind_{\frak g_0}^{\frak g} V(\lambda)$.
As a $\fg_0$-module $\Ind_{\frak g_0}^{\frak g} V(\lambda)\cong
V(\lambda)\otimes\Lambda {\frak g}_1$. 
For finite dimensional ${\frak g}_0$-modules
$L, V(\mu)$, the inequality $\Hom_{\fg_0}(V(\nu),V(\mu)\otimes 
L)\not=0$ implies $(\nu-\mu)\in\Omega(L)$.
Therefore
$(\lambda-\lambda')\in \Omega (\Lambda {\frak g}_1)$. 
Thus for any $\lambda\in\fh^*$ such that $V(\lambda)$ is simple,
there exists $\lambda'\in\lambda+\Omega (\Lambda {\frak g}_1)$ such that 
$\wdV(\lambda')$ is simple. Taking into account that
the set of $\lambda$'s such that  $\dim V(\lambda)<\infty$ 
is Zariski dense in $\fh^*$,
one concludes that the set of $\lambda$'s such that 
 $\dim \wdV(\lambda)<\infty$  is also Zariski dense in $\fh^*$.

The condition on weight to be atypical is polynomial
(see~\ref{typdef}) and so the set of $\lambda$'s 
such that $\wdV(\lambda)$ is typical finite dimensional 
is also Zariski dense.
A module $\wdV(\lambda)$ is a submodule of 
$\Coind_{\frak g_0}^{\frak g} V(\lambda)$ and so 
$\wdV^*\otimes \wdV(\lambda)$
is a submodule of $\wdV^*\otimes 
\Coind_{\frak g_0}^{\frak g} V(\lambda)$.
Again
$\Hom_{\fg_0}(V(\nu),\wdV^*\otimes 
\Coind_{\frak g_0}^{\frak g} V(\lambda))\not=0$ implies
$(\nu-\lambda)\in \Omega(\wdV^*\otimes \Lambda {\frak g}_1)$.
Therefore if $\wdV^*\otimes \wdV(\lambda)$ 
has a subquotient isomorphic
to $\wdV(\nu)$, then 
$(\nu-\lambda)\in \Omega(\wdV^*\otimes \Lambda {\frak g}_1)$.
Using~\Lem{typtsprod} and the fact that
the set $\Omega(\wdV^*\otimes \Lambda {\frak g}_1)$
is finite, we conclude that $R$
is Zariski dense in ${\frak h}^*$.

Frobenius reciprocity gives 
$$\Hom_{{\frak g}}\bigl(\wdV, \Hom(\wdV(\lambda-\mu),\wdV(\lambda))
\bigr)\cong
\Hom_{{\frak g}}\bigl(\wdV(\lambda-\mu),\wdV^*\otimes\wdV(\lambda)\bigr).$$

Take $\lambda\in R$ and denote by $m_{\nu}$ the multiplicity
of $\wdV(\nu)$
in the completely reducible module  $\wdV^*\otimes\wdV(\lambda)$.
The character formula~(\ref{tpcch}) gives
$$D\ch (\wdV^*\otimes \wdV(\lambda))=\sum_{w\in W} \sn (w)
e^{w.\lambda}\ch\wdV^* =\sum_{\nu} m_{\nu}  
\sum_{w\in W} \sn (w) e^{w.\nu}.$$
For typical finite dimensional modules $\wdV(\nu),\wdV(\nu')$ 
the equality $w.\nu=w'.\nu'$ implies $\nu=\nu'$ and $w=w'$
(see~\ref{comb}).
Therefore $\dim\Hom_{{\frak g}}\bigl(\wdV(\lambda-\mu),
\wdV(\lambda)\otimes\wdV^*\bigr)=m_{\lambda-\mu}$ is equal to
the coefficient of the term $e^{\lambda-\mu}$ in the expression
$\sum_{w\in W} \sn (w)e^{w.\lambda}\ch\wdV^*$. 
For ``sufficiently large'' $\lambda\in R$ this coefficient is equal
to $\dim \wdV^*|_{-\mu}$. More precisely, take
$\lambda\in R,w\not=\id$ and $\alpha\in\pi_0$
such that $s_{\alpha}w<w$. Then $w^{-1}\alpha\in\Delta_0^-$
and so $(w(\lambda+\rho),\alpha)=(\lambda+\rho,
w^{-1}\alpha)<0$ by~\ref{comb}. Then
$$(\lambda-\mu-w.\lambda,\alpha)=(\lambda+\rho-\mu,\alpha)
-(w(\lambda+\rho),\alpha)>(\lambda+\rho-\mu,\alpha).$$
As a consequence, for any $\lambda$ belonging to the set
$$R_1:=\{ \lambda\in R|\ (\lambda,\alpha)>(\xi+\mu-\rho,\alpha)\ \ 
\forall\xi\in \Omega(\wdV^*),\alpha\in \pi_0\}$$
the inclusion $(\lambda-\mu)\in (w.\lambda+\Omega(\wdV^*))$
implies $w=\id$. Hence 
$m_{\lambda-\mu}=\dim \wdV^*|_{-\mu}=\dim \wdV|_{\mu}$
for any $\lambda\in R_1$.

For any $\lambda\in R$ and $\alpha\in\pi_0$, 
the value $(\lambda,\alpha)$ belongs to ${\Bbb N}^+$, since 
$\fg_0$-module $V(\lambda)$ is finite dimensional.
Thus $R_1$
is obtained from $R$ by removing the points lying at
finitely many hyperplanes. Taking into account
that $R$ is Zariski dense in ${\frak h}^*$, one concludes that
$R_1$ is also Zariski dense.
This completes the proof.
\end{pf}

\begin{rem}{}
In~\Prop{propfdim} we prove a stronger assertion for the particular
case $\mu=0$.
\end{rem} 

\subsubsection{}
\begin{lem}{}
For all $\lambda,\nu\in {\frak h}^*$ and
all simple finite dimensional $\wdV$ one has
$$\dim\Hom_{\fg}\bigl(\wdV,\Hom(\wdM(\lambda-\mu),
\wdM(\lambda)^{\#})\bigr)=\dim\wdV|_{\mu}$$
\end{lem}
\begin{pf}
Frobenius reciprocity gives 
$$\Hom_{{\frak g}}\bigl(\wdV, \Hom(\wdM(\lambda-\mu),\wdM(\lambda)^{\#})
\bigr)\cong \Hom_{{\frak g}}\bigl(\wdM(\lambda-\mu),
\wdM(\lambda)^{\#}\otimes\wdV^*\bigr).$$
Using notation of~\ref{Vsmod} one has
\begin{equation}
\label{hoho}
\Hom_{{\frak g}}\bigl(\wdM(\lambda-\mu),
\wdM(\lambda)^{\#}\otimes\wdV^*\bigr)\cong
\Hom_{{\frak b}}\bigl({\Bbb C}_{\lambda-\mu},
\wdM(\lambda)^{\#}\otimes\wdV^*\bigr).
\end{equation}

All Verma modules are isomorphic as ${\frak n}^-$-modules.
Therefore their duals are isomorphic as  ${\frak n}^+$-modules.
Furthermore for any $\lambda'\in {\frak h}^*$ the ${\frak b}$-modules
$\wdM(\lambda)^{\#}$ and 
$\wdM(\lambda')^{\#}\otimes {\Bbb C}_{\lambda-\lambda'}$
are isomorphic.
Taking into account~(\ref{hoho}), one
concludes that 
$$k:=\dim\Hom_{{\frak g}}\bigl(\wdV, \Hom(\wdM(\lambda-\mu),
\wdM(\lambda)^{\#})\bigr)$$ 
does not depend
on $\lambda$ (for fixed $\wdV$ and $\mu$).

Recall that $\wdV(\lambda-\mu)$ is a quotient of $\wdM(\lambda-\mu)$
and $\wdV(\lambda)$ is a submodule of 
$\wdM(\lambda)^{\#}$. Consequently 
$k\geq \dim \Hom_{\frak g}\bigl(\wdV(\lambda-\mu),
\wdV(\lambda)\otimes\wdV^*\bigr)$ for all $\lambda\in {\frak h}^*$.
Using~\Lem{lemfdim}, one obtains $k\geq \dim \wdV|_{\mu}$.

To verify that $k\leq \dim \wdV|_{\mu}$ fix
$\lambda$ such that $\wdM(\lambda)^{\#}$ is simple.
Denote by $m_1$ (resp., $m_2$) a highest weight 
vector of $\wdM(\lambda-\mu)$ (resp., $\wdM(\lambda)^{\#}$).
Consider a map 
$$\begin{array}{l}\iota:\Hom_{{\frak g}}\bigl(\wdV, 
\Hom(\wdM(\lambda-\mu),
\wdM(\lambda)^{\#})\bigr)\to \wdV^*|_{-\mu}\ \text{ s.t}\\
\iota(\psi)(v)m_2=\psi(v)(m_1)\ \ \forall v\in \wdV|_{\mu},
\psi\in \Hom_{{\frak g}}\bigl(\wdV, \Hom(\wdM(\lambda-\mu),
\wdM(\lambda)^{\#})\bigr)\end{array}$$
Let us show that $\iota$ is a monomorphism.
Take a non-zero element
$\psi\in \Hom_{{\frak g}}\bigl(\wdV, \Hom(\wdM(\lambda-\mu),
\wdM(\lambda)^{\#})\bigr)$. The vector space 
$\psi(\wdV)(\wdM(\lambda-\mu))$ is 
a non-zero $\fg$-submodule of a simple module $\wdM(\lambda)^{\#}$;
thus it coincides with $\wdM(\lambda)^{\#}$. One has
$$\psi(\wdV)(M_1)=\psi(\wdV)(\cU({\frak n}^-)m_1)=
\cU({\frak n}^-)\psi(\wdV)(m_1)$$
since $\psi(\wdV)$ is $\ad \fg$-stable.
Therefore $\psi(\wdV)(m_1)$ contains the highest weight vector $m_2$
that is $\psi(v)(m_1)=m_2$ for a certain $v\in \wdV$. Obviously
one can choose $v$ be a weight vector; then $v\in \wdV|_{\mu}$
and $\iota(\psi)(v)=1$. Hence $\iota$ is a monomorphism
and so $k\leq \dim\wdV^*|_{-\mu}=\dim \wdV|_{\mu}$ as required.
\end{pf}

\subsubsection{}
\begin{rem}{dimhom}
If $\wdM(\lambda)$ is simple, the above lemma gives
$$\dim\Hom_{\fg}(\wdV,\Hom(\wdM(\lambda),\wdM(\lambda))=\dim\wdV|_0$$
for any simple finite dimensional $\wdV$. 
\end{rem}

\subsubsection{}
\begin{lem}{invfract}
Let $A$ be a polynomial algebra and $W$ be a finite group 
acting on $A$.
Assume that $p,p'\in A$ are such that $p'/p$ is $W$-invariant and 
let $q$ be a maximal $W$-invariant divisor of $p^{|W|}$.
Then there exist  $q'\in A^W$ such that $p'/p=q'/q$.
\end{lem}
\begin{pf}
Any non-zero polynomial has a unique factorization
into irreducible ones. Let  $a/b$ be a reduced form 
of the fraction $p'/p$. For any $s\in W$ one has $s(a)/s(b)=a/b$
and so $b/s(b)$ is a scalar. Since $W$ is finite,
$1=(b/s(b))^{|W|}=b^{|W|}/s(b^{|W|})$. Hence $b^{|W|}\in A^W$. Since $p$
is divisible by $b$,  $p^{|W|}$ is divisible by $b^{|W|}$. 
Therefore $q$ is divisible by $b^{|W|}$. Then
there exist $q'\in A$ such that $q'/q=a/b=p'/p$.
The $W$-invariance of both $q$ and $p'/p$ 
implies the $W$-invariance of $q'$.
\end{pf}

\subsection{}
In this subsection we present alternative proofs of~\Cor{cordetnn}
and~\Thm{fact?}.

{\em Another proof of~\Cor{cordetnn}.}

\subsubsection{}
Fix a subset $Y$ of $\Irr$ and set 
$X:=\{\nu\in\fh^*| \ \wdV(\nu)\in Y\}$.
For each $\nu\in X$ choose a $\Zg$-basic system $\theta_1^{\nu},\ldots,
\theta_{s(\nu)}^{\nu}$ of $\Hom_{\fg}(\wdV(\nu),{\Ug})$.
Let $z_{\nu}\in {\Zg}$ be such an element that 
$\Hom_{\fg}(\wdV(\nu),{\Ug}[z_{\nu}^{-1}])$ is a free 
${\Zg}[z_{\nu}^{-1}]$-module generated by 
$\theta_1^{\nu},\ldots,\theta_{s(\nu)}^{\nu}$
(such an element exists by~\Prop{fractPRV}). 
Let $S\subset \Zg\setminus\{ 0\}$ be a multiplicative
closed set containing $\{T^2;z_{\nu},\nu\in X\}$. Denote by
$A$ the localization of $\Zg$ on $S$ and
by ${\Ug}_A$ the localization of 
$\Ug$ on $S$. Both actions
$\ad\fg$ and $\ad'\fg$ can be canonically extended to ${\Ug}_A$.
Note that as $\ad\fg$-module (resp., $\ad'\fg$-module) ${\Ug}_A$
belongs to $\Fin$. For any $\nu\in X$ the localized algebra
${\Zg}[z_{\nu}^{-1}]$ is a subalgebra of $A$ and so
$\Hom_{\fg}(\wdV(\nu),{\Ug}_A)$ is a free 
$A$-module generated by $\theta_1^{\nu},\ldots,\theta_{s(\nu)}^{\nu}$.
Clearly ${\Ug}_A$ inherits the structure of a superalgebra and its
centre  is equal to $A$. Now we are ready to formulate the

\subsubsection{}
\begin{lem}{lemloc1}
There exist $\ad\fg$-submodules $H, N$ of ${\Ug}_A$ such that
$$\begin{array}{ll}
(i) & \text{ the multiplication map induces a monomorphism }
H\otimes A\to {\Ug}_A,\\
(ii)& {\Ug}_A=HA\oplus N,\\
(iii) & H \cong \oplus_{\nu\in X}\Inj(\wdV(\nu)^{s_{\nu}}),\\
(iv)& \Soc H =\sum_{\nu\in X}\sum_{i=1}^{s(\nu)}\theta_i^{\nu}(\wdV(\nu)),\\
(v)& \Hom_{\fg}(\wdV(\nu),N)=0\ \ \forall\nu\in X.
\end{array}$$
\end{lem}
\begin{pf}
Since  $\theta_1^{\nu},\ldots,
\theta_{s(\nu)}^{\nu}$ is a $\Zg$-basic system
of $\Hom_{\fg}(\wdV(\nu),{\Ug})$, the sum
$$L:=\sum_{\nu\in X}\sum_{i=1}^{s(\nu)}\theta_i^{\nu}(\wdV(\nu))$$
is direct.
By~\ref{TUg}, the $\ad\fg$-module $T{\Ug}$ is injective.
Since $T^2\in\Zg$, the $\ad\fg$-module
$T^{-1}{\Ug}$ is an injective submodule of ${\Ug}_A$.
Therefore $T^{-1}{\Ug}$ contains an injective envelope $H$ of $L$ and 
$H\cong \oplus_{\nu\in X}\Inj(\wdV(\nu)^{s_{\nu}})$. The multiplication 
map induces a ${\frak g}$-homomorphism $\phi:H\otimes A\to {\Ug}_A$. 

Recall that $\Hom_{\fg}(\wdV(\nu),{\Ug}[z_{\nu}^{-1}])$ is a free 
${\Zg}[z_{\nu}^{-1}]$-module generated by 
$\theta_1^{\nu},\ldots,\theta_{s(\nu)}^{\nu}$ for any $\nu\in X$. 
Therefore the module $\Hom_{\fg}(\wdV(\nu),{\Ug}_A)$
is a free $A$-module generated by 
$\theta_1^{\nu},\ldots,\theta_{s(\nu)}^{\nu}$.
This means that for any $\nu\in X$ the restriction 
of $\phi$ on $\sum_{i=1}^{s(\nu)}\theta_i^{\nu}(\wdV)\otimes A$
is a monomorphism and its image coincides with  the isotypical
component of $\wdV(\nu)$ in the socle of ${\Ug}_A$. Then
the restriction of $\phi$ on $L\otimes A$ is a monomorphism and 
$\Soc {\Ug}_A=\phi(L\otimes A)\oplus N'$ where 
$\Hom_{\frak g}(\wdV(\nu),N')=0$ for any $\nu\in X$.
Recall that $\Soc H=L$ and so $\Soc(H\otimes A)=L\otimes A$.
Therefore $\phi$ is a monomorphism
by~\ref{socphi}. From~\ref{suminj}, it follows
that  $H\otimes A$ is an injective module in $\Fin$. Therefore ${\Ug}_A$
contains a submodule $N$ such that  ${\Ug}_A=HA\oplus N$.  Then 
$\Soc {\Ug}_A=HA\oplus \Soc N$ and thus $\Soc N\cong N'$. This
completes the proof.
\end{pf}

\subsection{}
\label{nnzr}
In this subsection we prove that the $\Zg$-rank of the module
$\Hom_{\fg}(\wdV,\Ug)$ is equal to $\dim\wdV|_0$ for
any $\wdV\in\Irr$. Recall that, by~\Rem{frfr}, in order to
prove this assertion it is suffices to find  $\lambda\in {\frak h}^*$
such that the image of $\Soc {\Ug}$ in 
$\End_{\Bbb C}(\wdV(\lambda))$ contains an $\ad{\fg}$-submodule 
$L\cong\wdV^{\oplus \dim \wdV|_0}$.

To find this element $\lambda$, we use
Density Theorem stating the surjectivity of the natural map 
$\Ug\to\End(\wdV)$ for any $\wdV\in\Irr$. As it is shown in~\Lem{lemfdim},
for any  $\wdV\in\Irr$ there exists 
$\lambda$ such that $\End_{\Bbb C}(\wdV(\lambda))$
contains a submodule $\wdV^{\oplus \dim \wdV|_0}$. This would imply the
assertion when ${\Ug}=\Soc {\Ug}$ that is
for completely reducible Lie superalgebras (however for these algebras
the assertion immediately follows from Separation theorems).
The general case requires a certain extra work: one should choose $\lambda$ in
such a way that any copy of $\wdV$ inside the socle
of $\End_{\Bbb C}(\wdV(\lambda))$ lies in the image of $\Soc {\Ug}$.
This can be done with the help of~\Prop{propfdim}.

\subsubsection{}
\begin{prop}{propfdim}
For any $\wdV\in\Irr$ the set of $\lambda\in {\frak h}^*$ 
such that $\dim\wdV(\lambda)<\infty$ 
and there exists a monomorphism
$$\Inj(\wdV)^{\oplus \dim \wdV|_0}\to \End_{\Bbb C}(\wdV(\lambda)).$$
is a Zariski dense subset of $\fh^*$.
\end{prop}
\begin{pf}
Fix $\wdV\in\Irr$. For any 
$\lambda\in {\frak h}^*$ denote by $f_{\lambda}$ the natural
homomorphism ${\Ug}\to \End_{\Bbb C}(\wdV(\lambda))$ and
by $C_{\wdV} (\lambda)$ the isotypical component of $\wdV$
in the socle of $\End_{\Bbb C}(\wdV(\lambda))$. 

Let $X$ be the set of
$\nu\in {\frak h}^*$ such that $\dim\wdV(\nu)\in\Irr$ 
and $\Inj(\wdV(\nu))$ has a subquotient isomorphic to $\wdV$.
We claim that $X$ is finite. Indeed, $\Inj(\wdV(\nu))$ is a submodule
of $\Coind V(\nu)$. As ${\frak g}_0$-module, 
$\Coind V(\nu)\cong V(\nu)\otimes \Lambda {\frak g}_1$ and
so the weight of any ${\frak g}_0$-primitive weight 
vector of $\Coind V(\nu)$  belongs to the set 
$\nu+\Omega(\Lambda {\frak g}_1)$.
Thus the highest weight of $\wdV$ belongs to this set for any $\nu\in X$.
Therefore $X$ is finite.

Retain notation of~\ref{defuu}.
For each $\nu\in X$ choose a $\Zg$-basic 
system $\theta_1^{\nu},\ldots,\theta_{s(\nu)}^{\nu}$ of
$\Hom_{\fg}(\wdV,{\Ug})$. The collection
$\Psi\theta_1^{\nu},\ldots,\Psi\theta_{s(\nu)}^{\nu}$
is a $\Zg$-basic system of $\uu(\wdV)$ and so
it is also a $\Sh$-basic system of 
$\uu(\wdV){\Sh}$ (see~\Lem{lemfractPRV}).
Denote $p^{\nu}$ a minor of this $\Sh$-basic system.
Recall that $X$ is finite and take $z\in {\Zg}$ such that 
$\cP(z)=t^2 \prod_{\nu\in X}\prod_{w\in W} w.p^{\nu}$.
By~\Prop{fractPRV}, for any $\nu\in X$ the localized module 
$\Hom_{\frak g}(\wdV(\nu), {\Ug})[z^{-1}]$
is freely generated over ${\Zg}[z^{-1}]$ by  
$\{\theta_1^{\nu},\ldots,\theta_{s(\nu)}^{\nu}\}$.
Therefore one can apply~\Lem{lemloc1} to the set $X$ and the algebra
$A:={\Zg}[z^{-1}]$. This gives
$\Ug_A=(HA)\oplus N$ where $\Hom_{\fg}(\wdV(\nu),N)=0$
for any $\nu\in X$. Thus $\Soc N=\oplus_{i\in I} \wdV(\mu_i)$
where $\mu_i\not\in X$ for each $i\in I$.
By~\ref{injsoc}, $N$ is isomorphic to a submodule of $\Inj(\Soc N)
\cong\oplus_{i\in I} \Inj (\wdV(\mu_i))$. 
Each module $\Inj (\wdV(\mu_i))$ does not
have a subquotient isomorphic to $\wdV$ because $\mu_i\not\in X$.
Therefore $N$ also
does not have a subquotient isomorphic to $\wdV$. 

From~\Lem{lemfdim} it follows that the set
$$R:=\{\lambda\in {\frak h}^*|\ \dim \wdV(\lambda)<\infty,\ \ 
C_{\wdV}(\lambda)=\wdV^{\oplus\wdV|_0},\ \ \cP(z)(\lambda)\not=0\}$$
is Zariski dense in $\fh^*$. Let us show that for any $\lambda\in R$ there
exists an injective map required in the proposition.

Take $\lambda\in R$. The module $\wdV(\lambda)$ is finite dimensional and, by
Density Theorem, the map $f_{\lambda}$ is surjective. 
Since $z$ acts on $\wdV(\lambda)$ by
a non-zero scalar, one can extend $f_{\lambda}$ to the epimorphism
${\Ug}_A\to \End_{\Bbb C} (\wdV(\lambda))$. Taking into account that $N$
has no subquotient isomorphic to $\wdV$, one concludes
that $C_{\wdV}(\lambda)\subseteq f_{\lambda}(HA)=f_{\lambda}(H)$.  

Let us show that the restriction of $f_{\lambda}$ to $H$
is a monomorphism. By~\ref{socphi}, it is enough to verify
that the restriction of $f_{\lambda}$ to
$$\Soc H=\sum_{\nu\in X}\sum_{i=1}^{s(\nu)} \theta^{\nu}_i(\wdV(\nu))$$
is a monomorphism. Since $\wdV(\nu)$ are pairwise non-isomorphic
for different $\nu$, it suffices to check that
the restriction of $f_{\lambda}$ to
$\sum_{i=1}^{s(\nu)} \theta^{\nu}_i(\wdV(\nu))$ is a monomorphism
for each $\nu\in X$. This follows from~\Cor{cor2ass}.
Indeed, fix $\nu\in X$. Recall that $p_{\nu}$ is a non-zero 
$s(\nu)\times s(\nu)$-minor of
the matrix $\bigl(\cP(\theta_j(v_i))\bigr)_{i=1,r}^{j=1,s(\nu)}$ 
where $v_1,\ldots, v_r$ is basis of $\wdV(\nu)|_0$.
Since $\cP(z)(\lambda)\not=0$
one has $p^{\nu}(\lambda)\not=0$ and so, by~\Cor{cor2ass},
the restriction of $f_{\lambda}$ on 
$\sum_{i=1}^{s(\nu)} \theta^{\nu}_i(\wdV(\nu))$ is a monomorphism.

Combining the facts that $H$ is injective and that
$C_{\wdV}(\lambda)\subseteq f_{\lambda}(H)\cong H$,
one completes the proof.
\end{pf}

\subsubsection{}
\begin{cor}{cornz}
For any simple finite dimensional module $\wdV$ 
the $\Zg$-rank of $\Hom_{\fg}(\wdV,{\Ug})$ is equal to $\dim\wdV|_0$.
\end{cor}
\begin{pf}
Set $r:=\dim\wdV|_0$.
By~\Prop{propfdim} one can choose $\lambda\in {\frak h}^*$ such that
$\wdV(\lambda)$ is finite dimensional and there exists a monomorphism
$$\phi:\  L:=\Inj(\wdV)^{\oplus r}\longrightarrow
\End_{\Bbb C}(\wdV(\lambda)).$$
By Density Theorem the natural map 
$f_{\lambda}:{\Ug}\to\End_{\Bbb C}(\wdV(\lambda))$
is surjective. Since $L$ is projective in $\Fin$, $\Ug$ contains
a submodule $L'$ such that the restriction of $f_{\lambda}$ to $L'$
provides an isomorphism $L'\iso L$. In the light of~\Rem{frfr},
a basis $\theta_1,\ldots,\theta_r$ of $\Hom_{\fg}(\wdV,L)$
is a $\Zg$-basic system of $\Hom_{\fg}(\wdV,{\Ug})$.
\end{pf}

{\em Another proof of~\Thm{fact?}}

\subsubsection{}
\begin{thm}{}
An $\ad\fg$-submodule $L$ of $\Ug$ is a generic harmonic space
iff 
$$\begin{array}{ll}
(a) &  L\cong \displaystyle\oplus_{\wdV\in\Irr}
 \Inj(\wdV)^{\oplus\dim\wdV|_0}\\
\text{and} & \text{one of the following conditions holds}\\
(b) & \forall \wdV\in\Irr\ \   \ \eprv(\wdV;L)\not=0 \\ 
(c) & \text{ the multiplication map provides an embedding }
L\otimes \Zg\to\Ug.
\end{array}$$
\end{thm}
\begin{pf}
Applying~\Lem{lemloc1} for $Y:=\Irr$ and $S:=({\Zg}\setminus\{0\})$,
we conclude the existence of a generic harmonic space
satisfying the condition (a).
 
Let us show that all generic harmonic spaces are
pairwise isomorphic as $\ad\fg$-modules. Indeed, let
$L$ and $L'$  be generic harmonic spaces.
Since ${\Ug}$ is countably dimensional, one can choose
the corresponding sets $S,S'\subset{\Zg}\setminus\{0\}$ having
countable number of elements. Take a maximal ideal
$m$ of $\Zg$ such that $m\cap (S\cup S')=\emptyset$.
Then as $\ad\fg$-modules $L\cong \Ug/(m\Ug)\cong L'$.
Hence all generic harmonic spaces satisfy the condition (a).

Let $L$ be a generic harmonic space. Then the 
condition (c) obviously holds. Moreover, for
any $\wdV\in\Irr$ a basis $\theta_1,\ldots,\theta_r$
of $\Hom_{\fg}(\wdV,L)$ is a $\Zg$-basic system
of $\Hom_{\fg}(\wdV,\Ug)$. Therefore $\eprv(\wdV;L)\not=0$
and so (b) holds as well.

Fix $L$ satisfying (a) and (b). For $\wdV\in\Irr$
denote by $q(\wdV)$ a maximal 
$W.$-invariant divisor of $\bigl(\eprv(\wdV;L)\bigr)^{|W|}$.
Take $S\subset{\Zg}$
such that $\cP(S)$ consists of the elements $t^2q(\wdV)$ where
$\wdV\in\Irr$. 
It is easy to deduce from~\Prop{fractPRV}
that the multiplication map provides an isomorphism
$\Soc L\otimes {\Zg}[S^{-1}]\iso \Soc({\Ug}[S^{-1}])$.
Since $L$ is injective, $L\otimes {\Zg}[S^{-1}]$
is also injective and so 
$L\otimes {\Zg}[S^{-1}]\iso {\Ug}[S^{-1}]$. Hence
$L$ is a generic harmonic space.

Finally, fix $L$ satisfying (a) and (c). For
any $\wdV\in\Irr$ a basis 
of $\Hom_{\fg}(\wdV,L)$ contains $\dim\wdV|_0$ elements
due to the condition (a) and these elements  $\Zg$-linearly
independent due to condition (c). Hence these elements
form a $\Zg$-basic system
of $\Hom_{\fg}(\wdV,\Ug)$ that is $\eprv(\wdV;L)\not=0$.
Thus $L$ fulfills the condition (b) as well.
This completes the proof.
\end{pf}


\end{document}